\documentclass[aap,MSNbibl,seceqn,nameyear,dvips]{arximspdf}
\usepackage{mathbh}
\usepackage{graphicx}

\doi{10.1214/11-AAP831} %kopijuoti is PTS
\volume{22}
\issue{6}
\pubyear{2012}
\firstpage{2560}
\lastpage{2615}

\makeatletter

\newcommand{\widebar}{\overline}

\newcommand{\darkgreybullet}{\mbox{
\includegraphics{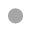}
}}
\newcommand{\darkgreybullett}{\mbox{
\includegraphics{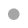}
}}

\newcommand{\xRightarrow}[1]{\stackrel{#1}{\Longrightarrow}}

\newtheorem{proposition}{Proposition}[section]
\newtheorem{corollary}[proposition]{Corollary}
\newtheorem{lemma}[proposition]{Lemma}
\newtheorem{theorem}{Theorem}%[section]

\newproclaim{definition}[proposition]{Definition}
\newproclaim{remark}[proposition]{Remark}
\newproclaim{example}[proposition]{Example}

\DeclareMathAlphabet{\mathpzc}{OT1}{pzc}{m}{it}
\newcommand{\smallu}{\mathpzc{u}}
\newcommand{\smallx}{\mathpzc{x}}
\newcommand{\smally}{\mathpzc{y}}
\newcommand{\ind}[1]{\mathbh{1}_{\{#1\}}} %Definition Indikatorfunktion

\newcommand{\R}{\mathbb{R}}
\newcommand{\M}{\mathbb{M}}
\newcommand{\N}{\mathbb{N}}

\newcommand{\ve}{\varepsilon}

\setattribute{abstract}{width}{285pt}

\makeatother

\begin{document}
\begin{frontmatter}

\title{Tree-valued Fleming--Viot dynamics with mutation and selection}
\runtitle{Tree-valued dynamics with selection}

\begin{aug}
\author[A]{\fnms{Andrej} \snm{Depperschmidt}\thanksref{t1,t2}\ead[label=e3]{depperschmidt@stochastik.uni-freiburg.de}\ead[label=u3,url]{http://www.stochastik.uni-freiburg.de/homepages/deppers/}},
\author[B]{\fnms{Andreas} \snm{Greven}\thanksref{t3}\ead[label=e1]{greven@mi.uni-erlangen.de}\ead[label=u1,url]{http://www.mi.uni-erlangen.de/\textasciitilde greven}}\\
\and
\author[A]{\fnms{Peter} \snm{Pfaffelhuber}\corref{}\thanksref{t1}\ead[label=e2]{p.p@stochastik.uni-freiburg.de}\ead[label=u2,url]{http://www.stochastik.uni-freiburg.de/homepages/pfaffelh/}}

\runauthor{A. Depperschmidt, A. Greven and P. Pfaffelhuber}
\affiliation{University of Freiburg, University of Erlangen and
University of Freiburg}
\address[A]{%
\hspace*{-5pt}A. Depperschmidt\\
\hspace*{-5pt}P. Pfaffelhuber\\
\hspace*{-5pt}Abteilung f{\"u}r mathematische Stochastik\\
\hspace*{-5pt}Albert-Ludwigs University of Freiburg\\
\hspace*{-5pt}Eckerstr. 1\\
\hspace*{-5pt}79104 Freiburg\\
\hspace*{-5pt}Germany\\
\hspace*{-5pt}\printead{e3}\\
\hspace*{-5pt}\hphantom{E-mail: }\printead*{e2}\\
\hspace*{-5pt}URL:\\
\hspace*{-5pt}\printead*{u3}\\
\hspace*{-5pt}\printead*{u2}}
\address[B]{%
\hspace*{-5pt}A. Greven\\
\hspace*{-5pt}Department Mathematik\\
\hspace*{-5pt}University of Erlangen\\
\hspace*{-5pt}Bismarckstr. 1$\tfrac12$\\
\hspace*{-5pt}91054 Erlangen\\
\hspace*{-5pt}Germany\\
\hspace*{-5pt}\printead{e1}\\
\hspace*{-5pt}URL:\\
\hspace*{-5pt}\printead*{u1}} %adresu isvedimo komanda gale!
\end{aug}

\thankstext{t1}{Supported by
the Federal Ministry of Education and Research, Germany
(BMBF) through FRISYS (Kennzeichen 0313921).}

\thankstext{t2}{Supported by the Hausdorff Center in Bonn.}

\thankstext{t3}{Supported by DFG Grant Gr 876/14-1-2.}

% HISTORY:
\received{\smonth{1} \syear{2011}}
\revised{\smonth{11} \syear{2011}}

% ABSTRACT
%
\begin{abstract}
The Fleming--Viot measure-valued diffusion is a Markov process
describing the evolution of (allelic) types
under mutation, selection and random reproduction. We enrich this
process by genealogical relations of individuals so that the random
type distribution as well as the genealogical distances in the
population evolve stochastically. The state space of this
tree-valued enrichment of the Fleming--Viot dynamics with mutation
and selection (TFVMS) consists of marked ultrametric measure
spaces, equipped with the marked Gromov-weak topology and a suitable
notion of polynomials as a separating algebra of test functions.

The construction and study of the TFVMS is based on a well-posed
martingale problem. For existence, we use approximating finite
population models, the tree-valued Moran models, while uniqueness
follows from duality to a function-valued process. Path properties
of the resulting process carry over from the neutral case due to
absolute continuity, given by a new Girsanov-type theorem on marked
metric measure spaces.

To study the long-time behavior of the process, we use a duality based
on ideas from Dawson and Greven [On the effects of migration in spatial
Fleming--Viot models with selection and mutation (2011c) Unpublished
manuscript] and prove ergodicity of the TFVMS if the Fleming--Viot
measure-valued diffusion is ergodic. As a further application, we
consider the case of two allelic types and additive selection. For
small selection strength, we give an expansion of the Laplace transform
of genealogical distances in equilibrium, which is a first step in
showing that distances are shorter in the selective case.
\end{abstract}

% KEYWORDS
%
\begin{keyword}[class=AMS]
\kwd[Primary ]{60K35}
\kwd{60J25}
\kwd[; secondary ]{60J68}
\kwd{92D10}.
\end{keyword}
\begin{keyword}
\kwd{Fleming--Viot process}
\kwd{tree-valued Fleming--Viot dynamics}
\kwd{measure-valued diffusion}
\kwd{metric measure space}
\kwd{resampling}
\kwd{genealogical tree}
\kwd{duality}
\kwd{coalescent}
\kwd{ancestral selection graph}
\kwd{Girsanov theorem}.
\end{keyword}

\end{frontmatter}

%s1 #&#
\section{Introduction}
\label{sintro}
Genealogies are fundamental in studying population models. In this
paper, we focus on the large population limit of constant size
populations evolving under \textit{resampling}, \textit{selection} and
\textit{mutation} in a stochastic fashion. The type distribution of this
limit is modeled by the Fleming--Viot measure-valued diffusion. Here,
resampling is the random reproduction of individuals, mutation is the
random change of (allelic) types of individuals and selection is the
dependence of offspring numbers on the types. By defining random
reproduction we obtain ancestral relations between individuals
described by a randomly evolving genealogy. In our approach, we model
both the genealogical and the type structure in the population.

Populations under selection are modeled either by finitely or by
infinitely many individuals (diffusion). An analysis of the former was
carried out using the biased voter model by \citet{NeuhauserKrone1997}
and \citet{KroneNeuhauser1997}. The large-population limit of the type
frequencies leads to the measure-valued Fleming--Viot dynamics; see,
for example, \citet{FlemingViot1978}, \citet{D93},
\citet{EthierKurtz1993}, Donnelly and Kurtz
(\citeyear{DonnellyKurtz1996,DonnellyKurtz1999}),
Dawson and Greven
(\citeyear{DawsonGreven1999b,DGsel0,DGselGartner2010,DGsel}). A main
tool in the mathematical analysis
of these models is historical information about the population in the
form of genealogical relations of individuals.

%These genealogies are fundamental in applications as
%well.

In applications, genealogies of a population sample are most
important. In particular, mutation rate estimators are based on the
average genealogical distance or the tree length of the genealogical
tree spanned by a sample of individuals
[\citet{Watterson1975}, \citet{Tajima1983}]. Moreover, the
enrichment of
population models by information on ancestral lines has become common
[e.g., \citet{KaplanDardenHudson1988}, \citet
{Kaplanetal1989}]. To cope with
the modeling needs in population genetics, many extensions and
generalizations of the Fleming--Viot dynamics have been given, for example,
the evolution under recombination [see, e.g.,
\citet{D93}, \citet{EthierKurtz1993},
Donnelly and Kurtz (\citeyear{DonnellyKurtz1996,DonnellyKurtz1999})],
as well as
the evolution of a spatially distributed population
[\citet{DGV95},
Dawson and Greven
(\citeyear{DawsonGreven1999b,DGsel0,DGselGartner2010,DGsel})] and general
exchangeable modes of exchange of types [Bertoin and Le~Gall
(\citeyear{BertoinLeGall2003,BertoinLeGall2005,BertoinLeGall2006})].

In order to understand the genealogical structure of population
models, consider the neutral case (i.e., no selection) and a fixed
time $t$ first. Since the resampling mechanism is completely
independent of allelic types, the genealogy can be constructed from
the present to the past using common ancestors of ancestral lines. In
the case of finite variance offspring distributions [and a weak
assumption on their third moments, \citet{MoehleSagitov2001}], the
result is Kingman's coalescent [\citet{Kingman1982a}].

As populations evolve, the underlying genealogies evolve as
well. Consequently, the resampling mechanism allows one to describe
genealogical information of individuals at all times. The main purpose
of the present paper is to give a new approach to studying ancestral
relationships under selection via evolving genealogies. In particular,
we extend the construction of the tree-valued Fleming--Viot dynamics
under neutrality carried out in
\citet{GrevenPfaffelhuberWinter2011}. Note that the resulting processes
are among the first tree-valued stochastic processes in the literature
[but see also Zambotti (\citeyear{MR1814427,MR1891060,MR1959795}),
\citet{EvaPitWin2006}, \citet{EvansWinter2006}, \citet
{mathPR0701657}].

The difficulty in understanding the genealogical structure of a
population under selection already arises for fixed time
genealogies. Most importantly, types and offspring distributions of
individuals are not independent in the selective case. To deal with
this dependence, three different approaches have been used.

First, \citet{KaplanDardenHudson1988}, \citet{Kaplanetal1989}
condition the
construction of the genealogy on the allelic frequency path; see also
\citet{KajKrone2003}, \citet{BartonEtheridgeSturm2004},
\citet{EtheridgePfaffelhuberWakolbinger2006}. If the allelic
frequency path
is known, and an allelic type is present with frequency $x\in[0,1]$ at
time $t$, the rate of coalescence of two lines of this type is
proportional to $1/x$. This construction leads to valuable insights,
for example, into the allelic types of ancestors of the population
[\citet{Taylor2007}].

Second, the ancestral selection graph from \citet{NeuhauserKrone1997}
and \citet{KroneNeuhauser1997} gives a two-step procedure to
derive the
genealogy of a population sample. This construction can, for example,
be used
to see that any ancestor has a higher fitness than a randomly chosen
individual [\citet{Fearnhead2002}]. [Other results derived from the
ancestral selection graph are, e.g., given in \citet{Fearnhead2000},
Slade (\citeyear{Slade2000a}, \citeyear{Slade2000b}) and
\citet{EtheridgeGriffiths2009}.] An important property of this second
approach is that the process generating the genealogy arises as a dual
process of the measure-valued Fleming--Viot process [\citet{Mano2009}].
A connection between the first two approaches has recently been found
in the case of strong balancing selection [\citet{pmid19371754}].

Third, the lookdown construction of \citet{DonnellyKurtz1996} and
\citet{DonnellyKurtz1999} establishes a particle representation of the
Fleming--Viot process with and without selection. Genealogies can as
well be read off from the lookdown process. In the neutral case, the
lookdown construction has, for example, been used to study the
evolution of
the time to the most recent common ancestor of the population
[\citet{PfaffelhuberWakolbinger2006}, \citet
{DelamsEtAl2010}]. In the selective
case, hardly any properties of the genealogies have been read off from
the lookdown process.

In the present paper, we extend the analysis of the neutral tree-valued
Fleming--Viot process from \citet{GrevenPfaffelhuberWinter2011} to include
mutation and selection. This leads to new tree-valued processes
describing the
joint evolution of the allelic type-frequencies and the underlying
genealogy. We encode random genealogies (trees) as random metric
spaces; see
\citet{Evans2000} for the first paper in this direction. In our
construction,
the genealogies evolve forward in time, but contain historical information
about the population. Allelic types are encoded by marks attached to elements
of the metric space.

%f1 #&#
%
\begin{figure}

\includegraphics{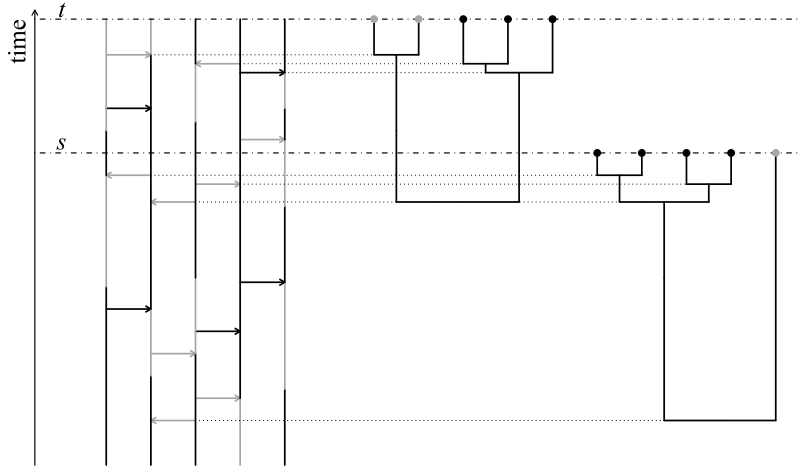}

\caption{Graphical construction of a tree-valued Moran model with
two types with mutation and selection. The fitter type is drawn by
the black line and the weaker type by the gray line. In the left part of
the figure the gray arrows are used independently of color of the involved
lines whereas the black arrows are only used if they start from a black
line. Changes of color along a single line are due to mutations. The right
part shows how the percolation structure on $S_N$ gives rise to a
genealogical tree, that is, a (pseudo-)metric space on the set of
leaves. The
leaves of the tree are marked by the types of the corresponding individuals.}
\label{figgr-costr}
\end{figure}

The starting point of our investigation is the continuous-time Moran
model with mutation and selection. This is a model of a population of
finitely many (distinct) individuals evolving under resampling,
mutation and selection and is best studied by its graphical
representation. At any fixed time, this representation generates a
genealogical tree marked with types; see also
Figure~\ref{figgr-costr}. In a straightforward way, this allows us to
introduce dynamics of genealogies with marks (types) as piecewise
deterministic Markov process with jumps. We show that the large
population limit of this collection of tree-valued Markov processes
exists and is the unique solution of a martingale problem
(Theorems~\ref{T1} and~\ref{T3}). The resulting process is an
enrichment of the measure-valued process and we call it the
\textit{tree-valued Fleming--Viot process with mutation and selection}
(\textit{TFVMS}). On the way, we develop the stochastic analysis for
tree-valued processes. In particular, we give a Girsanov-transform for
our processes and show that genealogies with and without selection can
be studied using a change of measure (Theorem~\ref{T2}).\looseness=1

We continue by showing that the function-valued dual for the
Fleming--Viot process [see, e.g., \citet{D93}] works in the
tree-valued setting. Using this duality and ideas from
Dawson and Greven
(\citeyear{DGsel0,DGselGartner2010,DGsel}), we obtain a stochastic
representation for
the expectation of functionals of sampled finite marked subtrees. As
an application we establish the long-time behavior and the ergodicity
of the TFVMS (Theorem~\ref{T4}), if the measure-valued Fleming--Viot
process is ergodic. We use this equilibrium to study an important
quantity in empirical population genetics in the case of two allelic
types and additive selection: the genealogical distance of two
randomly sampled individuals of the population. We compute the
Laplace transform of the genealogical distance of two sampled
individuals in the case where the selection coefficient is small
(Theorem~\ref{T5}). This result suggests that tree-lengths are shorter
under additive selection. This assertion is widely believed to be true
among biologists, but has never been proved.

Our construction gives a process on the space of marked trees, which
we can treat as marked metric measure spaces. For convenience, we
choose the space of types to be a compact metric space. For the
construction, we require knowledge of fundamental topological
properties of the marked metric measure spaces. While the case without
marks is treated in \citet{GPWmetric09}, topological properties
for the
case with marks are developed in \citet{DGP-topo2011}.

%s2 #&#
\section{Moran models with mutation and selection}
\label{ssbasicdiscrete}
In this section, we first describe a version of the Moran model with
mutation and selection (Section~\ref{ssgraph}), its graphical
construction (Section~\ref{ssgraphconst}) and then extend the
description to the tree-valued case (Section~\ref{sstmmms}). Finally,
we discuss various aspects of models including selection
(Section~\ref{ssbackground}).

%s2.1 #&#
\subsection{The dynamics of the Moran model}
\label{ssgraph}
Fix $N\in\mathbb N$, the population size of the Moran model. Every
individual carries an (allelic) type, element of a set $I$, and we
assume that
%
%e2.1 #&#
%
\begin{equation}\label{addx}
I \mbox{ is a compact metric space}
\end{equation}
for convenience. The individuals of the population are denoted by
$k,l,\ldots\in\{1,\ldots,N\}$. The initial configuration is
$(u_{1}(0),\ldots,u_N(0))$, where $u_k(0)\in I$ denotes the initial type
of individual $k$. The population evolves as a pure jump Markov
process, and the dynamics are given through the following mechanisms.

\begin{longlist}[$\blacktriangleright$]
\item[$\blacktriangleright$] \textit{Resampling} (also known as
\textit{pure genetic drift}):
every (unordered) pair $k\neq l$ is replaced at the
\textit{resampling rate}
%
%e2.2 #&#
%
\begin{equation} \label{eqgamma}
\gamma>0.
\end{equation}
Upon such a resampling event, $l$ is replaced by an offspring of $k$,
or $k$ is replaced by an offspring of $l$, each with probability
$\frac12$. In other words, for every ordered pair $k\neq l$,
individual $l$ is replaced by an offspring of $k$ at rate~$\frac
\gamma2$.
\item[$\blacktriangleright$] \textit{Mutation}: the type of every
individual changes from $u$
to $v$ at rate
%
%e2.3 #&#
%
\begin{equation}\label{eqvarthetabeta}
\vartheta\cdot\beta(u,dv),
\end{equation}
where $\vartheta\geq0$ (the \textit{mutation rate}) and $\beta(\cdot
,\cdot)$
is a stochastic kernel on $I$.
\end{longlist}
For selection, we have two different cases. (See also the discussion
in Section~\ref{ssbackground} on other forms of selection.)
Individuals are either haploid or diploid.

\begin{longlist}[$\blacktriangleright$]
\item[$\blacktriangleright$] \textit{Haploid selection}: every
(ordered) pair $k\neq l$ is
involved in a selection event at rate
%
%e2.4 #&#
%
\begin{equation} \label{eqalphachi}
\frac{\alpha}{N} \cdot\chi(u_k)
\end{equation}
for $\alpha\geq0$ (the \textit{selection coefficient}) and measurable
\textit{fitness function} $\chi\dvtx I\to[0,1]$. Upon a selective event,
individual $l$ is replaced by an offspring of individual $k$.
\item[$\blacktriangleright$] \textit{Diploid selection}: every
(ordered) triple of pairwise
distinct $k,l,m$ is involved in a selection event at rate
%
%e2.5 #&#
%
\begin{equation} \label{eqalphawtchi}
\frac{\alpha}{N^2} \cdot\chi'(u_k,u_m)
\end{equation}
for $\alpha\geq0$ and a symmetric $[0,1]$-valued function $\chi'$
with $\chi'(u,v)=\chi'(v,u)$, which denotes the \textit{fitness} of
the diploid $\{u,v\}$. Again, individual $l$ is replaced by an
offspring of individual $k$.
\end{longlist}
%
%re2.1 #&#
%
\begin{remark}[(Diploid selection)]
While the mechanism for haploid selection is intuitively clear, the
diploid case requires some explanation. Here, $N$ is the number of
haploid individuals, which are arranged in pairs to form
diploids. Since the formation of diploids according to the type
frequencies of the haploids acts on a fast timescale, we can assume
that the population is in Hardy--Weinberg equilibrium at all times,
meaning that the diploid individuals are random pairs of haploids,
and this formation is independent for all times.

Actually, to model diploid selection, we would have to say that every
quadruple $k,l,m,n$ of pairwise distinct individuals is involved in
a selective event at rate $\alpha\cdot\chi'(u_k, u_m)/N^3$ in which
the haploid $l$ from the diploid individual $\{l,n\}$ is replaced by
an offspring of haploid $k$ from the diploid individual
$\{k,m\}$. However, as the haploid individual $n$ is not affected by
such events, our definition above is appropriate.

Haploid and diploid selection leads to the same dynamics in special
cases. In the large population limit, we see that diploid selection
reduces to the haploid case for additive fitness, that is, if $\chi'$
is of the form $\chi'(u,v) = \chi(u) + \chi(v)$ for some function
$\chi$; see (\ref{eq207}) and (\ref{eq209}).
\end{remark}

%s2.2 #&#
\subsection{The graphical construction}
\label{ssgraphconst}
A useful construction of the Moran mod\-el is by means of a random
graph whose main benefit is to automatically generate \textit{ancestral
lines} explicitly. For instance, we use these ancestral lines in
order to bound the number of ancestors of the whole population
(Proposition~\ref{Pupp1}) and show tightness of a sequence of
tree-valued Moran models (see the proof of Theorem~\ref{T3}).
%
%de2.2 #&#
%
\begin{definition}[(Graphical construction of the Moran model)]\label
{defgraph}%
For fixed \mbox{$N\in\mathbb N$}, set
\[
U_N=\{1,\ldots,N\},
\]
and consider the following families of independent Poisson point
processes:
\begin{eqnarray*}% \label{eq303}
\eta_{\mathrm{res}} &:=& \{\eta_{\mathrm{res}}^{k,l}\dvtx k,l \in
U_N\}\qquad
\mbox{each $\eta_{\mathrm{res}}^{k,l}$ with rate } \frac\gamma2, \\
\eta_{\mathrm{mut}} &:=& \{\eta_{\mathrm{mut}}^{k}\dvtx k \in U_N\}
\qquad\mbox{each $\eta_{\mathrm{mut}}^{k}$ with rate } \vartheta
\end{eqnarray*}
and
\begin{eqnarray*}
% \label{eq303b}
\mbox{haploid selection: } \eta_{\mathrm{sel}}
&:=& \{\eta_{\mathrm{sel}}^{k,l}\dvtx k,l \in U_N\}\qquad
\mbox{each $\eta_{\mathrm{sel}}^{k,l}$ with rate } \frac\alpha N,\\
\mbox{diploid selection: } \eta_{\mathrm{sel}}
&:=& \{\eta_{\mathrm{sel}}^{k,l,m}\dvtx k,l,m \in U_N\}\qquad
\mbox{each $\eta_{\mathrm{sel}}^{k,l,m}$ with rate } \frac\alpha{N^2}.
\end{eqnarray*}
The graphical construction of the particle system defines a
percolation structure on the set $S_N:=U_N\times[0,\infty)$. If
$t\in\eta_{\mathrm{res}}^{k,l}$, we draw\vspace*{1pt} an arrow from $(k,t)$ to
$(l,t)$. If $t\in\eta_{\mathrm{sel}}^{k,l}$ in the haploid case, or
$t\in\eta_{\mathrm{sel}}^{k,l,m}$ in the diploid case, draw a
selective arrow from $(k,t)$ to $(l,t)$ in the haploid case and two
different selective arrows from $(k,t)$ to $(l,t)$ and from
$(m,t)$ to $(l,t)$.

Finally, consider the type process $(u_k(t))_{k\in U_N,
t\geq0}$, starting in $u_1(0),\ldots,\break u_N(0)$. Upon a resampling
event $t\in\eta_{\mathrm{res}}^{k,l}$, set $u_l(t)=u_k(t-)$. In
addition, we say that $(k,t-)$ is the ancestor of $(l,t)$ at time
$t-$. For $t\in\eta_{\mathrm{sel}}^{k,l}$, a selective event takes
place with probability $\chi(u_k(t-))$ in the haploid case. In this
case we set $u_l(t)=u_k(t-)$ and say that $(k,t-)$ is the ancestor
of $(l,t)$ at time $t-$. In the diploid case a selective event
$t\in\eta_{\mathrm{sel}}^{k,l,m}$ takes place with probability
$\chi'(u_k(t-), u_m(t-))$, and we set $u_l(t)=u_k(t-)$. In this case
$(k, t-)$ is ancestor of $(l,t)$ at time $t-$. Mutation events take
place at times $t\in\eta_{\mathrm{mut}}^{k}$ where we set $u_k(t)=v$
with probability $\beta(u_k(t-),dv)$.
\end{definition}
%
%ex2.3 #&#
%
\begin{example}[(Example with haploid selection and two types)]%
\label{extwotypes}
The left part of Figure~\ref{figgr-costr} illustrates the
graphical construction of the Moran model in the special case $N=5$,
haploid selection, $I=\{\bullet, {\darkgreybullet}\}$ and
$\chi=\ind{\bullet}$; that is, $\bullet$ is fit and
${\darkgreybullet}$ is unfit. Mutation from
${\darkgreybullet}$ to $\bullet$ and vice versa occurs at
two possibly different rates, denoted
$\vartheta_{\darkgreybullett}$ and
$\vartheta_{\bullet}$. Resampling arrows in
$\eta_{\mathrm{res}}$ are drawn in gray, while selective arrows
in\vadjust{\goodbreak}
$\eta_{\mathrm{sel}}$ are black. Thus, the gray arrows are always
used, whereas the black arrows are only used if they start from
black lines.
\end{example}
%
%re2.4 #&#
%
\begin{remark}[(Convergence to the Fleming--Viot process)]\label{remconvm}
Consider the graphical construction of a Moran
model of size $N$ with mutation and selection from Definition
\ref{defgraph}. For any $t$, the types $u_1(t),\ldots,u_N(t)\in I$
of individuals $1,\ldots,N$ at time $t$ can be read off. We define the
$N$th empirical type distribution process
$\zeta^N=(\zeta^N_t)_{t\geq0}$ by
%
%e2.6 #&#
%
\begin{equation}\label{eq495}
\zeta^N_t:= \frac1N \sum_{k=1}^N \delta_{u_k(t)}.
\end{equation}
It is well known that $\zeta^N \xRightarrow{N\to\infty}
\zeta$, where
$\zeta=(\zeta_t)_{t\geq0}$ is the measure-valued Fleming--Viot
process with mutation and selection; see, for example,
\citet{D93}, \citet{EthierKurtz1993}, \citet{MR1779100}.
In Example~\ref{exmFV}, we recall
its definition via a martingale problem.
\end{remark}

%s2.3 #&#
\subsection{The tree-valued Moran model}
\label{sstmmms}
We are now prepared to define the tree-valued stochastic process
arising from the Moran model with mutation and selection, in terms of
the graphical construction from Definition~\ref{defgraph}. For this
purpose we will need the notion of \textit{ancestors}. From
Figure~\ref{figgr-costr} it is clear that every $l\in U_N$ at time
$t$ has an ancestor $A_s(l,t) \in U_N$ at time $s<t$.
%
%de2.5 #&#
%
\begin{definition}[(Tree-valued Moran model with mutation and
selection)]\label{defTMMMS}
We use the same notation as in Definition
\ref{defgraph}. For every $(l,t)\in S_N$, define
the $U_N$-valued, piecewise constant process $(A_s(l,t))_{0\leq
s\leq t}$ that jumps from $k$ at time $s$ to $j$ at time $s-$, if
$(j,s-)$ is an ancestor of $(k,s)$ at time $s-$. We then say that
$A_s(l,t)$ is the ancestor of $(l,t)$ at time $s$.

The tree-valued Moran model of size $N$ with mutation and selection
takes values in triples $(U_N, r^N, \mu^N)$, where $r^N$ is a
pseudo-metric on $U_N$ [i.e., $r^N(k,l)=0$ is allowed for $k\neq l$]
and $\mu^N$ is a probability measure on $U_N\times I$.

Starting in a pseudo-metric $r_0^N$ on $U_N$, we define for $ k,l\in
U_N$ and $t \ge0$
%
%e2.7 #&#
%
\begin{equation}\label{eqdisttwo}\quad
r^N_t(k,l)
:= \cases{
2\bigl(t - \sup\{s\dvtx A_s(k,t) = A_s(l,t)\}\bigr),\vspace*{2pt}\cr
\hspace*{122.7pt}\qquad \mbox{if $A_0(k,t) = A_0(l,t)$},\vspace*{2pt}\cr
2t + r_0^N(A_0(k,t),A_0(l,t)),\qquad\mbox{else},}\nonumber
\end{equation}
a pseudo-metric $r^N_t$ on $U_N$, such that $r^N_t(k,l)$ is twice
the time to the most recent common ancestor of $k$ and
$l$. Finally, we define the sampling measure as
%
%e2.8 #&#
%
\begin{equation}\label{eqmut}
\mu^N_t:= \frac{1}{N} \sum_{k=1}^N \delta_{(k,u_k(t))}.
\end{equation}
Then the \textit{tree-valued Moran model with mutation and selection} is
given by
%
%e2.9 #&#
%
\begin{equation}\label{agx1}
((U_N, r^N_t, \mu^N_t))_{t\geq0}.\vadjust{\goodbreak}
\end{equation}
\end{definition}
%
%ex2.6 #&#
%
\begin{example}[(Example with two types)]\label{eqtwotypes}
Let us again consider Example~\ref{extwotypes}
and Figure~\ref{figgr-costr}. For any time $t$, a genealogical tree
can be read off for the individuals $(1,t),\ldots,(5,t)$, giving rise
to a (pseudo-)metric on $U_5$ based on genealogical distances. In
addition, the types $u_1(t),\ldots,u_5(t)$ are encoded in the graphical
representation as well and give rise to the empirical measure~$\zeta_t^5$.
\end{example}
%
%re2.7 #&#
%
\begin{remark}[(Trees as marked metric measure spaces, mark
functions)]\label{remtrees1}
(1)~Recall that an ultrametric space can be mapped isometrically
in a unique way onto the set of leaves of a rooted $\R$-tree,
justifying the name \textit{tree-valued}; see also Remark 2.2 in
\citet{GrevenPfaffelhuberWinter2011}.

(2) We call the states $(U_N, r^N_t, \mu^N_t)$ \textit{marked metric
measure spaces} (or \textit{mmm-spaces}); see also
Definition~\ref{defmmm}. To define an appropriate notion of
convergence, we will have to pass from $(U_N, r^N_t,\mu^N_t)$ to
equivalence classes (also defined in detail in Definition
\ref{defmmm}). Roughly speaking, $(U_N, r^N, \mu^N)$ and $(U_N,
{r'}^N,{\mu'}^N)$ are equivalent, if there is a bijection $\sigma$
on $U_N$ with $r^N(\sigma(i), \sigma(j))={r'}^N(i,j)$, and
${\mu'}^N$ is the image of $\mu^N$ under the reordering
$\sigma$. We will write
%
%e2.10 #&#
%
\begin{equation}\label{agx1d}
\mathcal U^N_t = \overline{(U_N, r^N_t, \mu^N_t)}
\end{equation}
for the equivalence class of $(U_N, r^N_t, \mu^N_t)$, and call
$\mathcal U^N=(\mathcal U_t^N)_{t\geq0}$ the \textit{tree-valued
Moran model with mutation and selection} (\textit{TMMMS}).

(3) For the tree-valued Moran model, $((U_N, r_t^N,
\mu^N_t))_{t\geq0}$, we can define a \textit{mark function},
$\kappa_t(k):= u_k(t)$. Moreover, resampling/selection and
mutation occur at different time points, which implies that
$\kappa_t$ is measurable with respect to the
Borel-$\sigma$-algebra of $(U_N, r^N_t)$ for all $t\geq0$, almost
surely. In particular, $\mu^N_t$ has the special form
%
%e2.11 #&#
%
\begin{equation}
\label{eq1023}
\mu^N_t(dx, du) = \Biggl(\frac1N \sum_{k=1}^N \delta_k(dx)\Biggr)
\cdot\delta_{\kappa_t(x)}(du).
\end{equation}
See Remark~\ref{remmark} for more on mark functions in the large
population limit.
\end{remark}

%s2.4 #&#
\subsection{Background on selection}
\label{ssbackground}
Since fitness is the fundamental concept in Darwin's \textit{Origin of
Species}, selection is the most important feature of population
models in biology. A vast amount of literature is devoted to this
topic. We briefly discuss aspects related to the tree-valued
processes.

\subsubsection*{Fertility, viability and state-dependent selection}
In a selective event of the Moran model described in
Section~\ref{ssgraph}, an individual replaces a randomly drawn
individual, independent of the fitness of the replaced individual.
Thus, we take the special form of \textit{fertility selection}
here; that is, individuals might have a \textit{fitness bonus} which
determines their chances to produce a higher number of offspring.
Sometimes, this is also called \textit{positive selection}.\vadjust{\goodbreak}

In the case of \textit{viability} or \textit{negative selection},
individuals have a \textit{fitness malus}, which determines their chances
to die and be replaced by the offspring of a randomly drawn
individual. In the case of viability selection acting on haploids, we
would have a fitness function $\widetilde\chi\dvtx I\to[0,1]$, and every
ordered pair $k\neq l$ is involved in a selective event at rate
$\alpha\cdot\widetilde\chi(u_l)/N$. Upon such an event, individual
$l$ is replaced by an offspring of individual $k$. Our main results,
Theorems~\ref{T1}--\ref{T5}, carry over to the situation of viability
selection.

Also the \textit{state-dependent selection} can be incorporated in our
model. For this, recall the empirical type distribution $\zeta^N$ of
the Moran model of size $N$ from Remark~\ref{remconvm}. Consider the
fitness function $\chi''\dvtx I \times\mathcal M_1(I) \to[0,1]$, that is,
$\chi''(u, \zeta)$ is the fitness of type $u$ if the type distribution
of the total population is $\zeta$. An offspring of individual $k$
replaces the individual $l$ at rate $\frac\alpha N\cdot\chi''(u_k,
\zeta)$. However, if
%
%e2.12 #&#
%
\begin{equation}
\label{eq1113}
\chi''(u,\zeta) = \int\chi'(u,v) \zeta(dv)
\end{equation}
for some $\chi'\dvtx I\times I \to[0,1]$ we find that an offspring of
individual $k$ replaces individual $l$ at selective events occurring
at rate
%
%e2.13 #&#
%
\begin{equation}
\label{eq1111}
\frac\alpha N\cdot\chi''(u_k, \zeta) = \frac\alpha N\cdot\int
\chi'(u_k,v) \zeta(dv) = \frac\alpha{N^2} \sum_{m=1}^N
\chi'(u_k,u_m).
\end{equation}
So, if (\ref{eq1113}) holds, (\ref{eqalphawtchi}) shows that
state-dependent selection is the same as diploid selection. Compare
also Section 7.6 in \citet{MR1779100}.

\subsubsection*{Kin selection}\label{remkinsel}
For measure-valued processes, selection is modeled
by a symmetric function $\widehat\chi'\dvtx I \times I\to\mathbb R$; see
Definition~\ref{defgraph}. In the TMMMS we encode both, the type
distribution and the genealogical tree in the process. This allows us to
treat diploid selection depending also on genealogical distance; that
is,
we can deal with fitness functions of the form
%
%e2.14 #&#
%
\begin{equation}\label{agx1e}
\chi\dvtx I \times I \times\R_+ \longrightarrow[0,1].
\end{equation}
Here, $\chi(u,v,r)$ is the fitness of a diploid individual with
genotype $\{u,v\}$ if the genealogical distance of the two haploids
forming the diploid individual is $r$. Equivalently, if $\smallu=
\overline{(U_N,r^N,\mu^N)}$ is the current state of the TMMMS, then
the offspring of the haploid individual $k\in U_N$ replaces individual
$l\in U_N$ at a selective event taking place at rate
%
%e2.15 #&#
%
\begin{equation}
\label{eq1121}
\frac\alpha N \cdot\sum_{m=1}^N \chi(u_k,u_m,r^N(k,m)).
\end{equation}
A special case of selection depending on genealogical distance is
\textit{kin selection} [e.g., \citet{pmid16593074}], leading to the
concept of \textit{inclusive fitness}
[Hamilton (\citeyear{pmid5875341,pmid5875340})]. The idea is that
the fitness of an
individual is higher if close relatives are around who can help to
raise offspring. Such an altruistic behavior can evolve since it might
also be beneficial for the helpers, because offspring of close
relatives is likely to carry similar genetic material. Such a scenario
can be modeled using a fitness function of the form (\ref{eq1121})
that is decreasing in its third coordinate, that is, in the genealogical
distance.

\subsubsection*{The ancestral selection graph of Krone and Neuhauser}
Genealogies under selection were studied in \citet{NeuhauserKrone1997}
and \citet{KroneNeuhauser1997} by introducing the \textit{ancestral
selection graph} (ASG). The construction can easily be explained
using Figure~\ref{figgr-costr}. Suppose that we are interested in the
genealogy at time $t$. The ASG produces the genealogy in a three-step
procedure from present to the past. Most importantly, when working
backward in time, it is not known in advance if a selective arrow is
used or not.
\begin{longlist}[(2)]
\item[(1)] Going from the top downward through the graphical
representation, consider first the resampling and selective
arrows. Two lines coalesce when a resampling event occurs between
them. If a line hits the tip of a selective arrow, a branching
event occurs. One line, the \textit{continuing line}, is followed in
order to get information on the ancestral line if the selective
arrow is not used, and the other line, the \textit{incoming line},
is followed if the selective arrow is used. Wait until time 0 and
stop the process.
\item[(2)] At time $0$, mark all individuals according to the initial
distribution, and superimpose the mutation process along the
graph, from time $0$ to time~$t$.
\item[(3)] Go through all selective arrows between times $0$ and
$t$. Follow the continuing line if the arrow does not go from a
black line to a gray line, because in this case, the selection
event is not realized. In the other cases, take the incoming
branch.
\end{longlist}
As a result, one obtains genealogical distances of the time $t$
population, together with their types. The main difference between
the ASG and our construction is that the ASG gives the genealogy
only at a single time, while we describe evolving
genealogies. However, our dual process in Section~\ref{sduality} is
reminiscent of the ASG.\vspace*{10pt}

\textit{Outline}:
The paper is organized as follows. In Section~\ref{Sresults}, we
state our main results on the TFVMS process. In Sections~\ref{schar}
and~\ref{sduality} we develop some tools which are not only needed
in the proofs of the main results, but are also of interest in their
own right. The techniques we use are a detailed analysis of the
\textit{generator} of TFVMS (Section~\ref{schar}) and \textit{duality} of
Markov processes (Section~\ref{sduality}). In
Section~\ref{secmart-probl-fixed} we state and prove important facts
concerning the Moran model. For instance, we give the generator
characterization of the finite population model (TMMMS) and discuss
properties of numbers of ancestors and descendants. Finally, the
proofs of our main results are
given in Sections~\ref{SPPP2} and~\ref{SproofAppl}.

We collect the most important notation needed in the paper in
the \hyperref[app]{Appendix}.

%s3 #&#
\section{Results}
\label{Sresults}
In this section we formulate our main results in the set-up of
and under assumptions listed in Sections~\ref{ssgraph} and
\ref{sstmmms}.
%In Remark~\ref{Rex} we discuss possible extensions
%to more general fitness functions $\chi$ and mutation mechanisms.
Our main point is to establish that the weak limit of the process
$((U_N, r^N_t, \mu^N_t))_{t\geq0}$ from
Definition~\ref{defTMMMS} as $N\to\infty$ exists, characterize it
intrinsically and to study its properties. The result is the
generalization of the convergence of the measure-valued Moran models
to the Fleming--Viot diffusion (see Remark~\ref{remconvm}) to the
level of marked genealogical trees.

Before we formulate the results, we have to specify the state space
and give a summary of its properties in Section~\ref{ssstatespace}.
Afterward, in Section~\ref{ssmp}, we give in Theorem~\ref{T1} the
construction of the TFVMS via a well-posed martingale problem. Theorem
\ref{T2} in Section~\ref{ssGirs} gives a Girsanov transformation
between the neutral and the selective tree-valued processes, and
Theorem~\ref{T3} from Section~\ref{ssconv} shows that the TFVMS
arises as weak limit of TMMMS. The long-time behavior of TFVMS is
studied in Theorem~\ref{T4} of Section~\ref{sslongtime}. Finally, an
application to genealogical distances of sampled individuals in
equilibrium is considered in Section~\ref{ssdist}, in Theorem
\ref{T5}.
%
%re3.1 #&#
%
\begin{remark}[(Notation)]\label{remnotation}
For product spaces $X\times Y\times\cdots,$ we denote by $\pi_X,
\pi_Y, \ldots$ the projection operators. For a Polish
space $E$, the function spaces $\mathcal B(E)$
and $\overline{\mathcal C}(E)$ denote the bounded measurable and
bounded continuous, real-valued functions on $E$, respectively. We
denote by $\mathcal M_1(E)$ the space of probability measures on
(the Borel sets of) $E$, equipped with the topology of weak
convergence, abbreviated by $\Rightarrow$. For $\mu\in\mathcal
M_1(E)$ and $\phi\in\mathcal B(E)$, we set $\langle\mu,\phi\rangle:=
\int\phi(x) \mu(dx)$. Moreover, for $\varphi\dvtx E\to E'$ (for some
other Polish space $E'$), the image measure of $\mu$ under $\varphi$
is denoted by $\varphi_\ast\mu$. For $A\subseteq\mathbb R$,
equipped with the Euclidean topology, we denote by $\mathcal C_E(A)$
($\mathcal D_E(A)$) the set of continuous (c\`adl\`ag) functions
$A\to E$, equipped with the topology of uniform convergence on
compact sets (the Skorohod topology).
\end{remark}

%s3.1 #&#
\subsection{State space}
\label{ssstatespace}
Here we introduce the set of isometry classes of marked ultrametric measure
spaces (denoted by $\mathbb U^I$) that will be the state space of both, the
TMMMS and the TFVMS. The starting point of our definition are results from
\citet{GPWmetric09} that are extended in \citet
{DGP-topo2011}. While
$I$ is a
compact metric space in all applications, the notions introduced in this
subsection are valid for any Polish space~$I$.
%
%de3.2 #&#
%
\begin{definition}[(mmm-spaces)]%
\label{defmmm}
(1) An \textit{$I$-marked metric measure space},
\mbox{\textit{$I$-mmm-space}} or \textit{mmm-space}, for short, is a triple
$(X, r, \mu)$ such that $(X,r)$ is a complete and separable metric
space and $\mu\in\mathcal M_1(X\times I)$. Without loss of
generality we assume that $X\subseteq\mathbb R$.

(2) An mmm-space $(X,r,\mu)$ is called \textit{compact} if
$(\operatorname{supp}((\pi_X)_\ast\mu),r)$ is compact. It is called
\textit{ultrametric}
if $(\operatorname{supp}((\pi_X)_\ast\mu),r)$ is ultrametric.

(3) Two mmm-spaces $(X, r_X, \mu_X)$ and $(Y, r_Y, \mu_Y)$ are
\textit{measure-preserving isometric} and $I$-preserving (or
\textit{equivalent}), if there exists a measurable map $\varphi\dvtx X
\to Y$ such that $r_{X}(x,x') = r_{Y}(\varphi(x), \varphi(x'))$
for all $x,x'\in\operatorname{supp}((\pi_X)_\ast\mu_X)$ and
$\widetilde\varphi_\ast\mu_X = \mu_Y$ for $\widetilde\varphi(x,u)
= (\varphi(x),u)$. The equivalence class of an mmm-space
$(X,r,\mu)$ is denoted by $\overline{(X,r,\mu)}$.

(4) We define
%
%e3.1 #&#
%
\begin{equation}
\label{eqPP001}
\mathbb M^I:=\{\overline{(X,r,\mu)}\dvtx (X,r,\mu) \mbox{
mmm-space}\}.
\end{equation}
Moreover,
%
%e3.2 #&#
%
\begin{eqnarray}
\label{eqPP002}
\mathbb M^I_c &:=& \{\overline{(X,r,\mu)}\dvtx (X,r,\mu) \mbox{ compact
mmm-space}\},\nonumber\\
\mathbb U^I &:=& \{\overline{(X,r,\mu)}\dvtx (X,r,\mu) \mbox{ ultrametric
mmm-space}\},\\
\mathbb U^I_c &:=& \mathbb M^I_c \cap\mathbb U^I.
\nonumber
\end{eqnarray}
Generic elements of $\mathbb M^I$ ($\mathbb U^I$) are denoted
by $\smallx, \smally,\ldots(\smallu,\ldots)$.
\end{definition}
%
%re3.3 #&#
%
\begin{remark}[(Pseudo-metrics)] \label{rempseudo}%
Occasionally, we will encounter pseudo-metric
spaces $(X,r)$ [i.e., $r(x_1,x_2)=0$ for $x_1\neq x_2$ is
possible]. The notion of the equivalence class from
Definition~\ref{defmmm} carries over to marked pseudo-metric
measure spaces. Moreover, in the equivalence class
$\overline{(X,r,\mu)}$ of a marked pseudo-metric measure space
$(X,r,\mu)$, we always find an mmm-space $(X',r',\mu')$, such that
the topology on $X$ generated by $r$ is in 1--1 correspondence to the
topology on $X'$ generated by $r'$. That is, the open subsets of $X$
can be mapped onto the open subsets of $X'$ and vice versa. In
particular, it is no restriction to use marked metric measure spaces
instead of marked pseudo-metric measure spaces.
\end{remark}

In order to define an appropriate topology on $\mathbb M^I$, we
introduce the notion of the marked distance matrix distribution.
%
%de3.4 #&#
%
\begin{definition}[(Marked distance matrix distribution)]\label
{defmdistmat} %
Let $(X,r,\mu)$ be an mmm-space,
$\smallx:=\overline{(X,r,\mu)}\in\mathbb M^I$ and
%
%e3.3 #&#
%
\begin{equation} \label{eq22}
R^{(X,r)}\dvtx
\cases{ (X\times I)^{\mathbb N} \to
\mathbb R_+^{{\mathbb N\choose2}}\times I^{\mathbb N}, \cr
((x_i,u_i)_{i\geq1}) \mapsto
((r(x_i, x_j))_{1\leq i<j}, (u_k)_{k\geq1}).}
\end{equation}
The \textit{marked distance matrix distribution of $\smallx=
\overline{(X,r,\mu)}$} is given by
%
%e3.4 #&#
%
\begin{equation}\label{eq23}
\nu^{\smallx}:= \bigl(R^{(X,r)}\bigr)_\ast\mu^{\mathbb N } \in\mathcal
M_1\bigl(\mathbb R^{{\mathbb N\choose2}} \times I^{\mathbb N}\bigr).
\end{equation}
\end{definition}
%
%re3.5 #&#
%
\begin{remark}[(Distance matrix distribution is exchangeable)]%
\label{remdistwell}
(1) Note that $(R^{(X,r)})_\ast\mu^{\mathbb N}$ in the above
definition does not depend on the particular element $(X,r,\mu)$
of $\smallx=\overline{(X,r,\mu)}$. In particular, $\nu^\smallx$ is
well defined. Moreover, by Theorem~1 in \citet{DGP-topo2011}, we
have $\smallx=\smally$ if and only if $\nu^\smallx=\nu^\smally$.

(2) Let
%
%e3.5 #&#
%
\begin{equation}\label{eqSigma}
\Sigma:=\{\sigma\dvtx \mathbb N\to\mathbb N \vert\mbox{$\sigma$ is
injective}\}%\\
\end{equation}
be the set of injective maps on $\mathbb N$. For $\sigma\in\Sigma$, set
%
%e3.6 #&#
%
\begin{equation}
\label{eq29}
R_\sigma\dvtx \cases{ \mathbb R_+^{{\mathbb N\choose2}} \times
I^{\mathbb N} \to\mathbb R_+^{{\mathbb N\choose2}} \times
I^{\mathbb N},\vspace*{2pt}\cr
((r_{ij})_{1 \leq i<j}, (u_k)_{k\geq1})
\mapsto\bigl(\bigl(r_{\sigma(i)\wedge\sigma(j), \sigma(i)\vee\sigma(j)}\bigr),
\bigl(u_{\sigma(k)}\bigr)_{k\geq1}\bigr).}
\end{equation}
Then, for $\smallx\in\mathbb M^I$, the measure $\nu^\smallx$ is
exchangeable in the sense that
%
%e3.7 #&#
%
\begin{equation}
\label{eq28a}
(R_\sigma)_\ast\nu^\smallx= \nu^\smallx.
\end{equation}
\end{remark}
%
%de3.6 #&#
%
\begin{definition}[(Marked Gromov-weak topology)]\label{defmGw}
Let $\smallx, \smallx_1,\smallx_2,\ldots\in\mathbb
M^I$. We say that $\smallx_n \to\smallx$ as $n\to\infty$ in the
\textit{marked Gromov-weak topology} if
%
%e3.8 #&#
%
\begin{equation}\label{de1}
\nu^{\smallx_n} \xRightarrow{n\to\infty} \nu^\smallx
\end{equation}
in the weak topology on $\mathcal M_1(\mathbb R_+^{{\mathbb
N\choose2}} \times I^{\mathbb N})$, where, as usual, ${\mathbb
R}_+^{{\mathbb N\choose2}} \times I^{\mathbb N}$ is equipped with
the product topology of ${\mathbb R}_+$ and $I$, respectively.
\end{definition}

Several topological facts on the marked Gromov-weak topology were
established in \citet{DGP-topo2011}. One of the most important, showing
that $\mathbb M^I$ is a space suitable for probability theory, is
that the space $\mathbb M^I$ is Polish [Theorem 2 in
\citet{DGP-topo2011}]. Before we state our results, we need to
introduce several function spaces on $\mathbb M^I$.
%
%de3.7 #&#
%
\begin{definition}[(Polynomials)]%
\label{defpolyn}
(1) We denote by
%
%e3.9 #&#
%
\begin{eqnarray}\label{BDF1}
\mathcal B_n&:=&\mathcal B_n\bigl(\mathbb R_+^{{\mathbb
N\choose2}}\times I^{\mathbb N}\bigr),\qquad
\overline{\mathcal C}_n:=\overline{\mathcal C}_{n}\bigl(\mathbb
R_+^{{\mathbb
N\choose2}}\times I^{\mathbb N}\bigr),\nonumber\\[-8pt]\\[-8pt]
\overline{\mathcal C}{}^1_n&:=&\overline{\mathcal C}{}^1_{n}\bigl(\mathbb
R_+^{{\mathbb N\choose2}}\times I^{\mathbb N}\bigr),\nonumber
\end{eqnarray}
the sets of bounded measurable (continuous, continuous and
continuously differentiable with respect to all variables in
$\mathbb R_+^{{\mathbb N\choose2}}$) functions $\phi$ on
$\mathbb R_+^{{\mathbb N\choose2}} \times I^{\N}$, such that
$(\underline{\underline r}, \underline
u)\mapsto\phi(\underline{\underline r}, \underline u)$ depends on
the first ${n\choose2}$ variables in $\underline{\underline r}$ and
the first $n$ in $\underline u$\vspace*{1pt} only. (If $n=0$, the spaces
consist of constant functions.)

(2) A function $\Phi\dvtx \mathbb M^I\to\mathbb R$ is a
\textit{polynomial} if, for some $n\in\mathbb N$, there exists
$\phi\in\mathcal B_n$, such that for all $\smallx\in\mathbb M^I$,
%
%e3.10 #&#
%
\begin{equation}\label{eqphi}
\Phi(\smallx):= \Phi^{n,\phi} = \langle\nu^\smallx,
\phi\rangle= \int\phi(\underline{\underline r}, \underline u)
\nu^\smallx(d\underline{\underline r}, d\underline u).
\end{equation}

(3) The \textit{degree} of a polynomial $\Phi$ is the smallest
number $n$
for which there exists $\phi\in\mathcal B_n$ such that (\ref{eqphi})
holds.

(4) Writing $\overline{\mathcal C}{}^{0}_n:=\overline{\mathcal C}_n$,
we set
%
%e3.11 #&#
%
\begin{eqnarray}
\label{eq201}
\Pi&:=& \bigcup_{n=0}^\infty\Pi_n,\qquad \Pi_n
:= \{\Phi^{n,\phi}\dvtx \phi\in\mathcal B_n \},\nonumber\\[-8pt]\\[-8pt]
\Pi^k &:=& \bigcup_{n=0}^\infty\Pi^k_n,\qquad \Pi_n^k:=
\{\Phi^{n,\phi}\dvtx \phi\in\overline{\mathcal C}{}^k_{n}\},\qquad
k=0,1.
\nonumber
\end{eqnarray}
\end{definition}

We use the sets of polynomials as domains for the generator of the
TFVMS process. In this context, we require that $\Pi^1$ is an algebra
that separates points, a result proved in Proposition 4.1 in
\citet{DGP-topo2011}.

%s3.2 #&#
\subsection{Martingale problem}
\label{ssmp}
In this subsection, we define the TFVMS dynamics by a \textit{well-posed
martingale problem}. First we recall the notion of martingale
problems that we use here; see \citet{EthierKurtz86}. Throughout
the following, $I$ is assumed to be a compact metric space (and hence
Polish).

%de3.8 #&#
%
\begin{definition}[(Martingale problem)]\label{D01}
Let $E$ be a Polish space, $\mathbf P_0 \in\mathcal
M_1(E)$, $\mathcal F \subseteq{\mathcal B}(E)$ and $\Omega$ a
linear operator on ${\mathcal B}(E)$ with domain $\mathcal F$. The
law $\mathbf{P}$ of an $E$-valued stochastic process $\mathcal
X=(X_t)_{t\geq0}$ is called a solution of the $(\mathbf
P_0,\Omega,\mathcal F)$-martingale problem if $X_0$ has distribution
$\mathbf P_0$, $\mathcal X$ has paths in the space ${\mathcal
D}_E([0,\infty))$, almost surely, and for all $F\in\mathcal F$,
%
%e3.12 #&#
%
\begin{equation}\label{13def}
\biggl(F(X_t)-\int_0^t \Omega F(X_s) \,ds \biggr)_{t\geq0}
\end{equation}
is a $\mathbf{P}$-martingale with respect to the canonical
filtration. Moreover, the $(\mathbf P_0,\Omega,\mathcal
F)$-martingale problem is said to be well-posed if there is a unique
solution $\mathbf{P}$.
\end{definition}

As an example we now give the martingale problem characterization of
the classical Fleming--Viot diffusion to prepare for the tree-valued
process.

%ex3.9 #&#
%
\begin{example}[(The measure-valued Fleming--Viot process)] \label
{exmFV}%
We recall the classical Fleming--Viot
measure-valued diffusion $\zeta=(\zeta_t)_{t\geq0}$ with mutation
and selection. It arises as the large population limit of the
process describing the evolution of type frequencies $\zeta^N =
(\zeta^N_t)_{t\geq0}$ in the Moran models introduced in
Section~\ref{ssbasicdiscrete}. The state space is $\mathcal
M_1(I)$, and $\zeta_t$ describes the distribution of allelic types
in the population at time $t$.

The process can be characterized in various ways by a martingale
problem, for example, by a second order differential operator on
$\widebar{\mathcal C}(\mathcal M_1(I))$ with domain
$\widebar{\mathcal C}{}^2 (\mathcal M_1(I))$, with an appropriate
definition of the derivative. However, our choice of an operator on
polynomials reveals best the connection to the tree-valued process.

Define\vspace*{1pt} the set of polynomials $\mathcal F$ on $\mathcal
M_1(I)$ by letting $\mathcal F = \bigcup_{n=0}^\infty\mathcal F_n$,
where $\mathcal F_n$ is the set of functions $\widehat\Phi\dvtx
\mathcal M_1(I)\to\mathbb R$ with $\widehat\Phi(\zeta) = \langle
\zeta^{\otimes\mathbb N}, \widehat\phi\rangle$ for some\vspace*{1pt}
$\widehat\phi\in\overline{\mathcal C}(I^{\mathbb N})$, depending only
on the first $n$ variables. Define the linear operator on
$\overline{\mathcal C}(\mathcal M_1(I))$ with domain $\mathcal F$
%
%e3.13 #&#
%
\begin{equation}
\label{eq202}
\widehat\Omega= \widehat\Omega^{\mathrm{res}} +
\widehat\Omega^{\mathrm{mut}} + \widehat\Omega^{\mathrm{sel}}.
\end{equation}
Here, for $\widehat\Phi\in\mathcal F_n$ with $\widehat\Phi(\zeta)
= \langle
\zeta^{\otimes\mathbb N},\widehat\phi\rangle$, the different terms
are given
as follows:
\begin{longlist}[(2)]
\item[(1)] For \textit{resampling rate} $\gamma>0$, the \textit{resampling
operator} is defined by
%
%e3.14 #&#
%
\begin{equation}
\label{eq203}
\widehat\Omega^{\mathrm{res}} \widehat\Phi(\zeta) = \frac\gamma2
\sum_{k,l=1}^n\langle\zeta^{\otimes\mathbb N}, \widehat\phi\circ
\widehat\theta_{k,l} - \widehat\phi\rangle,
\end{equation}
where the \textit{replacement operator} $\widehat\theta_{k,l}$ is the map
which replaces the $l$th component of an infinite sequence by the $k$th;
that is, for $\underline u = (u_1,u_2,\ldots)$,
%
%e3.15 #&#
%
\begin{eqnarray}\label{gen2}
\widehat\theta_{k,l}(\underline u) &:=& \underline
u_l^{u_k},\nonumber\\[-8pt]\\[-8pt]
\underline u_l^v &:=& (u_1,\ldots,u_{l-1},v,u_{l+1},\ldots).
\nonumber
\end{eqnarray}
\item[(2)] For \textit{mutation rate} $\vartheta\geq0$, the \textit{mutation
operator} is defined by
%
%e3.16 #&#
%
\begin{equation}
\label{eq203b}
\widehat\Omega^{\mathrm{mut}} \widehat\Phi(\zeta) = \vartheta
\sum_{k\geq1} \langle\zeta^{\otimes\mathbb N}, \widehat
B_k\widehat\phi\rangle,
\end{equation}
where, for some stochastic kernel $\beta(\cdot,\cdot)$ on $I$,
%
%e3.17 #&#
%
\begin{eqnarray}
\label{eq205}
\widehat B_k\widehat\phi &:=& \widehat\beta_k
\widehat\phi-\widehat\phi,\nonumber\\[-8pt]\\[-8pt]
(\widehat\beta_k \widehat\phi)(\underline u) &:=& \int
\widehat\phi(\underline u_k^v)\beta(u_k, dv).
\nonumber
\end{eqnarray}
That is, $\widehat B$ is the bounded generator of a Markov jump
process on $I$ with c\`adl\`ag paths. It is always possible to
write
%
%e3.18 #&#
%
\begin{equation}
\label{eq1001}
\beta(u, dv) = z \widebar\beta(dv) + (1-z)
\widetilde\beta(u,dv)
\end{equation}
for some $z\in[0,1]$, $\widebar\beta\in\mathcal M_1(I)$ and a
stochastic kernel $\widetilde\beta(\cdot,\cdot)$ on $I$. We refer
to the
case $z=1$ as \textit{parent-independent mutation} or the
\textit{house-of-cards} model. The latter was introduced in
\citet{Kingman1978} who argued that mutations might destroy the
fragile fitness advantage, which was built up during evolution, and
lead to a replacement with an independent type. In this case,
%
%e3.19 #&#
%
\begin{equation} \label{eq206} \beta(u,dv) = \widebar\beta(dv)\qquad
\mbox{does not depend on }u\in I.
\end{equation}
For $z\in(0,1]$, we say that mutation has a parent-independent
component.\vspace*{8pt}
\item[(1)] For \textit{selection intensity} $\alpha\geq0$, the
\textit{selection operator} is given by
%
%e3.20 #&#
%
\begin{equation}\label{eq207}
\widehat\Omega^{\mathrm{sel}} \widehat\Phi(\zeta) =
\alpha\sum_{k\geq1} \langle\zeta^{\otimes\mathbb N}, \widehat
\phi\cdot{\widehat\chi}'_{k,n+1} - \widehat\phi\cdot
{\widehat\chi}'_{n+1,n+2}\rangle.
\end{equation}
Here, the \textit{fitness function}
%
%e3.21 #&#
%
\begin{equation}
\label{eq206b}
\widehat\chi'\dvtx I \times I \to[0,1]
\end{equation}
is measurable and symmetric in both coordinates, and
${\widehat{\chi}'}_{k,l}$ acts on the $k$th and $l$th
coordinate. The special case for $\chi^\prime$, when there exists a
function
%
%e3.22 #&#
%
\begin{equation}\label{eq208}
\widehat\chi\dvtx I\to[0,1] \qquad\mbox{with } \widehat\chi'(u,v) =
\widehat\chi(u) + \widehat\chi(v)
\end{equation}
is called \textit{additive selection} or \textit{haploid
selection}. In this case,
%
%e3.23 #&#
%
\begin{equation}\label{eq209}
\widehat\Omega^{\mathrm{sel}} \widehat\Phi(\zeta) =
\alpha\sum_{k\geq1} \langle\zeta^{\otimes\mathbb N}, \widehat
\phi\cdot\widehat\chi_{k} - \widehat\phi\cdot\widehat
\chi_{n+1}\rangle,
\end{equation}
where $\widehat\chi_k$ acts on the $k$th coordinate. Note that
selective events lead to replacements of individuals similar to
resampling events [see also (\ref{eq11}) and (\ref{eq12}) in the
case of Moran models]. However, the replacement operator
$\widehat\theta_{k,l}$ does not appear in (\ref{eq207})
and (\ref{eq209}). The reason (in the haploid case) is that the
chance that the $k$th individual reproduces through a resampling
event depends only on the fitness difference to a randomly chosen
individual from the population. See also (\ref{eq803}),
(\ref{eq804}) and (\ref{eq620}).
\end{longlist}

Given $\mathbf P_0\in\mathcal M_1(\mathcal M_1(I))$, it was shown in
\citet{EthierKurtz1993} [see also \citet{D93}] that the
$(\mathbf
P_0, \widehat\Omega, \mathcal F)$-martingale problem is
well-posed. We refer to the solution as the (measure-valued)
\textit{Fleming--Viot process with mutation and selection}, FVMS. This
is a strong Markov process with continuous paths and hence a
diffusion.

More general generators were considered in \citet{DawsonMarch1995},
where state-dependent resampling and mutation rates were
allowed. Selection intensities depending on the state of the FVMS
were considered in \citet{DonnellyKurtz1999} and unbounded selection
operators are studied in \citet{EthierShiga2000}. In all these cases
well-posedness of the corresponding martingale problem was shown.
\end{example}
%
%de3.10 #&#
%
\begin{definition}[(Generator of TFVMS)]\label{defgentvms}
We use the same notation as in
Example~\ref{exmFV}. The generator of TFVMS is the linear operator
on $\Pi$ with domain~$\Pi^1$, given by
%
%e3.24 #&#
%
\begin{equation}\label{dgpGir2}
\Omega:= \Omega^{\mathrm{grow}}+\Omega^{\mathrm{res}} +
\Omega^{\mathrm{mut}} + \Omega^{\mathrm{sel}}.
\end{equation}
Here, for $\Phi^{n,\phi} \in\Pi^1_{n}$ the different terms are
given as follows:
\begin{longlist}[(2)]
\item[(1)] We define the \textit{growth operator} by
%
%e3.25 #&#
%
\begin{equation}\label{eqomega1}
\Omega^{\mathrm{grow}}\Phi(\smallu)
:= \langle\nu^{\smallu},\langle\nabla_{\underline{\underline
r}}\phi,
\underline{\underline2} \rangle\rangle
\end{equation}
with
%
%e3.26 #&#
%
\begin{equation}
\label{eq27}
\langle\nabla_{\underline{\underline r}}\phi, \underline
2\rangle= 2 \sum_{1 \leq i<j} \frac{\partial\phi}{\partial
r_{ij}}(\underline{\underline r}, \underline u).
\end{equation}
\item[(2)] We define the \textit{resampling operator} by
%
%e3.27 #&#
%
\begin{equation}
\label{eqomega2}
\Omega^{\mathrm{res}}\Phi(\smallu):=
\frac{\gamma}{2}\sum_{k,l=1}^n \langle
\nu^\smallu,\phi\circ\theta_{k,l} -\phi\rangle
\end{equation}
with $\theta_{k,l}( \underline{\underline r}, \underline u
) = (\underline{ \underline{\widetilde r}},
\widehat\theta_{k,l}(\underline{u}) )$ [recall
$\widehat\theta_{k,l}$ from (\ref{gen2})] and
%
%e3.28 #&#
%
\begin{equation}
\label{pp11b}
\widetilde r_{ij}:=
\cases{r_{ij}, &\quad if $i,j\neq l$,
\cr
r_{i\wedge k, i\vee k}, &\quad if $j=l$,
\cr
r_{j\wedge k, j\vee k}, &\quad if $i=l$.}
\end{equation}
As an example,
%
%e3.29 #&#
%
\begin{equation}
{\fontsize{10.6pt}{12.5pt}\selectfont{
\theta_{1,3}(\underline{\underline r}, \underline u) =
\left(\pmatrix{ 0 & r_{12} & 0 & r_{14} & r_{15} &
\cdots\cr& 0 & r_{12} & r_{24} & r_{25} & \cdots\cr& & 0
& r_{14} & r_{15} & \cdots\cr& & & 0 & r_{45} & \cdots
\cr
& & & & \cdots& \cdots}, (u_1, u_2,
u_1, u_4, u_5,\ldots)\right).}}\hspace*{-36pt}
\end{equation}
\item[(1)] For the \textit{mutation operator}, let $\vartheta,
\beta(\cdot,\cdot)$ be as in Example~\ref{exmFV}, and set
%
%e3.30 #&#
%
\begin{equation}
\label{eqomega3}
\Omega^{ \mathrm{mut}}\Phi(\smallu):= \vartheta\sum_{k=1}^n
\langle\nu^\smallu, B_k\phi\rangle,
\end{equation}
such that
%
%e3.31 #&#
%
\begin{eqnarray}
\label{eq33}
B_k\phi &:=& \beta_k\phi- \phi, \nonumber\\[-8pt]\\[-8pt]
(\beta_k\phi)(\underline{\underline r}, \underline u)
&:=&
\int\phi(\underline{\underline r}, \underline u_k^v)
\beta(u_k, dv).
\nonumber
\end{eqnarray}
\item[(2)] For \textit{selection}, consider
%
%e3.32 #&#
%
\begin{equation}\label{dgpGir4}
\chi'\dvtx I \times I \times\mathbb R_+ \to[0,1]
\end{equation}
with $\chi'(u,v,r) = \chi'(v,u,r)$ for all $u,v\in I, r\in\mathbb
R_+$; recall (\ref{agx1e}). In our main results, we require that
$\chi'\in\overline{\mathcal C}{}^{0,0,1}(I\times I \times\mathbb
R_+)$; that is, $\chi'$ is continuous and continuously differentiable
with respect to its third coordinate. Then with
%
%e3.33 #&#
%
\begin{equation}\label{dgpGir5}
\chi'_{k,l} ( \underline{\underline r}, \underline u):=
\chi'(u_k, u_l, r_{k\wedge l, k\vee l}),
\end{equation}
we set
%
%e3.34 #&#
%
\begin{equation}\label{eqpp15}
\Omega^{ \mathrm{sel}}\Phi(\smallu):= \alpha\sum_{k=1}^n \langle
\nu^\smallu, \phi\cdot\chi'_{k,n+1} - \phi\cdot\chi'_{n+1,n+2}
\rangle.
\end{equation}
If $\chi'(u,v,r)$ does not depend on $r$, and if there is $\chi\dvtx I
\to
[0,1]$ such that
%
%e3.35 #&#
%
\begin{equation}\label{eq208b}
\chi'(u,v,r)=\chi(u)+\chi(v)
\end{equation}
[compare (\ref{eq208})], we say that selection is additive and
conclude that with
%
%e3.36 #&#
%
\begin{equation}
\label{eqpp15c}
\chi_k(\underline{\underline r}, \underline u) = \chi(u_k).
\end{equation}
We obtain
%
%e3.37 #&#
%
\begin{equation}\label{eqpp15b}
\Omega^{ \mathrm{sel}}\Phi(\smallu):= \alpha\cdot
\sum_{k=1}^n \langle\nu^\smallu, \phi\cdot\chi_{k} -
\phi\cdot\chi_{n+1} \rangle.
\end{equation}
\end{longlist}
\end{definition}

Now, we are ready to give our first main result.
%
%th3.1 #&#
%
\begin{theorem}[(Martingale problem is well posed)]\label{T1}
Let $\mathbf P_0\in\mathcal M_1(\mathbb U^I)$, $\Pi^1$
be as in (\ref{eq201}) and $\Omega$ as in (\ref{dgpGir2}).
\begin{enumerate}[(2)]
\item[(1)] The $(\mathbf P_0, \Omega, \Pi^1)$-martingale is well posed.
The unique solution $\mathcal U:=(\mathcal U_t)_{t\geq0}$ is
called the tree-valued Fleming--Viot dynamics with mutation and
selection (TFVMS).
\item[(2)] The process $\mathcal U$ has the following properties:
\begin{enumerate}[(a)]
\item[(a)]$\mathbf P(t\mapsto\mathcal U_t \mbox{ is continuous})=1$;
\item[(b)]$\mathbf P(\mathcal U_t \in\mathbb U^I_c \mbox{ for all }t>0)=1$;
\item[(c)]$\smallu\mapsto\mathbf E[f(\mathcal U_t)|\mathcal U_0=\smallu
]$ is
continuous for all $f\in\overline{\mathcal C}(\mathbb U^I)$, that
is, $\mathcal
U$ has the Feller property;
\item[(d)]$\mathcal U$ is strong Markov;
\item[(e)] for $\Phi=\Phi^{n,\phi}\in\Pi^1$, the quadratic variation of the
process $\Phi(\mathcal U) = (\Phi(\mathcal U_t))_{t\geq0}$ is given
by
%
%e3.38 #&#
%
\begin{equation} \label{eq222} [ \Phi(\mathcal U) ]_t = \gamma
\sum_{k,l=1}^n \int_0^t \langle\nu^{\mathcal U_s}, \phi\cdot
(\phi\circ\rho_1^n)\cdot\theta_{k,n+l} - \phi\cdot
(\phi\circ\rho_1^n)\rangle \,ds,
\end{equation}
where
%
%e3.39 #&#
%
\begin{equation}
\label{eqrho}
\rho_1^n(\underline{\underline r}, \underline u) =
((r_{i+n,j+n})_{1\leq i<j}, (u_{i+n})_{i\geq1})
\end{equation}
denotes the $n$-shift of the sample sequence.
\end{enumerate}
\end{enumerate}
\end{theorem}
%
%re3.11 #&#
%
\begin{remark}[(Mark function)]\label{remmark}
We will show in forthcoming work that states of the
TFVMS only take special forms:
\begin{longlist}[(2)]
\item[(1)] Consider an mmm-space $\smallu=
\overline{(U,r,\mu)} \in\mathbb U^I$. We say that $\smallu$ has a
\textit{mark function} if there is an $U$-valued random variable
$X$ and $\kappa\dvtx U\to I$ [both measurable with respect to the
Borel-$\sigma$-algebra of $(U,r)$] such that $(X, \kappa(X))$
has the distribution $\mu$. In other words, $\smallu$ has a
mark function if there is a measurable function $\kappa\dvtx U\to I$
with
%
%e3.40 #&#
%
\begin{equation}
\label{eqmark100}
\mu(dx,du) = ((\pi_U)_\ast\mu)(dx) \cdot
\delta_{\kappa(x)}(du).
\end{equation}
As argued in Remark~\ref{remtrees1}, the TMMMS always admits
states in $\mathbb U^I$ which have a mark function. It turns out
that the same holds for the TFVMS as well.
\item[(2)] Another path property we will address are atoms of the measure
$\mu$. Consider the TFVMS $\mathcal U = (\mathcal U_t)_{t\geq0}$
with $\mathcal U_t = \overline{(U,r,\mu)}$. Then, $(\pi_U)_\ast
\mu$ has an atom if and only if $\mu^{\otimes2}\{(x,y)\dvtx r(x,y)=0\}
>0$. We
shall show that $\mathcal U$ only takes values in the space of
mmm-spaces $\smallx= \overline{(X,r,\mu)}$ with the property that
$(\pi_U)_\ast\mu$ has no atoms. Note that only the projection
$(\pi_U)_\ast\mu$ can be free of atoms since it is well known
that $(\pi_I)_\ast\mu$ is atomic for all $t\geq0$, almost
surely; see, for example, Theorem 10.4.5 in \citet{EthierKurtz86}.
\end{longlist}
\end{remark}

%s3.3 #&#
\subsection{Girsanov theorem for the TFVMS}
\label{ssGirs}
One possibility to establish the existence and uniqueness of
martingale problems and to analyze its properties is to show that
solutions of different martingale problems are absolutely continuous
to each other for finite time horizons. Uniqueness as well as several
other properties (e.g., path properties) then carry over from one
martingale problem to the other. The densities of the solutions of the
martingale problems are calculated by the Cameron--Martin--Girsanov
theorem for real-valued semimartingales [see Theorem 16.19 in
\citet{Kallenberg2002}] and Dawson's Girsanov theorem for
measure-valued processes [\citet{D93}, Section 7.2]. Here, we
carry out
the corresponding program for TFVMS by considering two martingale
problems with different selection strength.
%
%re3.12 #&#
%
\begin{remark}[(Notation)]%
For $\alpha\in\mathbb R_+$, we write $\Omega_\alpha$ and
$\Omega_\alpha^{\mathrm{sel}}$ for the operators defined in
(\ref{dgpGir2}) and (\ref{eqpp15}), respectively, when we want to
stress the value of the selection coefficient $\alpha$.
\end{remark}
%
%th3.2 #&#
%
\begin{theorem}[(Girsanov Transform for the TFVMS processes)]\label{T2}
Let $\alpha,\alpha'\in\mathbb R_+$, $\mathbf
P_0\in\mathcal M_1(\mathbb U^I)$, and using $\chi'_{1,2}$ from
(\ref{dgpGir5}) define $\Psi\in\Pi^1$ by
%
%e3.41 #&#
%
\begin{equation}
\label{eqPsi}
\Psi(\smallu):= \frac{\alpha' - \alpha}{\gamma} \cdot\langle
\nu^{\smallu}, \chi'_{1,2}\rangle.
\end{equation}
Let $\mathbf P \in\mathcal M_1(\mathcal C_{\mathbb
U^I}(\mathbb R_+))$ be a solution of the $(\mathbf P_0,
\Omega_\alpha, \Pi^1)$-martingale problem, $\mathcal U=(\mathcal
U_t)_{t\geq0}$ the canonical process with respect to $\mathbf P$,
$(\mathcal F_t)_{t\geq0}$ its canonical filtration and
%
%e3.42 #&#
%
\begin{equation}\label{eqap4}
\mathcal M = (M_t)_{t\geq0} = \biggl(\Psi(\mathcal U_t) -
\Psi(\mathcal U_0) - \int_0^t \Omega_\alpha\Psi(\mathcal U_s)
\,ds\biggr)_{t\geq0}.
\end{equation}
Then, $\mathcal M$ is a $\mathbf P$-martingale and the probability measure
$\mathbf Q$, defined by
%
%e3.43 #&#
%
\begin{equation}\label{eqdQdP}
\frac{d\mathbf Q}{d\mathbf P}\bigg|_{\mathcal F_t} = e^{M_t -
(1/2) [ \mathcal M]_t}
\end{equation}
solves the $(\mathbf P_0, \Omega_{\alpha'}, \Pi^1)$-martingale problem.
\end{theorem}

%s3.4 #&#
\subsection{Convergence of Moran models}
\label{ssconv}
Our next task is to relate the Fleming--Viot process to the finite
population models and their evolving genealogies on the level of
trees, that is, mmm-spaces.
%
%de3.13 #&#
%
\begin{definition}[(TMMMS)] \label{defTMMMS2}
Recall the process $(U_N, r^N_t, \mu^N_t)_{t\geq
0}$ from Definition~\ref{defTMMMS}, started in a random mmm-space
$(U_N, r^N_0, \mu^N_0)$. The fitness function is either given as in
Definition~\ref{defgraph} or by (\ref{agx1e}). The
\textit{tree-valued Moran model with mutation and selection} (\textit{TMMMS})
is given by
%
%e3.44 #&#
%
\begin{equation}\label{eq497}
\mathcal U^N = (\mathcal U_t^N)_{t\geq0},\qquad \mathcal U_t^N =
\overline{(U_N, r^N_t, \mu^N_t)}.
\end{equation}
\end{definition}
%
%th3.3 #&#
%
\begin{theorem}[(Convergence to TFVMS)]\label{T3}
Let $\mathcal U^N$ be the TMMMS, started in~$\mathcal U_0^N$, and
$\mathcal U$ be the TFVMS, started in $\mathcal U_0$. If $\mathcal
U^N_0\xRightarrow{N\to\infty}\mathcal U_0$, weakly
with respect to the Gromov-weak topology, then
%
%e3.45 #&#
%
\begin{equation}\label{ag2}
\mathcal U^N \xRightarrow{N\to\infty} \mathcal U,
\end{equation}
weakly with respect to the Skorohod topology on ${\mathcal
D}_{\mathbb{U}^I}([0,\infty))$.
\end{theorem}

%s3.5 #&#
\subsection{Long-time behavior}
\label{sslongtime}
We now determine under which conditions the TFVMS has a unique
invariant measure and is ergodic.
This is not always the case, since already for the measure-valued
process there are examples where the process is nonergodic. (A
trivial example is $\vartheta=0$, but cases when mutation has several
invariant distributions are also possible.)

Recall\vspace*{1pt} the measure-valued Fleming--Viot process $\zeta=(\zeta
_t)_{t\geq
0}$ from Example~\ref{exmFV} and the projection $\pi_I$ on $I$ from
Remark~\ref{remnotation}. Given $\mathcal U_t =
\overline{(U_t,r_t,\mu_t)}$, $t\geq0$, define the process
%
%e3.46 #&#
%
\begin{equation} \label{addx2}
\widetilde\zeta
:= (\widetilde\zeta_t)_{t\geq0}:= ((\pi_I)_\ast\mu_t)_{t\geq0},
\end{equation}
and note that $(\widetilde\zeta_t)_{t\geq0} \stackrel d =
(\zeta_t)_{t\geq0}$ if $\chi'(u,v,r) = \widehat\chi'(u,v)$, that is, if
the fitness is independent\vspace*{2pt} of the genealogical distance.
Hence, existence of a unique equilibrium for $\widetilde\zeta$ is
always implied by existence of a unique equilibrium for $\mathcal U$.
Theorem~\ref{T4} shows that the opposite is also true. The proof of
Theorem~\ref{T4} is based on duality, introduced in Section
\ref{sduality}.
%
%th3.4 #&#
%
\begin{theorem}[(Long-time behavior)] \label{T4}
\textup{(a)} Let $\mathcal U=(\mathcal U_t)_{t\geq0}$ be
the TFVMS
with $\mathcal U_0=\smallu$ and $\widetilde\zeta$ be as above. Then
there exists an $\mathbb U_c^I$-valued random variable $\mathcal
U_{\infty}$ with
%
%e3.47 #&#
%
\begin{equation}\label{eq498}
\mathcal U_t \xRightarrow{t\to\infty} \mathcal U_{\infty},
\end{equation}
if and only if $\widetilde\zeta$ has a unique equilibrium
distribution.\vspace*{-4pt}

\begin{longlist}
\item[(b)] The law of $\mathcal U_{\infty}$ is the unique invariant distribution
of $\mathcal U$. It depends on all the model parameters but is independent
of the initial state.
\end{longlist}

In particular, if mutation and selection are present, $\vartheta>0$,
$\alpha>0$ and mutation has a parent-independent component (i.e.,
(\ref{eq1001}) holds for some $z\in(0;1]$), then (\ref{eq498})
holds.
\end{theorem}
%
%re3.14 #&#
%
\begin{remark}[(Conditions for ergodicity of $\zeta$)] \label
{remconvmoreg}%
Various results about ergodicity of the measure-valued Fleming--Viot
process have been obtained, which carry over to the TFVMS by
Theorem~\ref{T4}. For example, under neutral evolution, $\alpha=0$
(or $\chi'=0$), ergodicity has been shown if the Markov pure jump
process on $I$ with generator (\ref{eq205}) has a unique
equilibrium distribution [\citet{D93}]. In the case $\alpha>0$ and
$\chi'\neq0$, ergodicity of $\zeta$ in the
case of no parent-independent component in the mutation operator
[i.e., $z=0$ in (\ref{eq1001})] have been shown in
\citet{EthierKurtz1998} using coupling techniques. Using different
techniques, \citet{EthierKurtz1998} also prove an ergodic theorem
for a version of the infinitely-many-alleles model with symmetric
overdominance. In \citet{ITATSUSeiichi2002-03-31} a perturbative
approach is used to prove ergodicity of measure-valued Fleming--Viot
processes with weak selection under ergodicity assumption on the
mutation process. In \citet{DGsel} a set-valued dual [see also
\citet{DGsel0}] allows one to prove ergodic theorems, even if the
population is distributed on geographic sites if mutation has a
parent-independent part.
\end{remark}

%s3.6 #&#
\subsection{Application: Distance between two individuals}
\label{ssdist}
It is widely believed that genealogical distances under additive
selection are smaller than under neutrality. The heuristics are that
beneficial alleles spread quicker through the population than neutral
ones by their fitness advantage. Hence, after the allele has spread,
randomly chosen individuals have a more recent last common ancestor
than under neutrality. In other words, genealogical distances are
shorter. However, shorter distances under selection are actually
difficult to ascertain, because there is no monotonicity of
genealogical distances in the selection coefficient $\alpha$ since the
state of the process is due to an intricate interaction between the
mutation and the selection. (Note that, as $\alpha\to\infty$ the
genealogies look essentially neutral since fixation on the fittest
types takes place.) We cannot prove that genealogical distances are
shorter under additive selection yet, but we make a first step in that
direction.

Namely, we apply our machinery to the comparison of pairwise
genealogical distances in the selective and in the neutral case. We
give a concrete example how genealogical distances change under
selection in the case of two alleles and if the selection coefficient
is small.

In order to make the comparison of distances precise, we proceed as
follows. Let $\mathcal U_{\infty}^\alpha$ be the unique invariant
$\mathbb U^I$-valued random variable from Theorem~\ref{T4} (if it
exists). Let $R^\alpha_{12}$ denote the distance of two randomly
chosen points from $\mathcal U_{\infty}^\alpha$. Hence,
%
%e3.48 #&#
%
\begin{equation}\label{agx2}
R^\alpha_{12} \mbox{ has distribution } A\mapsto
\mathbf E[(r_{12})_\ast\nu^{\mathcal U_{\infty}^\alpha}(A)]
\end{equation}
for Borel-sets $A\subseteq\mathbb R_+$, and $r_{12}$ denotes the
function $\underline{\underline r}\mapsto r_{12}$. In other words, the
distribution of $R_{12}$ is the first moment measure of the random
probability distribution $(r_{12})_\ast\nu^{\mathcal
U_{\infty}^\alpha}$. For $\alpha>0$, the issue is now to decide
whether $R^\alpha_{1,2} < R^0_{1,2}$ in stochastic order.
%
%re3.15 #&#
%
\begin{remark}[(Laplace-transform order and Landau symbol)]
\label{remR12}
(1) For two random variables $X, Y$, we say that $X\leq Y$ in the
\textit{Laplace-transform order} if $\mathbf E[e^{-\lambda X}] \geq
\mathbf E[e^{-\lambda Y}]$ for all $\lambda>0$. Note that this
does not necessarily imply that $X\leq Y$ stochastically.

(2) In the next theorem, we use the Landau symbol $\mathcal O(\cdot)$.
In particular, for functions $g$ and $h$, both\vspace*{1pt} depending
on $\alpha$, we write $g(\alpha) = h(\alpha) + \mathcal O(\alpha^3)$ as
$\alpha\to0$ if $\limsup_{\alpha\to
0}|(g(\alpha)-h(\alpha))/\alpha^3|<\infty$.
\end{remark}

The following theorem is dealing with the same case as
Example~\ref{extwotypes}.
%
%th3.5 #&#
%
\begin{theorem}[(Distance of two randomly sampled individuals)]
\label{T5}
Let $I=\{\bullet, {{\darkgreybullet}}\}$,
$\chi(u)=1_{\{u=\bullet\}}$. Assume that the mutation rate is
$\overline\vartheta/2$ and for the mutation stochastic kernel
$\beta(\cdot,\cdot)$,
%
%e3.49 #&#
%
\begin{equation}\label{eq499}
\frac{\overline\vartheta}{2}\cdot\beta(u,dv) =
\frac{\vartheta_\bullet}2\ind{v={\darkgreybullet}} +
\frac{\vartheta_{{\darkgreybullett}}} 2\ind{v={\bullet}}
\end{equation}
for some $\vartheta_\bullet, \vartheta_{\darkgreybullett}>0$
with $\overline\vartheta= \vartheta_\bullet+
\vartheta_{\darkgreybullett}$, that is, $\bullet$ mutates to
${\darkgreybullet}$ at rate $\vartheta_\bullet/2$ and from
${\darkgreybullet}$ to $\bullet$ at rate
$\vartheta_{{\darkgreybullett}}/2$. In addition, selection is
additive, that is, (\ref{eqpp15b}) holds for some $\alpha>0$ and
$\mathcal U^\alpha_{\infty}:=\mathcal U_{\infty}$ is as in
Theorem~\ref{T4}. [Note that $\beta(u,dv)$ does not depend on~$u$,
and therefore (\ref{eq1001}) holds with $z=1$.] Let
$R_{12}^\alpha$ be as in (\ref{agx2}).

Then as $\alpha\to0$, for $\lambda>0$,
%
%e3.50 #&#
%
\begin{equation}
\label{eq501}
\mathbf E[ e^{-\lambda R^\alpha_{12}}] =
\frac{\gamma}{\gamma+2\lambda} + f \alpha^2 + \mathcal
O(\alpha^3),
\end{equation}
where
$f:=f(\gamma,\vartheta_\bullet,\vartheta_{\darkgreybullett},\lambda)$
is given by
\[
f = \frac{8\gamma
\vartheta_\bullet\vartheta_{{\darkgreybullett}}(2\gamma+
2\lambda+ \overline\vartheta) \lambda
}{\overline\vartheta(\gamma+ \overline\vartheta)(\gamma+
2\lambda+ \overline\vartheta)(6\gamma+ 2\lambda+
\overline\vartheta)(\gamma+ 2\lambda)^2(6\gamma+ 4\lambda+
\overline\vartheta)}.
\]
In particular, $R_{12}^\alpha\leq R_{12}^0$ in the
Laplace-transform order for small $\alpha$ and
%
%e3.51 #&#
%
\begin{equation}
\label{eq502}
\mathbf E[R^\alpha_{12}] = \frac{1}{\gamma}\biggl( 2 -
\frac{8\vartheta_\bullet\vartheta_{\darkgreybullett}(2\gamma
+ \overline\vartheta)}{\overline\vartheta(\gamma+
\overline\vartheta)^2(6\gamma+
\overline\vartheta)^2}\alpha^2\biggr) + \mathcal O(\alpha^3).
\end{equation}
\end{theorem}
%
%re3.16 #&#
%
\begin{remark}[{[Distances under selection and connection to
\citet{KroneNeuhauser1997}]}]\label{remR12NK} %
(1) Under neutrality, $R_{12}^0$ is
exponentially distributed with rate $\gamma/2$, thus $\mathbf
E[e^{-\lambda R^{0}_{12}}] = \frac{\gamma}{\gamma+2\lambda}$. Note
that for small $\alpha$, the Laplace transform differs from the
neutral case only in second order in $\alpha$. The fact that the
first order is the same as under neutrality was already obtained
by \citet{KroneNeuhauser1997} for a finite Moran model. Our proof
in Section~\ref{SPPP3} can be extended to obtain higher order
terms. However, it is an open problem to show $R_{12}^\alpha<
R_{12}^0$ stochastically for small $\alpha$ since the
Laplace-transform order is weaker than the stochastic order.

(2) The order $R^{\alpha_2}_{12} < R^{\alpha_1}_{12}$ cannot be
expected to hold for all values $\alpha_1<\alpha_2$. The reason
is that for large values of $\alpha$, most individuals in the
population carry the fit type $\bullet$ and therefore, the
genealogy is close to the Kingman coalescent with
pair-coalescence-rate $\gamma$.
\end{remark}

\textit{Outline of the proof section}: before we come to the
proofs of the Theorems~\ref{T1}--\ref{T5}, we develop three main
technical tools. These are an analysis of the generator for the TFVMS
(Section~\ref{schar}), duality (Section~\ref{sduality}) and an
investigation of the tree-valued Moran model with mutation and
selection (Section~\ref{secmart-probl-fixed}). The proofs of Theorems
\ref{T1}--\ref{T4} are given in Section~\ref{SPPP2} and the
application, Theorem~\ref{T5}, is proved in Section~\ref{SproofAppl}.

%s4 #&#
\section{Infinitesimal characteristics}
\label{schar}
The TFVMS is a strong Markov process with continuous paths, and
therefore may be called a \textit{tree-valued diffusion}. Since
generators of diffusions are typically second order differential
operators, it
is natural to ask in which sense the same is true for the TFVMS with
the generator $\Omega$ from (\ref{dgpGir2}). Here it is useful to
work with an abstract
concept of order of linear operators. The distinction of first
and second order terms is also the key to the proof of the
Girsanov-type result, Theorem~\ref{T2}.

%s4.1 #&#
\subsection{First and second order operators}
We recall some basic facts about linear operators, which are related
to differential operators. For their connection to Markov processes
see \citet{FukushimaStroock1986} and Section~VIII.3 of
\citet{RevuzYor2010}.
%
%de4.1 #&#
%
\begin{definition}[(First and second order operators)]\label{deffirstO}
Let $\Omega$ be a linear operator with domain
$\mathcal D$ and $\Pi\subseteq\mathcal D$ an algebra. We say that
$\Omega$ is \textit{first order} (with respect to $\Pi$) if for all
$\Phi\in\Pi$,
%
%e4.1 #&#
%
\begin{equation}\label{dgpGir14}
\Omega\Phi^2 - 2 \Phi\cdot\Omega\Phi= 0.
\end{equation}
We say that $\Omega$ is \textit{second order} if it is not first order,
and for all $\Phi\in\Pi$
%
%e4.2 #&#
%
\begin{equation}\label{dgpGir15}
\Omega\Phi^3 + 3 \Phi^2 \cdot\Omega\Phi- 3 \Phi\cdot\Omega
\Phi^2 = 0.
\end{equation}
\end{definition}
%
%re4.2 #&#
%
\begin{remark}[(Diffusions in $\mathbb R^d$ and higher order
operators)]%
\label{remdiff}
(1) A diffusion process on $\mathbb R^d$ has a generator
%
%e4.3 #&#
%
\begin{equation}
\label{eqdiffRd}\qquad
\Omega= \Omega_1 + \Omega_2,\qquad \Omega_1:= \sum_{i=1}^d
\mu_i(\underline x) \,\frac{\partial}{\partial x_i},\qquad
\Omega_2 = \sum_{i,j=1}^d \sigma^2_{ij}(\underline x)
\,\frac{\partial^2}{\partial x_i \,\partial x_j}
\end{equation}
with domain $\mathcal D = \mathcal C^2_b(\mathbb R^d)$, for a
vector $(\mu_i)_{i=1,\ldots,d}$ and a positive definite matrix
$(\sigma_{ij})_{1\leq i,j\leq d}$, which are continuous functions
on $\R^d$. It can be easily checked that $\Omega_1$ is a first
order operator, and $\Omega_2$ is a second order operator with
respect to $\mathcal D$, according to
Definition~\ref{deffirstO}. Hence, the above definitions of first
and second order operators extend the usual notions for
differential operators.\vspace*{-6pt}

\begin{longlist}[(2)]
\item[(2)] The operator defined through the left-hand side of
(\ref{dgpGir14}) is connected to the \textit{square field operator},
also called \textit{op\'erateur carr\'e du champ}, which
is given by
%
%e4.4 #&#
%
\begin{equation}
\Gamma(\Phi, \Psi):= \Omega\Phi\Psi- \Phi\Omega\Psi- \Psi
\Omega\Phi.
\end{equation}
In particular, a straightforward calculation (similar to the
proof of Lemma~\ref{lsecond0} below) shows that $\Omega$ is
second order if and only if $\Gamma$ is a \textit{derivation} [in the
sense of
\citet{BakryEmery1985}, i.e., $\Gamma(\Phi\Psi, \Lambda) = \Phi
\Gamma(\Psi, \Lambda) + \Psi\Gamma(\Phi, \Lambda)$ for all $\Phi,
\Psi, \Lambda\in\Pi$].
\item[(3)] Typically, higher order operators do not arise if $\mathcal
D$ is a
subset of continuous functions, and $\Omega$ is the generator of a Markov
process $(X_t)_{t\geq0}$ with continuous paths. The reason is that
$(\Phi(X_t) - \int_0^t \Omega\Phi(X_s)\,ds)_{t\geq0}$ is a
continuous martingale and therefore $(\Phi(X_t))_{t\geq0}$ can only have
quadratic variation, which means that $\Omega$ is at most second order;
see Proposition~\ref{Pqv} below.
\end{longlist}
\end{remark}

First and second order operators satisfy some further relations when
applied to products or powers, which we derive next.
%
%le4.3 #&#
%
\begin{lemma}[(First order operators)]\label{lfirstO}
If a linear operator $\Omega$ is first order with respect
to the algebra $\Pi$, then
%
%e4.5 #&#
%
\begin{equation}\label{dgpGir16}
\Omega(\Phi\cdot\Psi) - \Phi\cdot\Omega\Psi- \Psi\cdot
\Omega\Phi= 0.
\end{equation}
In particular, (\ref{dgpGir15}) holds.
\end{lemma}
\begin{pf}
Equation (\ref{dgpGir16}) follows immediately once we compute $\Omega
(\Phi
+ \Psi)^2$ and use linearity of $\Omega$. Furthermore, (\ref{dgpGir15})
follows by using $\Psi=\Phi^2$ and (\ref{dgpGir14}) in~(\ref{dgpGir16}).
\end{pf}
%
%le4.4 #&#
%
\begin{lemma}[(Second order operators)]\label{lsecond0}
If a linear operator $\Omega$ is first or second
order with respect to the algebra $\Pi$, then for all $\Phi,\Psi\in
\Pi$
%
%e4.6 #&#
%
\begin{equation}\label{dgpGir16b}
\Omega\Psi\Phi^2 + 2 \Psi\Phi\cdot\Omega\Phi+
\Phi^2\cdot\Omega\Psi- \Psi\cdot\Omega\Phi^2 - 2 \Phi\cdot
\Omega\Psi\Phi= 0.
\end{equation}
In particular, for any $\Phi\in\Pi$,
%
%e4.7 #&#
%
\begin{equation}\label{dgpGir16c}
\Omega\Phi^4 + 8\Phi^3\cdot\Omega\Phi- 6
\Phi^2\cdot\Omega\Phi^2=0.
\end{equation}
\end{lemma}
\begin{pf}
Applying %the definition of a second order operator
(\ref{dgpGir15})
to $(\Psi+\Phi)^3$ and $(\Psi-\Phi)^3$, and summing up, gives
%
%e4.8 #&#
%
\begin{eqnarray}
\label{dgpGir16ca}\quad
0 & = & 2 \Omega\Psi^3 + 6\Omega\Psi\Phi^2 + 6\Psi^2 \cdot
\Omega\Psi+ 12 \Psi\Phi\cdot\Omega\Phi+ 6
\Phi^2\cdot\Omega\Psi\nonumber\\
&&{}
- 6\Psi\cdot\Omega\Psi^2 - 6\Psi\cdot\Omega\Phi^2 - 12
\Phi\cdot\Omega\Psi\Phi\\
& = & 6\Omega\Psi\Phi^2 + 12
\Psi\Phi\cdot\Omega\Phi+ 6 \Phi^2\cdot\Omega\Psi-
6\Psi\cdot\Omega\Phi^2 - 12 \Phi\cdot\Omega\Psi\Phi,
\nonumber
\end{eqnarray}
which implies (\ref{dgpGir16b}). To show (\ref{dgpGir16c}), we use
(\ref{dgpGir16b}) with $\Psi=\Phi^2$ and obtain
%
%e4.9 #&#
%
\begin{eqnarray}
\label{dgpGir16d}\quad
0 & = & \Omega\Phi^4 + 2 \Phi^3 \cdot\Omega\Phi+
\Phi^2\cdot\Omega\Phi^2 - \Phi^2\cdot\Omega\Phi^2 - 2 \Phi
\cdot\Omega\Phi^3 \nonumber\\[-8pt]\\[-8pt]
& = &\Omega\Phi^4 + 8\Phi^3\cdot
\Omega\Phi- 6 \Phi^2\cdot\Omega\Phi^2,
\nonumber
\end{eqnarray}
since $\Omega$ is at most second order.
\end{pf}

%s4.2 #&#
\subsection{Order of operators: Application to Markov processes}
In this subsection we use the concepts of the last subsection to
compute processes of quadratic variation and covariation for
functionals of a Markov process.
%, using the generator from the martingale problem from
% \eqref{dgpGir2}.
%
%pr4.5 #&#
%
\begin{proposition}[(Path continuity\vspace*{1pt} of second order martingale
problems)]\label{Pqv}
Let $E$ be a Polish space, $\Omega= \Omega^{(1)} +
\Omega^{(2)}$ be a linear operator on $\mathcal B(E)$ with domain
$\mathcal D\subseteq\overline{\mathcal C}(E)$, where $\Omega^{(1)}$
is a
first order operator, and $\Omega^{(2)}$ is a second order
operator. Assume that $\mathcal D$ contains a countable algebra
$\Pi$ that separates points in~$E$.

Assume that $\mathcal X = (X_t)_{t\geq0}$ is a solution of the
$(\mathbf P_0,\Omega, \mathcal D)$-martingale problem for $\mathbf
P_0\in\mathcal M_1(E)$ (with paths in $\mathcal D_E([0,\infty))$).
Then, $\mathcal X$ has the following path properties:
\begin{longlist}[(2)]
\item[(1)] $\mathcal X$ has paths in $\mathcal C_E([0,\infty))$, almost
surely;
\item[(2)] for $\Phi\in\Pi$, the process $\Phi(\mathcal X) =
(\Phi(X_t))_{t\geq0}$ is a continuous semimartingale with
quadratic variation given by
%
%e4.10 #&#
%
\begin{equation} \label{eq2} [ \Phi(\mathcal X) ]_t = \int_0^t
\Omega^{(2)} \Phi^2(X_s) - 2\Phi(X_s) \cdot\Omega^{(2)}
\Phi(X_s)\,ds.
\end{equation}
\end{longlist}
\end{proposition}
%
%co4.6 #&#
%
\begin{corollary}[(Covariation)]%
\label{corcov}
Under the assumptions of Proposition~\ref{Pqv}, let
$\Phi,
\Psi\in\Pi$. The covariation of the processes\vadjust{\goodbreak} $\Phi(\mathcal X)=
(\Phi(X_t))_{t\geq0}$ and $\Psi(\mathcal X)= (\Psi(X_t))_{t\geq0}$ is
given by
\begin{eqnarray*}
[\Phi(\mathcal X), \Psi(\mathcal X)]_t &=& \int_0^t \Omega^{(2)}
(\Phi\Psi)(X_s) - \Phi(X_s) \Omega^{(2)} \Psi(X_s)\\[-2pt]
&&\hspace*{12.5pt}{} - \Psi(X_s)
\Omega^{(2)} \Phi(X_s)\,ds.
\end{eqnarray*}
\end{corollary}
\begin{pf}
This is a simple consequence of (\ref{eq2}) and polarization.
\end{pf}
%
%re4.7 #&#
%
\begin{remark}[{[Connection to \citet{BakryEmery1985}]}] The path
continuity of functionals of $\mathcal X$ was already studied by
\citet{BakryEmery1985} using similar techniques. They show that
$(\Phi(X_t))_{t\geq0}$ is continuous for all $\Phi\in\Pi$ if and
only if the
square field operator is a derivative [or if and only if $\Omega$ is a second
order operator; see Remark~\ref{remdiff}, item (2)]. We extend their
result, since Proposition~\ref{Pqv} gives a sufficient condition
for path continuity of the process $\mathcal X$ (rather than of
functionals of $\mathcal X$). In order to show continuity of
$\mathcal X$, we must require that the domain of $\Omega$ contains a
countable algebra that separates points.
\end{remark}
%
%re4.8 #&#
%
\begin{remark}[(Usual assumption on $\mathcal D$)]\label{remsepsep}
Usually, in order to guarantee that a solution of a
martingale problem has paths in $\mathcal D_E([0,\infty))$, one
requires that $\mathcal D(\Omega)$ is separating and contains a
countable subset that separates points; see \citet{EthierKurtz86},
Theorem 4.3.6.
\end{remark}
\begin{pf*}{Proof of Proposition~\ref{Pqv}}
The proof consists of three steps. First, we show that
$\Phi(\mathcal X)$ is continuous, almost surely, for all
$\Phi\in\Pi$. To have a self-contained proof, we give the full
argument here. However, note that continuity of $\Phi(\mathcal X)$
follows from Proposition 2 in \citet{BakryEmery1985}. Second, we
establish that $t\mapsto X_t$ is almost surely continuous. Third, we
prove (\ref{eq2}).\vspace*{10pt}

\textit{Step} 1: \textit{$\Phi(\mathcal X)$ has continuous paths}: we use
similar arguments as in the proof of Theorem 1.1 and Corollary 1.2
in \citet{FukushimaStroock1986} as well as Kolmogorov's criterion
[e.g., Proposition 3.10.3 in \citet{EthierKurtz86}]. Setting $\Psi_y(x)
:= \Phi(x)-\Phi(y)$ and using that $\mathcal X$ solves the
martingale problem for $\Omega$, we see that
%
%e4.11 #&#
%
\begin{equation}\label{dgpGir210}
\mathbf E\bigl[\bigl(\Phi(X_t)-\Phi(X_s)\bigr)^2\bigr] = \mathbf
E[\Psi^2_{X_s}(X_t)] = \int_s^t \mathbf
E[\Omega\Psi^2_{X_s}(X_r)]\,dr \leq C(t-s)\hspace*{-35pt}
\end{equation}
for some $C<\infty$ by the boundedness of $\Omega\Psi^2$. Moreover,
by Lemma~\ref{lsecond0}, (\ref{dgpGir16c}), using (\ref{dgpGir210})
and some $C'<\infty$,
%
%e4.12 #&#
%
\begin{eqnarray}
\label{dgpGir20}
&&
\mathbf E\bigl[\bigl(\Phi(X_t)-\Phi(X_s)\bigr)^4\bigr] \nonumber\\
&&\qquad= \mathbf
E[\Psi^4_{X_s}(X_t)] \nonumber\\
&&\qquad=
\int_s^t \mathbf E\bigl[
\Psi^2_{X_s}(X_r) \bigl( 6 \Omega\Psi^2_{X_s}(X_r) - 8
\Psi_{X_s}(X_r)\cdot\Omega\Psi_{X_s}(X_r)\bigr)\bigr]\,dr \\
&&\qquad\leq
C' \int_s^t \mathbf E\bigl[\bigl(\Phi(X_r)-\Phi(X_s)\bigr)^2\bigr]
\,dr \leq C' \int_s^t (r-s)\,dr \nonumber\\
&&\qquad\leq C' (t-s)^2,
\nonumber
\end{eqnarray}
and continuity of $\Phi(\mathcal X)$ follows.

Before we carry the continuity of $t\mapsto\Phi(X_t)$ for all
$\Phi\in\Pi$ over to continuity of $t\mapsto X_t$, we recall a basic
topological fact:
%
%re4.9 #&#
%
\begin{remark}\label{remsepsep2}
If $\Pi\subseteq\mathcal C(E)$ separates points and
$x,x_1,x_2,\ldots\in K$, where $K\subseteq E$ is compact. Then,
$x_n\stackrel{n\to\infty}{\longrightarrow} x$ in $E$ if and only if
$\Phi(x_n)\stackrel{n\to\infty}{\longrightarrow} \Phi(x)$ for all
$\Phi\in\Pi$.

The direction ``$\Rightarrow$'' is trivial, since all $\Phi$'s are
continuous. For ``$\Leftarrow$,'' note that $\{x_1,x_2,\ldots\}$ is
relatively compact by assumption. Take any convergent subsequence
$x_{n_k}\stackrel{k\to\infty}{\longrightarrow} y$. Clearly, for all
$\Phi\in\Pi$,
we have $\Phi(y) = \break\lim_{k\to\infty} \Phi(x_{n_k}) =
\lim_{n\to\infty} \Phi(x_{n}) = \Phi(x)$ and hence, $x=y$ since
$\Pi$ separates points.
\end{remark}

\textit{Step} 2: \textit{$\mathcal X$ has continuous paths}: next
we show that $t \mapsto X_t$ is continuous as a function on
$[0,T]\cap\mathbb Q$ for all $T>0$. Since $E$ is Polish, $\mathbf P$
is regular and we can choose an increasing sequence of compact subsets
of $K_1, K_2,\ldots\subseteq E$ with
%
%e4.13 #&#
%
\begin{equation}
\mathbf P(X_t\in K_n \mbox{ for all }0\leq t\leq T) > 1-\frac1n.
\end{equation}
Then set
%
%e4.14 #&#
%
\begin{equation}
\Omega_n:= \{\omega\dvtx X_t(\omega)\in K_n \mbox{ for all }0\leq
t\leq T\}.
\end{equation}
Moreover, take $\Omega'$ with $\mathbf P(\Omega')=1$ and
$\Phi(\mathcal X)$ is continuous on $\Omega'$ for all
$\Phi\in\Pi$. Set $\widetilde\Omega:= \Omega'\cap
\bigcup_{n=1}^\infty\Omega_n$, and note that this set has
probability $1$.

Let $\omega\in\Omega'\cap\Omega_n$ for some $n$ and $t\in\mathbb Q
\cap[0,T]$. Then, for any $t_1,t_2,\ldots$ with
$t_k\stackrel{k\to\infty}{\longrightarrow} t$, $X_{t_1}(\omega),
X_{t_2}(\omega),\ldots\in K_n$ we have
$\Phi(X_{t_k}(\omega))\stackrel{k\to\infty}{\longrightarrow}\Phi
(X_t(\omega))$ for
all \mbox{$\Phi\in\Pi$}, and $X_{t_k}(\omega)\stackrel{k\to\infty
}{\longrightarrow}
X_t(\omega)$ follows as in Remark~\ref{remsepsep2}. Consequently,
$t\mapsto X_t(\omega)$ is continuous for all $t\in\mathbb Q\cap
[0,T]$ and hence is continuous for all $t\in[0,T]$, because
$\mathcal X$ has sample paths in $\mathcal D_E([0,\infty))$ by
assumption. Since $T$ was arbitrary, continuity of sample paths
$t\mapsto X_t$ follows.\vspace*{10pt}

\textit{Step} 3: proof of (\ref{eq2}). Now, we show that the right-hand
side of (\ref{eq2}) is the conditional quadratic variation of
$\Phi(\mathcal X)$. First note that since $\Omega^{(1)}$ is first
order,
%
%e4.15 #&#
%
\begin{equation}\label{d1}
\Omega\Phi^2 - 2\Phi\cdot\Omega\Phi= \Omega^{(2)} \Phi^2 -
2\Phi\cdot\Omega^{(2)}\Phi.
\end{equation}
We use martingales $(M_\Phi(t))_{t\geq0}$ with
%
%e4.16 #&#
%
\begin{equation}
M_\Phi(t):= \Phi(X_t) - \int_0^t \Omega\Phi(X_s) \,ds.
\end{equation}
Now we decompose the square of the martingale

%e4.17 #&#
%
\begin{equation}
\label{dgpGir18}\quad
(M_\Phi(t))^2 = \Phi^2(X_t) - 2 M_\Phi(t)\cdot\int_0^t\Omega
\Phi(X_s)\,ds - \biggl(\int_0^t\Omega\Phi(X_s)\,ds\biggr)^2.\hspace*{-25pt}
\end{equation}
Next using
%(\ref{dgpGir18}) for $\Phi^2$ and
partial integration we have
%
%e4.18 #&#
%
\begin{eqnarray}
\label{dgpGir18b}\quad
(M_\Phi(t))^2 &=& M_{\Phi^2}(t) + \int_0^t \Omega\Phi^2(X_s)\,ds - 2
\int_0^t
M_\Phi(s) \cdot\Omega\Phi(X_s)\,ds \nonumber\\[-8pt]\\[-8pt]
&&{}- 2 \int_0^t
\Omega\Phi(X_s) \,dM_\Phi(s)\,ds - \biggl(\int_0^t\Omega\Phi
(X_s)\,ds\biggr)^2.
\nonumber
\end{eqnarray}
With (\ref{d1}) we get finally
%
%e4.19 #&#
%
\begin{eqnarray} \label{dgpGir18c}
(M_\Phi(t))^2
&=& \biggl( M_{\Phi^2}(t) - 2\int_0^t \Omega\Phi(X_s)
\,dM_\Phi(s) \,ds
\biggr) \nonumber\\[-8pt]\\[-87pt]
&&{} + \int_0^t \Omega^{(2)} \Phi^2 (X_s) -
2\Phi(X_s)\cdot\Omega^{(2)}\Phi(X_s)\,ds.
\nonumber
\end{eqnarray}
Clearly, this is the decomposition of the
submartingale $M_\Phi^2$ into its martingale part and its predictable part
of finite variation, and (\ref{eq2}) follows.
%See e.g.\ \citet{RevuzYor2010}, Proposition VIII.3.3.
\end{pf*}

%s4.3 #&#
\subsection{Operators for the tree-valued FV process}
We apply the concepts of the last subsection to the different
components of the generator for the TFVMS process.
%
%pr4.10 #&#
%
\begin{proposition}[(Order of generator terms of the TFVMS
process)]\label{PfirstSecond}
\textup{(1)}~The operators $\Omega^{\mathrm{grow}}$, $\Omega^{\mathrm{sel}}$
and $\Omega^{\mathrm{mut}}$ are first-order operators with respect
to $\Pi^1$.

\textup{(2)} The operator $\Omega^{\mathrm{res}}$ is a second-order operator
with respect to $\Pi^0$. Moreover, for $\Phi= \Phi^{n,\phi}\in\Pi^0_n$
and with $\rho_1^n$ from (\ref{eqrho}),
%
%e4.20 #&#
%
\begin{eqnarray}\label{eqresQv}
&&
\Omega^{\mathrm{res}} \Phi^2(\smallu) - 2\Phi(\smallu) \cdot
\Omega^{\mathrm{res}} \Phi(\smallu) \nonumber\\[-8pt]\\[-8pt]
&&\qquad
= \gamma\sum_{k,l=1}^n
\langle\nu^\smallu, \phi\cdot(\phi\circ\rho_1^n)\cdot
\theta_{k,n+l} - \phi\cdot(\phi\circ\rho_1^n)\rangle.
\nonumber
\end{eqnarray}
\end{proposition}
\begin{pf}
Let $\Phi^\phi\in\Pi^1_{n}$. Then, using $\rho_1^n$ from
(\ref{eqrho}), we show that $\Omega^{\mathrm{grow}}, \Omega
^{\mathrm{sel}}$
and $\Omega^{\mathrm{mut}}$ are first-order operators by calculating
\begin{eqnarray*}
\Omega^{\mathrm{grow}} \Phi^2(\smallu) & = & \bigl\langle\nu
^\smallu,
\langle\nabla_{\underline{\underline
r}}\phi\cdot(\phi\circ\rho_1^n), \underline2 \rangle
\bigr\rangle\\[-1pt]
&=& \bigl\langle\nu^\smallu,
\langle\nabla_{\underline{\underline r}}\phi,
\underline{\underline2} \rangle\cdot\phi\circ
\rho_1^n\bigr\rangle+ \bigl\langle\nu^\smallu, \phi\cdot
\langle\nabla_{\underline{\underline r}} (\phi\circ\rho_1^n),
\underline{\underline2}\rangle\bigr\rangle\\[-1pt]
& = &
2\Phi(\smallu)\cdot
\Omega^{\mathrm{grow}} \Phi(\smallu),\\[-1pt]
\Omega^{\mathrm{sel}} \Phi^2(\smallu)
% & = \Omega^{\mathrm{sel}}
% \langle\nu^\smallu, (\phi,\phi)_n\rangle
% \\
& = & \alpha\sum_{k=1}^{2n} \langle\nu^{\smallu}, \phi\cdot
(\phi\circ\rho_1^n)\cdot\chi'_{k,2n+1} -
\phi\cdot(\phi\circ\rho_1^n)\cdot\chi'_{2n+1, 2n+2}\rangle\\[-1pt]
& = & 2\alpha\sum_{k=1}^{n} \langle\nu^{\smallu}, \phi\cdot
\chi'_{k,n+1} \cdot(\phi\circ\rho_1^{n+1}) -
\phi\cdot\chi'_{n+1, n+2} \cdot(\phi\circ\rho_1^{n+2})\rangle\\[-1pt]
& = & 2\alpha\langle\nu^{\smallu}, \phi\rangle\cdot
\sum_{k=1}^{n} \langle
\nu^{\smallu}, \phi\cdot\chi'_{k,n+1} - \phi\cdot\chi
'_{n+1,n+2} \rangle\\[-1pt]
& = & 2\Phi(\smallu) \cdot\Omega^{\mathrm{sel}} \Phi(\smallu),\\[-1pt]
\Omega^{\mathrm{mut}} \Phi^2(\smallu)
& = &\sum_{k=1}^{2n} \bigl\langle
\nu^\smallu, B_k\bigl(\phi\cdot(\phi\circ\rho_1^n)\bigr)\bigr\rangle= 2
\sum_{k=1}^{n} \langle\nu^\smallu, (B_k\phi) \cdot
(\phi\circ\rho_1^n)\rangle\\[-1pt]
& = & 2\langle\nu^\smallu,
\phi\rangle\cdot\sum_{k=1}^{n} \langle\nu^\smallu,
B_k\phi\rangle= 2\Phi(\smallu) \cdot\Omega^{\mathrm{mut}}
\Phi(\smallu).
\end{eqnarray*}
For $\Omega^{\mathrm{res}}$, Corollary 2.15 in
\citet{GrevenPfaffelhuberWinter2011} shows
(\ref{eqresQv}). Informally, the second-order term, as given in
(\ref{eqresQv}), arises by interactions between two samples, drawn
independently from $\smallu$.

In order to establish $\Omega^{\mathrm{res}}$ as a second-order
operator, observe that all interactions between three independently
drawn samples are due to interactions between pairs of samples. A
formal calculation showing that $\Omega^{\mathrm{res}}$ is second order
is as follows:
\begin{eqnarray*}
&&
-3\Phi(\smallu) \Omega^{\mathrm{res}} \Phi^2(\smallu) + 3
\Phi^2(\smallu) \Omega^{\mathrm{res}} \Phi(\smallu) \\
&&\qquad = - 3
\Phi(\smallu) \bigl(\Omega^{\mathrm{res}} \Phi^2(\smallu) - 2
\Phi(\smallu) \Omega^{\mathrm{res}} \Phi(\smallu)\bigr)
- 3
\Phi^2(\smallu) \Omega^{\mathrm{res}} \Phi(\smallu) \\
&&\qquad = - 3
\Phi(\smallu) \gamma\sum_{k,l=1}^n \langle\nu^\smallu,
\phi\cdot(\phi\circ\rho_1^n)\cdot\theta_{k,n+l} - \phi\cdot
(\phi\circ\rho_1^n)\rangle\\
&&\qquad\quad{}- 3 \Phi^2(\smallu) \Omega^{\mathrm{res}}
\Phi(\smallu),
\end{eqnarray*}
where we used (\ref{eqresQv}) in the last step. Furthermore,
\begin{eqnarray*}
\Omega^{\mathrm{res}} \Phi^3(\smallu) &=& \frac\gamma2
\sum_{k,l=1}^{3n} \bigl\langle\nu^\smallu, \bigl(\phi\cdot(\phi
\circ\rho_1^n) \cdot(\phi\circ\rho_1^{2n})\bigr)\circ\theta_{k,l}\\
&&\hspace*{70.5pt}{} - \phi\cdot(\phi\circ\rho_1^n) \cdot
(\phi\circ\rho_1^{2n}) \bigr\rangle\\
&=& \frac{3\gamma}2
\sum_{k,l=1}^{n} \langle\nu^\smallu, (\phi\circ\theta_{k,l})
\cdot(\phi\circ\rho_1^n) \cdot(\phi\circ\rho_1^{2n})\\
&&\hspace*{78.6pt}{} - \phi
\cdot(\phi\circ\rho_1^n) \cdot(\phi\circ\rho_1^{2n}) \rangle
\\
&&{} + \frac{6\gamma}2 \sum_{k,l=1}^{n} \bigl\langle
\nu^\smallu, \bigl(\bigl(\phi\cdot(\phi\circ\rho_1^n)\bigr)\circ
\theta_{k,n+l}\bigr) \cdot(\phi\circ\rho_1^{2n})\\
&&\hspace*{100.4pt}{}  - \bigl(\phi\cdot
(\phi\circ\rho_1^n) \cdot(\phi\circ\rho_1^{2n})\bigr) \bigr\rangle\\
&=& 3 \Phi^2(\smallu) \Omega^{\mathrm{res}} \Phi(\smallu)\\
&&{} + 3
\Phi(\smallu) \gamma\cdot\sum_{k,l=1}^n \bigl\langle\nu^\smallu,
\bigl(\phi\cdot(\phi\circ\rho_1^n)\bigr)\cdot\theta_{k,n+l} - \phi
\cdot(\phi\circ\rho_1^n)\bigr\rangle.
\end{eqnarray*}
Summing the last two displays, we see that $\Omega^{\mathrm{res}}$ is
second order with respect to~$\Pi^0$, according to
Definition~\ref{deffirstO}.
\end{pf}

%s5 #&#
\section{Duality}
\label{sduality}
One of the main tools in studying the long-time behavior of a Markov
process is to construct and to study a dual process $\Xi$ in the
limit $t \to\infty$. In this section, we define a dual process
of the TFVMS process, which takes values in functions. Its state space
is the following separable metric space [recall (\ref{BDF1})]:
%
%e5.1 #&#
%
\begin{equation}\label{eq904}
\Upsilon:=\bigcup_{n=0}^\infty\overline{\mathcal C}{}^1_n,
\end{equation}
and the \textit{duality function} $H(\cdot, \cdot)$ is
%
%e5.2 #&#
%
\begin{equation}\label{eq905}
H\dvtx \cases{ \M_I \times\Upsilon\to\mathbb R, \cr
(\smallu, \xi) \mapsto H(\smallu, \xi):= \langle\nu^\smallu,
\xi
\rangle.}
\end{equation}
We next define the Markov process $\Xi$. The formal duality
result is given in Proposition~\ref{propdual}.
%
%de5.1 #&#
%
\begin{definition}[(The function-valued dual process
$\Xi$)] \label{defdual} %
The process $\Xi=(\Xi_t)_{t\geq0}$ is a
piecewise deterministic jump process with state space
$\Upsilon$. Recall that the mutation transition kernel has the form
(\ref{eq1001}) for some $z\in[0,1]$. Here are the evolution rules:

\begin{longlist}[(2)]
\item[(1)] Between jumps the process evolves according to the semigroup
%
%e5.3 #&#
%
\begin{equation}
\label{eqdua0}
(S_t \xi)(\underline{\underline r}, \underline u) = \xi(s_t
\underline{\underline r}, \underline u)
\end{equation}
with
%
%e5.4 #&#
%
\begin{equation}
( s_t (r_{ij}))_{1\leq i<j}:= (r_{ij}+2t)_{1\leq i<j}.
\end{equation}
\item[(2)] To describe the resampling transition, we define
%
%e5.5 #&#
%
\begin{equation}
\label{eqsig1}
(\overline\sigma_l((\underline{\underline r}, \underline u))) =
\bigl(\bigl(r_{i-\ind{i> l}, j-\ind{j> l}}, u_{i-\ind{i> l}}\bigr)\bigr).
\end{equation}
Then for $n\geq1$, the process jumps from the state
$\xi\in\overline{\mathcal C}{}^1_n(\mathbb R_+^{{\mathbb
N\choose2}}\times I^{\mathbb N})$ to
%
%e5.6 #&#
%e5.7 #&#
%e5.8 #&#
%
\begin{eqnarray}
\label{eqdua1}
&\displaystyle \Theta_{kl}\xi:= \xi\circ
\theta_{kl}\circ\overline\sigma_{l}
\qquad\mbox{at rate }\frac\gamma2,
k,l=1,\ldots,n,&\\
\label{eqdua2}
&\displaystyle \widetilde\beta_k\xi\qquad\mbox{at rate }\vartheta(1-z),
k=1,\ldots,n,&\\
\label{eqdua3}
&\displaystyle \widebar\beta_k\xi\circ\overline\sigma_{k} \qquad\mbox{at rate
}\vartheta z, k=1,\ldots,n,&
\end{eqnarray}
with\vspace*{1pt} $\theta_{kl}$ from before (\ref{pp11b}),
$\widetilde\beta_k\xi$ and $\widebar\beta_k\xi$ as in
(\ref{eq1001}). Since $\overline\beta_k\xi$ does not depend on
the $k$th variable, we note that
$\langle\nu^\smallu,\overline\beta_k\xi\circ\overline\sigma
_k\rangle=
\langle\nu^\smallu,\overline\beta_k\xi\rangle$ for
$\smallu\in\mathbb U$; see also (\ref{eq1001}), (\ref{eq33}) and
Remark~\ref{remdualbeh} [item (3)].
\item[(3)]For haploid and diploid selection, (\ref{eqpp15b}) and
(\ref{dgpGir4}), respectively, we use an operation
%
%e5.9 #&#
%
\begin{equation}
\label{eqsig2}
(\sigma_k((\underline{\underline r}, \underline u))) =
\bigl(\bigl(r_{i+\ind{i\geq k}, j+\ind{j\geq k}}, u_{i+\ind{i\geq k}}\bigr)\bigr),
\end{equation}
which arises\vspace*{1pt} by deleting the $k$th column and line from
$\underline{\underline r}$ and the $k$th entry from $\underline
u$. Then we introduce jumps from $\xi$ to (in the haploid and
diploid case, resp.)
%
%e5.10 #&#
%e5.11 #&#
%
\begin{eqnarray}
\label{eqdua3a}
&\xi\cdot\chi_{k} + (\xi\circ\sigma_k)\cdot(1-\chi_{k})
\qquad\mbox{at rate }\alpha, k=1,\ldots,n,&\\
\label{eqdua3b}
&\xi\cdot\chi'_{k,n+2} +
(\xi\circ\sigma_k)\cdot(1-\chi'_{k,n+2}) \qquad\mbox{at rate }\alpha,
k=1,\ldots,n,&
\end{eqnarray}
with $\chi_k$ as in (\ref{eqpp15c}), $\chi'_{k,n+2}$ as in
(\ref{dgpGir5}). (These transitions are reminiscent of the dual
process $(\eta_t, \mathcal G_t^{++})_{t\geq0}$ from
\citet{DGsel0}. In particular, they differ from the construction
given in \citet{DawsonGreven1999b}. See Remark~\ref{remdualbeh}
[item (2)] for the advantage of our construction.)
\item[(4)] If $\xi\in\overline{\mathcal C}{}^1_0$ is constant, it
stays in
$\xi$ for all times.
\end{longlist}
\end{definition}

%
%re5.2 #&#
%
\begin{remark}[(Behavior of $\Xi$ and underlying birth and death
process)]
\label{remdualbeh}
(1)~To better understand what is going on, look at the form of the
function after the transition. For example,
for (\ref{eqdua1}),
%
%e5.12 #&#
%
\begin{eqnarray}
\label{eq601}\qquad
(\Theta_{kl}\xi)(\underline{\underline r}, \underline u) & = &
\xi(\theta_{kl}(r_{ij})_{i,j=1,2,\ldots,l-1,l,l,l+1,\ldots},
(u_i)_{i=1,\ldots,l-1,l,l,l+1,\ldots}))\nonumber\\[-8pt]\\[-8pt]
& = &
\xi((r_{ij})_{i,j=1,2,\ldots,l-1,k,l,l+1,\ldots},
(u_i)_{i=1,\ldots,l-1,k,l,l+1,\ldots})).
\nonumber\\[-16pt]\nonumber
\end{eqnarray}

\begin{longlist}[(2)]
\item[(2)] In order to show that $\Xi$ is dual to the TFVMS
(Proposition~\ref{propdual}), we could as well have used a
transition from $\xi$ to $\xi\circ\theta_{kl}$ instead of
(\ref{eqdua1}), to $\xi\cdot\chi_k + \xi\cdot(1-\chi_{n+1})$ and
to $\xi\cdot(\chi'_{k,n+1} + (1-\chi'_{n+1,n+2}))$ instead of
(\ref{eqdua3a}) and (\ref{eqdua3b}), respectively. However, the
above formulation has two advantages:
\begin{itemize}[$\blacktriangleright$]
\item[$\blacktriangleright$] By (\ref{eq601}), we see that $\Theta
_{kl}\xi\in
\overline{\mathcal C}{}^1_{n-1}$ for $\xi\in\overline{\mathcal
C}{}^1_{n}$.
\item[$\blacktriangleright$] We can show that $t\mapsto\|{\Xi
_t}\|_{\infty}$ is
nonincreasing (see Proposition~\ref{lnoninc}).
\end{itemize}
\item[(3)] For the process $\Xi$, consider the process $(N_t)_{t\geq0}$,
where $N_t=n$ if $\Xi_t \in\widebar{\mathcal C}{}^1_n$. In the case
of selection acting on haploids, the process jumps from $n$ to
%
%e5.13 #&#
%
\begin{eqnarray}
\label{eq200pp}
&n-1 \qquad\mbox{at rate } \gamma\pmatrix{n\cr2} + \vartheta z \cdot
n,&\nonumber\\[-9pt]\\[-9pt]
&n+1 \qquad\mbox{at rate } \alpha n.&
\nonumber
\end{eqnarray}
Note that the additional rate $\vartheta z \cdot n$ of decrease
comes from the choice of transitions
$\xi\to\overline\beta_k\xi\circ\overline\sigma_k$ instead of
$\xi\to\overline\beta_k\xi$. The process $(N_t)_{t\geq0}$ plays
(for $z=0$) again an important role in Section~\ref{sschar} in
estimating the numbers of ancestors of the total population.\vspace*{-2pt}
% In the diploid case we obtain upon birth two new
% variables.
\end{longlist}
\end{remark}

We can now state the duality relation between $\mathcal U$ and $\Xi$.\vspace*{-2pt}
%
%pr5.3 #&#
%
\begin{proposition}[(Duality relation)]\label{propdual} %
Let $\mathcal U=(\mathcal U_t)_{t\geq0}$ be the
tree-valued Fleming--Viot process and $\Xi=(\Xi_t)_{t\geq0}$ the
function-valued process from Definition~\ref{defdual}.\vspace*{-2pt}
\begin{longlist}[(2)]
\item[(1)] The set of functions $\{\smallu\mapsto H(\smallu, \xi)\dvtx
\xi\in\Upsilon\}$ from (\ref{eq905}) is separating on~$\mathbb
M^I$.
\item[(2)] The processes $\mathcal U$, started in $\mathcal U=\smallu$,
and $\Xi$, started in $\Xi_0=\xi$, are dual to each other, that
is, for $H$ from (\ref{eq905}) and $t\geq0$,
%
%e5.14 #&#
%
\begin{equation}\label{eqdual}
\mathbf E_\smallu[H(\mathcal U_t, \xi)] = \mathbf
E_\xi[H(\smallu, \Xi_t)].\vspace*{-2pt}
\end{equation}
\end{longlist}
\end{proposition}
\begin{pf}
For (1) we just note that $\{\smallu\mapsto\langle\nu^\smallu,
\xi\rangle\dvtx \xi\in\Upsilon\} = \Pi^1$ which is separating by
Proposition 4.1 in \citet{DGP-topo2011}. For~(2) we have to show
that [\citet{EthierKurtz86}, Proposition~4.4.7]
%
%e5.15 #&#
%
\begin{equation}\label{eq517}
(\Omega\langle\cdot, \xi\rangle) (\nu^\smallu) =
(\Omega_{\mathrm{dual}} \langle\nu^\smallu, \cdot\rangle)
(\xi),\qquad \smallu\in\mathbb U^I, \xi\in\Upsilon,
\end{equation}
where $\Omega$ is the generator of $\mathcal U$, and
$\Omega_{\mathrm{dual}}$ is the generator of the dual process~$\Xi$. We begin by calculating the left-hand side. For $\xi\in
\overline{\mathcal C}{}^1_n$, in the case of diploid selection (here
the operators act on the first argument), we obtain
%
%e5.16 #&#
%
\begin{eqnarray}
\label{eqdualproof1}\qquad
\Omega^{\mathrm{grow}} \langle\nu^\smallu, \xi\rangle
& = &
\langle\nu^\smallu, \langle\nabla_{\underline{\underline
r}} \xi, \underline{\underline2}\rangle\rangle,\nonumber\\[-2pt]
\Omega^{\mathrm{res}} \langle\nu^\smallu, \xi\rangle
& = &
\frac\gamma2\sum_{k,l=1}^n \langle\nu^\smallu, \xi\circ
\theta_{k,l} - \xi\rangle= \frac\gamma2\sum_{k,l=1}^n \langle
\nu^\smallu, \xi\circ\theta_{kl}\circ\overline\sigma_l - \xi
\rangle,\nonumber\\[-2pt]
\Omega^{\mathrm{mut}} \langle\nu^\smallu, \xi\rangle
& = &
\vartheta z \sum_{k=1}^n \langle\nu^\smallu, \widebar
\beta_k\xi\circ\overline{\sigma}_k - \xi\rangle+ \vartheta
(1-z) \sum_{k=1}^n\langle\nu^\smallu, \widetilde
\beta_k\xi- \xi\rangle,\\[-2pt]
\Omega^{\mathrm{sel}} \langle\nu^\smallu, \xi\rangle
& = &
\alpha\sum_{k=1}^n \langle\nu^\smallu, \xi\cdot\chi'_{k,n+1} -
\xi\cdot\chi'_{n+1,n+2}\rangle\nonumber\\[-2pt]
& = &\alpha\sum_{k=1}^n
\langle\nu^\smallu, \xi\cdot\chi'_{k,n+2} -
(\xi\circ\sigma_k)\cdot\chi'_{k,n+2}\rangle
\nonumber
\end{eqnarray}
due to the exchangeability of $\nu^\smallu$, where we have used that
$ \langle\nu^\smallu, \widebar\beta_k\xi\rangle=
\langle\nu^\smallu,\widebar\beta_k\xi\circ\overline\sigma
_k\rangle$
in $\Omega^{\mathrm{mut}}$. Summing both sides of all terms in the
last display exactly gives the left-hand side of (\ref{eq517}). An
analogous calculation shows this in case of haploid selection.

Next we calculate the righ-hand side of (\ref{eq517}). The generator of
the Markov process $\Xi$ is easy to write down for functions of the
form $\widebar{\mathcal C}{}^1 \ni\xi\mapsto\langle\nu, \xi\rangle$
and $\nu\in\mathcal M_1(\mathbb R_+^{{\mathbb N\choose2}}\times
I^{\mathbb N})$. Let $\xi\in\widebar{\mathcal C}{}^1_n$ for some
$n=0,1,2,\ldots.$

First, consider the semigroup $(S_t)_{t\geq0}$. Its generator is
given by
%
%e5.17 #&#
%
\begin{equation}\label{eq604}
\langle\nu, \xi\rangle\mapsto\langle\nu,
\langle\nabla_{\underline{\underline r}}\xi, \underline{\underline
2}\rangle\rangle.
\end{equation}
The other parts of the dynamics of $\Xi$ are pure jump. Hence, the
generator of $\Xi$ acts on the above functions in the following way:
\begin{eqnarray*}
% \label{eq502a}
\Omega_{\mathrm{dual}} \langle\nu, \xi\rangle
& = &\langle\nu,
\langle\nabla_{\underline{\underline r}}\xi,
\underline{\underline2}\rangle\rangle+
\frac{\gamma}{2}\mathop{\sum_{{k,l=1}}}_{{k\neq l}}^n (\langle
\nu, \Theta_{kl}\xi\rangle- \langle\nu,\xi\rangle) \\
&&
{} + \vartheta z\sum_{k=1}^n (\langle\nu,
\widebar\beta_k\xi\circ\overline\sigma_k\rangle- \langle\nu,
\xi\rangle) + \vartheta(1-z)\sum_{k=1}^n (\langle\nu,
\widetilde\beta_k\xi\rangle- \langle\nu, \xi\rangle) \\
&&
{} + \alpha\sum_{k=1}^n
\bigl(\langle\nu,\xi\cdot\chi'_{k,n+2} +
(\xi\circ\sigma_k)\cdot(1 - \chi'_{k,n+2})\rangle-
\langle\nu,\xi\rangle\bigr)
\end{eqnarray*}
in the case of diploid selection. An analogous expression holds for
haploid selection. Combining the last display with
(\ref{eqdualproof1}) gives (\ref{eq517}).
\end{pf}

The following is fundamental in using the dual process for the analysis
of the
long-time behavior of $\mathcal U$.
%
%pr5.4 #&#
%
\begin{proposition}[(Long-time behavior of $\Xi$)]
\label{lnoninc}
Let $\Xi=(\Xi_t)_{t\geq0}$ be the dual process
from Definition~\ref{defdual}. Then, the following assertions hold:
\begin{longlist}[(2)]
\item[(1)] $t\mapsto\|{\Xi_t}\|_{\infty}$ is a.s. nonincreasing;
\item[(2)] if $z\in(0,1]$, then $\Xi_t$ converges to a random variable
$\Xi_{\infty}$ which is a.s. bounded by $\|{\Xi_0}\|_{\infty}$;
\item[(3)] there is an a.s. finite time $T>0$ such that $\Xi_T$ does not
depend on $\underline{\underline r}$.
\end{longlist}
\end{proposition}
\begin{pf}
(1) By a restart argument and right-continuity of $(\Xi_t)_{t\geq
0}$, it suffices to show that $\|{\Xi_t}\|_{\infty}\leq
\|{\Xi_0}\|_{\infty}$, almost surely. For this, we consider all
transitions of the dual process. Between jumps it evolves according
to the semigroup $(S_t)_{t\geq0}$ and, given $\Xi_0 = \xi$,
%
%e5.18 #&#
%
\begin{equation}\label{eq503b}
\|{ S_t\xi}\|_{\infty}= \sup_{(\underline{\underline r},
\underline u)} \bigl|\xi\bigl((r_{ij}+2t)_{1\leq i<j},
\underline
u\bigr)\bigr|
\leq\|{\xi}\|_{\infty}.
\end{equation}
If $\Xi_{t-}=\xi$ and a jump occurs at time $t$, we have one of the
following cases:
%
%e5.19 #&#
%e5.20 #&#
%
\begin{eqnarray}\label{eq504}
&&\hspace*{-12pt}
\|{\Xi_t}\|_{\infty} = \|{\Theta_{kl}\xi}\|_{\infty}=
\|{\xi\circ\theta_{kl}\circ
\overline\sigma_l}\|_{\infty}\leq\|{\xi}\|_{\infty},\nonumber\\
&&\hspace*{-17pt}\|{\beta_k \xi}\|_{\infty}
= \sup_{(\underline{\underline r},
\underline u)} \biggl| \int\xi(\underline{\underline r},
\underline u_k^v)\beta_k(\underline u, dv)\biggr| \leq
\|{\xi}\|_{\infty},\nonumber\\
&&
\|{\xi\cdot\chi_k + (\xi\circ
\sigma_k)\cdot(1-\chi_k)}\|_{\infty}\\
&&\qquad\leq\|{\xi}\|_{\infty}\cdot\|{\chi_k + (1-\chi_k)}\|_{\infty}=
\|{\xi}\|_{\infty},\nonumber\\
&&\Vert\xi\cdot\chi'_{k,n+2}  +
(\xi\circ\sigma_k)\cdot(1-\chi'_{k,n+2})\Vert_{\infty}\nonumber\\
&&\qquad\leq\|{\xi}\|_{\infty}\cdot\|{\chi'_{k,n+2} +
(1-\chi'_{k,n+2})}\|_{\infty}
= \|{\xi}\|_{\infty}.\nonumber
\end{eqnarray}
Hence, all transitions of $\Xi$ do not increase $\|{\Xi_\bullet}\|$,
and the result follows.

(2) Considering all possible transitions, it is clear that for
$\xi\in\overline{\mathcal C}{}^1_n$ (see also
Remark~\ref{remdualbeh}),
%
%e5.21 #&#
%
\begin{eqnarray}
\label{eq602}
&\displaystyle (S_t\xi)\in\overline{\mathcal C}{}^1_n, \qquad(\Theta_{kl}\xi) \in
\overline{\mathcal C}{}^1_{n-1},\qquad \beta_k\xi\in\overline
{\mathcal
C}_n^1,&\nonumber\\
&\displaystyle \xi\cdot\chi_k +
(\xi\circ\sigma_k)\cdot(1-\chi_k)\in\overline{\mathcal
C}{}^1_{n+1},&\\
&\displaystyle \xi\cdot\chi'_{k,n+2} +
(\xi\circ\sigma_k)\cdot(1-\chi'_{k,n+2}) \in\overline{\mathcal
C}{}^1_{n+2}.&\nonumber
\end{eqnarray}
Moreover, in the case $z>0$ and $\xi\in\overline{\mathcal C}_1$, we
have $\widebar\beta\xi\in\overline{\mathcal C}_0$. Recall from
Remark~\ref{remdualbeh} [item (3)] that the process $(N_t)_{t\geq0}$
with $N_t=n$ if $\Xi_t\in\overline{\mathcal C}{}^1_n$ decreases at a
quadratic rate and increases at a linear rate. In particular, there
is an almost surely finite stopping time $T$ with
$\Xi_T\in\overline{\mathcal C}_0$; that is, $\Xi_T$ is constant with
$|\Xi_T| \leq\|{\Xi_0}\|_{\infty}$; see~(1).\vspace*{1pt}

(3) Note that any $\xi\in\widebar{\mathcal C}{}^1_1$ does not
depend on $\underline{\underline r}$. As in (2), $T = \inf\{t\geq
0\dvtx\Xi_t \in\widebar{\mathcal C}{}^1_1\}$ is almost surely finite,
and we are done.
\end{pf}

%s6 #&#
\section[The tree-valued Moran model]{The tree-valued Moran model with
mutation and selection}
\label{secmart-probl-fixed}
In this section, we study the tree-valued process introduced in
Section~\ref{sstmmms}. In Section~\ref{ssgenMM}, we give the
generator of the TMMMS from Definition~\ref{defTMMMS}, show
convergence to the generator of TFVMS in Section~\ref{secconvMMFV}
and obtain some characteristics of the TMMMS in Section~\ref{sschar}.

%s6.1 #&#
\subsection{The martingale problem for the TMMMS}
\label{ssgenMM}
Recall the TMMMS $\mathcal U^N = (\mathcal U_t^N)_{t\geq0}$ with
$\mathcal U_t^N:=\overline{(U_N, r^N_t, \mu^N_t)}$ from
Definition~\ref{defTMMMS2}. Its state space is
%
%e6.1 #&#
%
\begin{equation}
\mathbb U_N^I:=\mathbb M_N^I\cap\mathbb U^I,\qquad \mathbb
M_N^I:=\{\overline{(X,r,\mu)}\in\mathbb M^I\dvtx N\mu\in\mathcal
N(X\times I)\},
\end{equation}
where $\mathcal N(X\times I)$ is the set of counting measures on
$X\times I$. Note that $\mathbb U_N^I$ is Polish as a closed
subspace of the Polish space $\mathbb U^I$.

In order to construct the TMMMS via its generator, we need to define
its domain. The construction we use here is similar to the approach
taken in Sections~\ref{ssstatespace} and~\ref{ssmp}, the main
difference being that we have to sample individuals from finite
populations without replacement. Compare analogous concepts from
Definition~\ref{defmdistmat}.
%
%de6.1 #&#
%
\begin{definition}[(Finite marked distance matrix distribution)]
Let $\smallx= \overline{(X,r,\mu)}\in\mathbb M^I_N$.
\begin{longlist}[(2)]
\item[(1)] The sampling without replacement from $\mu$ uses the measure
%
%e6.2 #&#
%
\begin{eqnarray}
\label{eq801}
\mu^{\otimes\downarrow N}(d\underline x, d\underline u) :\!&=&
\mu(dx_1, du_1) \cdot\frac{\mu-(1/N)
\delta_{x_1,u_1}}{1-1/N}(dx_2, du_2) \nonumber\\
&&\times{} \cdots
\frac{\mu- (1/N) \sum_{i=1}^{N-1}\delta_{x_i,
u_i}}{1-({N-1})/{N}}(dx_N, du_N)\\
&\in&\mathcal M_1 (X^N
\times I^N)
\nonumber
\end{eqnarray}
for $(\underline{x}, \underline u) \in X^N\times I^N$.
\item[(2)] We define
%
%e6.3 #&#
%
\begin{equation}
\label{eqN-distmat}
R^{N,(X,r)}\dvtx
\cases{
(X \times I)^N \to\mathbb R_+^{{ N\choose2}} \times I^N, \cr
((x_i,u_i)_{1\le i \le N}) \mapsto(
(r(x_i,x_j))_{1\le i < j \le N}, (u_k)_{1 \le k \le N}),}
\end{equation}
and let $\nu^{N,\smallx}$ denote the corresponding marked distance
matrix distribution
%
%e6.4 #&#
%
\begin{equation}
\label{eqN-mdmdistr}
\nu^{N,\smallx}:=\bigl(R^{N,(X,r)}\bigr)_\ast\mu^{\otimes\downarrow N} \in
\mathcal M_1 \bigl(\mathbb R_+^{{ N\choose2}} \times I^N\bigr).
\end{equation}
\end{longlist}
\end{definition}
%
%re6.2 #&#
%
\begin{remark}[(Marked distance distribution is well defined on
$\mathbb U^I$)]%
(1) As in Remark~\ref{remdistwell}, for $\smallx=
\overline{(X,r,\mu)} \in\mathbb M_N^I$, the marked distance
matrix distribution $\nu^{N,\smallx}$ does not depend on the
representative $(X,r,\mu)$ and hence is well defined.

(2) Let $\smallx= \overline{(X,r,\mu)} \in\mathbb M^I
\setminus
\mathbb M^I_N$. Then, $\mu^{\otimes\downarrow N}$ can still be
defined as in (\ref{eq801}), but is a signed measure. The same
holds for $\nu^{N,\smallx}$.
\end{remark}

Now we can define the domain and range of the generator of the TMMMS.

%de6.3 #&#
%
\begin{definition}[(Polynomials on $\mathbb U_N^I$)]%
A function $\Phi\dvtx \mathbb U_N^I\to\mathbb R$ is a \textit{polynomial}
if there exists $\phi\in\mathcal B(\mathbb R_+^{{ N\choose
2}}\times I^N)$ such that
%
%e6.5 #&#
%
\begin{equation}
\label{eqN-polyn}
\Phi_N^\phi(\smallu) = \langle\nu^{N,\smallu}, \phi\rangle=
\int_{\mathbb R_+^{{ N\choose2}} \times I^N}
\phi(\underline{\underline r}, \underline u )
\nu^{N,\smallu}(d\underline{\underline r}, d \underline u).
\end{equation}
In this case, we set $\Phi_N^\phi:=\Phi$. As the space of all
polynomials of this form is not an algebra, we define
%
%e6.6 #&#
%e6.7 #&#
%
\begin{eqnarray}
\label{eq7}
\Pi_N &:=& \mbox{algebra generated by } \bigl\{\Phi_N^\phi\dvtx \phi\in
\mathcal
B\bigl(\mathbb R_+^{{ N\choose2}} \times I^N\bigr)\bigr\}, \\
\label{eq6}
\Pi_N^1 &:=& \mbox{algebra generated by } \bigl\{\Phi_N^\phi\dvtx \phi
\in\mathcal C_b^1\bigl(\mathbb R_+^{{ N\choose2}} \times I^N\bigr)\bigr\},
\end{eqnarray}
where differentiability in $\mathcal C_b^1(\mathbb R_+^{{ N\choose2}}
\times I^N)$ is only required for the coordinates in $\mathcal
C_b^1(\mathbb R_+^{{ N\choose2}})$.
\end{definition}

For the definition of the generator of the TMMMS recall the notation
introduced in Definition~\ref{defgentvms} and (\ref{agx1e}).
%
%de6.4 #&#
%
\begin{definition}[(Generator of the TMMMS)]%
The generator of the TMMMS with population size $N$ is the linear
operator $\Omega^N$ on $\Pi_N$ with domain $\Pi_N^1$ given by
%
%e6.8 #&#
%
\begin{equation}
\label{eq802}
\Omega^N:=\Omega^{\mathrm{grow},N} +
\Omega^{\mathrm{res},N}+\Omega^{\mathrm{mut},N}+\Omega^{\mathrm{sel},N}.
\end{equation}
The growth and resampling operators are given by
%
%e6.9 #&#
%e6.10 #&#
%
\begin{eqnarray}
\label{eq8}
\Omega^{\mathrm{grow},N} \Phi_N^\phi(\smallu) &:=& \langle
\nu^{N,\smallu},\langle\nabla_{\underline{\underline r}} \phi,
\underline{\underline2}\rangle
\rangle, \\
\label{eq9}
\Omega^{\mathrm{res},N} \Phi_N^\phi(\smallu) &:=& \frac\gamma
2\sum_{k,l=1}^N (\langle\nu^{N,\smallu}, \phi\circ
\theta_{k,l}\rangle- \langle\nu^{N,\smallu}, \phi\rangle).
\end{eqnarray}
The mutation operator is given by
%
%e6.11 #&#
%
\begin{equation}
\label{eq10}
\Omega^{\mathrm{mut},N} \Phi_N^\phi(\smallu):= \vartheta
\sum_{k=1}^N \langle\nu^{N,\smallu}, B_k \phi\rangle.
\end{equation}
The selection operators in the cases of haploid and diploid
selection are given by
%
%e6.12 #&#
%e6.13 #&#
%
\begin{equation}
\label{eq11}
\Omega^{\mathrm{sel},N} \Phi_N^\phi(\smallu) := \frac\alpha N
\sum_{k,l=1}^N \langle\nu^{N,\smallu}, \chi_k(\phi\circ
\theta_{k,l} - \phi) \rangle
\end{equation}
and
\begin{equation}
\label{eq12}
\Omega^{\mathrm{sel},N} \Phi_N^\phi(\smallu) := \frac\alpha{N^2}
\sum_{k,l,m=1}^N \langle\nu^{N,\smallu}, \chi_{k,m}'(\phi\circ
\theta_{k,l} - \phi) \rangle,
\end{equation}
respectively.
\end{definition}
%
%re6.5 #&#
%
\begin{remark}[(Interpretation of generator terms)]\label{remgentermsN}
Clearly, the generator terms
$\Omega^{\mathrm{grow},N}$ and $\Omega^{\mathrm{res},N}$ describe tree
growth and resampling; see also Section 5.1 of
\citet{GrevenPfaffelhuberWinter2011} for the case without marks. The
terms $\Omega^{\mathrm{res},N}$ and $\Omega^{\mathrm{mut},N}$ describe
resampling and mutation arising from the Poisson processes
$\eta_{\mathrm{res}}$ and $\eta_{\mathrm{mut}}$ from
Definition~\ref{defgraph}, respectively. For selection, recall
$\eta_{\mathrm{sel}}$ from that definition. In the case of haploid
selection, $l$ is replaced by an offspring of $k$ at rate $\alpha
\chi(u_k) /N$, for $k,l=1,\ldots,N$, which easily translates into
(\ref{eq11}). The case of diploid selection is similar.
\end{remark}
%
%pr6.6 #&#
%
\begin{proposition}[(Well-posedness of TMMMS martingale problem)]%
\label{proptvmd}
Let $N \in\N$, $\mathbf P_0^N \in\mathcal M_1 (\mathbb U_N^I)$,
$\Pi_N^1$ as in (\ref{eq6}) and $\Omega^N$ as in
(\ref{eq802}). Then, the $(\mathbf P_0^N, \Omega^N$, $
\Pi_N^1)$-martingale problem has exactly one solution, the
tree-valued Moran model with mutation and selection.
\end{proposition}
\begin{pf}
Existence is straight-forward from the graphical
construction (see Definition~\ref{defgraph} and
Remark~\ref{remgentermsN}). In particular, the TMMMS solves the
$(\mathbf P_0^N, \Omega^N$, $\Pi_N^1)$-martingale problem. To get
well-posedness,
note that the $(\mathbf P_0^N,\Omega^{\mathrm{grow},N},
\Pi^1_N)$-martingale problem is well posed. Furthermore
$B:=\Omega^{\mathrm{res},N}+\Omega^{\mathrm{mut},N}+\Omega
^{\mathrm{sel},N}$ is
a bounded jump operator (since the population is finite). Hence,
uniqueness follows from Theorem 4.10.3 in \citet{EthierKurtz86}.
\end{pf}

%s6.2 #&#
\subsection{Convergence of generators}
\label{secconvMMFV}
Here, we prove that the sequence of generators $\Omega^N$ of the TMMMS
defined in (\ref{eq802}) converges (uniformly) to the generator
$\Omega$ for
the TFVMS from (\ref{dgpGir2}).
%
%pr6.7 #&#
%
\begin{proposition}[(Generator convergence)] \label{propconv2TVFV}%
For any $\Phi\in\Pi^1$ there is a sequence
$(\Phi_N)_{N \in\N}$ such that $\Phi_N \in\Pi^1_N$ for all $N$ and
%
%e6.14 #&#
%e6.15 #&#
%
\begin{eqnarray}
\label{eq13}
\lim_{N \to\infty} \sup_{\smallu\in\mathbb U^I} |
\Phi_N(\smallu)
- \Phi(\smallu)|&=&0, \\
\label{eq14}
\lim_{N \to\infty} \sup_{\smallu\in\mathbb U^I} |\Omega^N
\Phi_N(\smallu) - \Omega\Phi(\smallu)|&=&0.
\end{eqnarray}
\end{proposition}
\begin{pf}
Let $\Phi\in\Pi^1$. Then, by definition of $\Pi^1$, $\Phi=
\Phi^{n,\phi}$ for some $n \in\N$ and $\phi\in\overline{\mathcal
C}{}^1$. We define $\widetilde\nu^{N, \smallu}:=(\iota_N)_\ast
\nu^{N,\smallu}$ for
%
%e6.16 #&#
%
\begin{equation}
\label{eq15}\qquad
\iota_N\dvtx
\cases{
\mathbb R^{{ N\choose2}}\times I^N \to\mathbb R^{{ \N\choose
2}}\times
I^\N,\cr
((r_{i,j})_{1 \le i,j\le N},(u_\ell)_{1 \le\ell\le N})
\mapsto((r_{i \simeq N,j \simeq N})_{1 \le i < j},(u_{\ell
\simeq
N})_{1 \le\ell}),}
\end{equation}
where $i \simeq N:=1+((i-1) \mbox{ mod } N)$. We define
$\Phi_N\in\Pi^1$ by setting
%
%e6.17 #&#
%
\begin{equation}\label{eq805}
\Phi_N(\smallu)=\langle\nu^{N,\smallu}, \phi\circ\iota_N
\rangle
= \langle\widetilde{\nu}^{N,\smallu}, \phi\rangle.
\end{equation}
Then there is a constant $C=C(n,\phi) >0$, such that for all $N \ge n$,
%
%e6.18 #&#
%
\begin{equation}
\label{eq16}
\sup_{\smallu\in\mathbb U^I} |
\Phi_N(\smallu)-\Phi(\smallu) |= \sup_{\smallu\in
\mathbb U^I} |\langle\widetilde\nu^{N,\smallu} -
\nu^\smallu, \phi\rangle|\le\frac{C}N.
% \norm{\phi}_{\infty}.
\end{equation}
To show (\ref{eq14}) for $\Phi\in\Pi^1$ in the case $\alpha=0$,
note that $\Omega_{0} \Phi(\smallu)= \langle\nu^\smallu, \psi
\rangle$ and $\Omega^N_0 \Phi_N (\smallu)=\langle\widetilde
\nu^{N,\smallu}, \psi\rangle$ for some $\psi\in\widebar{\mathcal
C}{}^1_n$. Thus, in that case, (\ref{eq14}) follows from
(\ref{eq16}).

It remains to prove the convergence of the selection operators in
haploid and diploid selection cases. We give the proof in the
haploid case;\vadjust{\goodbreak} the diploid case is similar. For $N \ge n$, we have
%
%e6.19 #&#
%
\begin{eqnarray}
\label{eq803}
\Omega^{\mathrm{sel},N} \Phi_N(\smallu) & = & \frac\alpha N
\sum_{k,l=1}^n \langle
\widetilde\nu^{N,\smallu}, \chi_k(\phi\circ\theta_{k,l} - \phi
) \rangle\nonumber\\
&&{} + \frac\alpha N \sum_{k=1}^N \sum_{l=n+1}^N \langle
\widetilde\nu^{N,\smallu}, \chi_k(\phi\circ\theta_{k,l} - \phi)
\rangle\\
&&{} + \frac\alpha N \sum_{k=n+1}^N \sum_{l=1}^n \langle
\widetilde\nu^{N,\smallu}, \chi_k(\phi\circ\theta_{k,l} -
\phi) \rangle.\nonumber
\end{eqnarray}
Here the first summand on the right-hand side is of order $N^{-1}$,
and the second vanishes. Thus, we need to consider only the last
summand. Define the swapping operator $\tau_{k,l}$ through the
permutation $\sigma_{k,l}:= (1,\ldots,k-1,l,k+1,\ldots
,l-1,k,l+1,\ldots,n)$
by $\tau_{k,l} (\underline{\underline r}, \underline u):=
R_{\sigma_{k,l}}^N$ [with an obvious extension of the operator
$R_\sigma$ from (\ref{eq29}) to finite $N$]. Observe that for
$k>n$, and $l\leq n$ by exchangeability of $\nu^{N,\smallu}$, since
$\phi$ only depends on the first $n$ indices,
%
%e6.20 #&#
%
\begin{eqnarray}
\label{eq804}
\langle\widetilde\nu^{N, \smallu},
\chi_k(\phi\circ\theta_{k,l}) \rangle& = &\langle\widetilde
\nu^{N, \smallu}, (\chi_l\cdot\phi)\circ\theta_{k,l}
\rangle\nonumber\\[-8pt]\\[-8pt]
&=& \langle\widetilde\nu^{N, \smallu},
(\chi_l\cdot\phi)\circ\tau_{k,l}\rangle
= \langle\widetilde
\nu^{N, \smallu}, \chi_l\cdot\phi\rangle.
\nonumber
\end{eqnarray}
Hence, for constants $C=C(n,\alpha,\chi,\phi)$ not depending on
$\smallu$ and possibly changing from line to line, by
exchangeability of $\nu^{N,\smallu}$ and (\ref{eq803}),
%
%e6.21 #&#
%
\begin{eqnarray}
\label{eq620}
&&
| \Omega^{\mathrm{sel},N} \Phi_N(\smallu) - \Omega^{\mathrm{sel}}
\Phi(\smallu) | \nonumber\\
&&\qquad \leq\Biggl| \frac{\alpha
(N-n)} N
\sum_{k=1}^n \langle\widetilde\nu^{N,\smallu}, \chi_k \phi-
\chi_{n+1} \phi \rangle\nonumber\\[-8pt]\\[-8pt]
&&\hspace*{71pt}{} - \alpha\sum_{k=1}^n \langle
\nu^{\smallu}, \chi_k \phi- \chi_{n+1} \phi \rangle\Biggr| +
\frac{C}{N} \nonumber\\
&&\qquad \leq\frac{C}{N}\nonumber
\end{eqnarray}
by the argument leading to (\ref{eq16}). Since $C$ does not depend
on $\smallu$, (\ref{eq14}) follows.
\end{pf}

%s6.3 #&#
\subsection{Bounds on the number of ancestors, descendants and
pairwise distances}
\label{sschar}
Here we provide bounds needed to prove the compact containment
condition for the TMMMS. We use the notation from
Definitions~\ref{defgraph},~\ref{defTMMMS} and
\ref{defTMMMS2}. Most importantly, $\mathcal U^N = (\mathcal
U_t^N)_{t\geq0}$ with $\mathcal U_t^N =
\overline{(U_N,r^N_t,\mu^N_t)}$ is the TMMMS, and we use $A_s(l,t)$ to
denote the ancestor of $(l,t)$ at time $s$.

The key to compact containment conditions for tree-valued processes
arising in the context of population models is to control the number
of ancestors times $\ve>0$ in the past and the number of descendants
of some given subpopulation uniformly in the relevant parameter (here
$N$); see Section~\ref{subproof1}. For both we provide the needed
estimates here.

The following birth and death process, more precisely its infimum, serves
as an upper bound on the number of ancestors in the Moran model with mutation
and selection.
%
%de6.8 #&#
%
\begin{definition}[(The processes $\mathcal J$ and $\mathcal J^\ast
$)]\label{equpp}
Let $\mathcal J = (J_t)_{t\geq0}$ be the homogeneous
Markov jump process which jumps
%
%e6.22 #&#
%
\begin{eqnarray}
\label{eq806}
&\mbox{from $j$ to $j+1$}\qquad \mbox{at rate
$j\alpha$},&\nonumber\\[-8pt]\\[-8pt]
&\mbox{from $j$ to $j-1$}\qquad \mbox{at rate
$\gamma\pmatrix{j\cr2}$}.&
\nonumber
\end{eqnarray}
Moreover, we define $\mathcal J^\ast= (J_t^\ast)_{t\geq0}$ by
$J_t^\ast:= \inf_{0\leq s\leq t} J_s. $
\end{definition}
%
%pr6.9 #&#
%
\begin{proposition}[(An upper bound for the number of ancestors)]
\label{Pupp1}
Let $\mathcal U^N = (\mathcal U_t^N)_{t\geq0}$ be
the TMMMS as well as $\mathcal J^\ast= (J_s^\ast)_{s\geq0}$ from
Definition~\ref{equpp}, started in $J_0^\ast= J_0=j\in\mathbb
N$. For $0\leq s\leq t$ and $n_1,\ldots,n_j\in U_N$ pairwise different,
set
%
%e6.23 #&#
%
\begin{equation}
A^{j,N}_{s,t}:= \# \{A_s(n_i, t)\dvtx i=1,\ldots,j\}.
\end{equation}
Then
%
%e6.24 #&#
%
\begin{equation}\label{agx3}
A^{j,N}_{s,t} \leq J_{t-s}^\ast\qquad \forall 0\leq s \leq t,
N \in\N\qquad\mbox{stochastically}.
\end{equation}
\end{proposition}
\begin{pf}
Look at the graphical construction of the Moran model with mutation
and selection at time $t$. Following the ancestral lines of
$n_1,\ldots,n_j$ backward, two things might occur at some time $s$: at
a resampling arrow between two ancestral lines, these ancestral
lines have a common ancestor, and $A_{s,t}^{j,N}$ decreases by
one. The rate of such an event is proportional to $\gamma$ and the
number of pairs. If an ancestral line hits the tip of a selective
arrow, there are two possible ancestors, one of which is the real
one depending on the types of the two. The process $\mathcal J$
counts both of them which certainly gives an upper bound for the
number of ancestors. This proves that $A_{s,t}^{j,N}\leq
J_{t-s}$ stochastically. Moreover, the number of ancestors can never
increase when going back in time, and hence, $A_{s,t}^{j,N}\leq
J_{t-s}^\ast$ follows.
\end{pf}
%
%co6.10 #&#
%
\begin{corollary}[(The number of ancestors of the total population)]
\label{corupp2}
For $0\leq s<t$,
\[
\mathbf E[A^{N,N}_{s,t}] \leq\frac{(\gamma+ 2 \alpha)
e^{(\gamma/2 + \alpha) (t-s)} N}{2 \alpha+ \gamma+
\gamma(e^{(\gamma/2 + \alpha) (t-s)}-1)N}
\stackrel{N\to\infty}{\longrightarrow}\frac{(\gamma+ 2 \alpha)
e^{(\gamma/2 +
\alpha) (t-s)} }{\gamma(e^{(\gamma/2 + \alpha) (t-s)} -1)}.
\]
\end{corollary}
\begin{pf}
Set $J_0=N$. Writing $y(s)=\mathbf E[J_s]$ and using the backward
equation, we have
%
%e6.25 #&#
%
\begin{equation}\label{eq806a}
\dot y(s) = \alpha\mathbf E[J_s] - \gamma\mathbf
E\biggl[\pmatrix{J_s\cr2}\biggr] \le\biggl(\frac12\gamma+ \alpha
\biggr)y(s)
- \frac12 \gamma y(s)^2,
\end{equation}
where we used Jensen's inequality in the last step. The solution of
the initial value problem
%
%e6.26 #&#
%
\begin{equation}\label{eq807}
\dot z = \bigl(\tfrac12\gamma+ \alpha\bigr)z - \tfrac12\gamma
z^2,\qquad z(0)=N
\end{equation}
is given by
%
%e6.27 #&#
%
\begin{equation}\label{eq808}
z(s) = \frac{(\gamma+ 2 \alpha) e^{(\gamma/2 + \alpha) s} N}{2
\alpha+ \gamma+ \gamma(e^{(\gamma/2 + \alpha) s}-1)N}.
\end{equation}
The last three equations together with Proposition~\ref{Pupp1} give
the assertion.
\end{pf}

% For

% \begin{corollary}
% For $\smallu\in\mathbb M^I$, let $\mathcal U^\smallu= (\mathcal
% U_t^\smallu)_{t\geq0}$ be the TFVMS, started in $\mathcal
% U_0^\smallu= \smallu$. Then, the set $\{A_{s,t}(\mathcal
% U_t^\smallu): \smallu\in\mathbb M^I, 0\leq s\leq t\}$ is tight.
% \end{corollary}

% \begin{pf}
% xxx Approximation
% \end{pf}

Our next task is to bound the frequency of descendants.

%. This notion
% is easily defined for the TMMMS.
%
%de6.11 #&#
%
\begin{definition}[(Frequency of descendants in TMMMS and
filtration)]
\label{deffreqanc}
(1)~Let $\mathcal U^N:=(\mathcal U^N_t)_{t\ge0}$ be the TMMMS
with population size $N$ defined by the graphical
construction. For $s\leq t$ and $\mathcal V\subseteq U_N$, we
define
%
%e6.28 #&#
%
\begin{equation}\label{eq914}
D_t^N(\mathcal V,s):=\{l \in U_N\dvtx A_s(l,t)\in\mathcal
V\},
\end{equation}
the set of descendants of $\mathcal V$ at time $t$.

(2) For the TMMMS $\mathcal U^N = (\mathcal U_t^N)_{t\geq0}$,
recall the Poisson processes $\eta^{\mathrm{res}},
\eta^{\mathrm{mut}},\break
\eta^{\mathrm{sel}}$ on $U_N \times\mathbb R_+$ and $S_N(t) = U_N
\times(-\infty, t]$ from Definition~\ref{defgraph}. We
define the filtration $(\mathcal A^N_t)_{t\geq0}$ by $\mathcal
A_t^N = \sigma(\eta^{\mathrm{res}}|_{S_N(t)},
\eta^{\mathrm{mut}}|_{S_N(t)}, \eta^{\mathrm{sel}}|_{S_N(t)})$.
\end{definition}
%
%le6.12 #&#
%
\begin{lemma}[(Bounds on the frequency of descendants)]%
\label{lemfr-desc}
For $0<\varepsilon\leq T$ there is
$\delta>0$ such that for $0 \le s \le T$ and any sequence $(\mathcal
V^N)_{N\in\mathbb N}$ of $\mathcal A_s^N$-measurable subsets of
$U_N$, we have
%
%e6.29 #&#
%
\begin{equation} \label{eq17} \limsup_{N \to\infty} \mu
_s^N(\mathcal
V^N) \leq\delta \quad\Longrightarrow\quad \limsup_{N \to\infty}
\mathbf{P}\Bigl(\sup_{s \le t \le T} \mu_t(D_t^N(\mathcal V^N,s)) >
\varepsilon\Bigr) \le\varepsilon.\hspace*{-35pt}
\end{equation}
\end{lemma}
\begin{pf}
By time-homogeneity of the TMMMS, it suffices to show the assertion
for $s=0$. We restrict ourselves to the haploid case. The extension
to the diploid case is straightforward. The proof is based on a
\textit{coupling argument} that we describe next.

For $N\in\mathbb N$, consider the graphical construction of
$\mathcal U^N = (\mathcal U_t^N)$, given by means of the Poisson
processes $(\eta^{\mathrm{res}}, \eta^{\mathrm{mut}},
\eta^{\mathrm{sel}})$. Moreover, let $\mathcal V^N$
satisfy the assumption on the left-hand side of (\ref{eq17}). We
define a process $\overline{\mathcal U}^N = (\overline{\mathcal
U}_t^N)_{t\geq0}$ with $\overline{\mathcal U}{}^N_t =
\overline{(U_N, \overline r^N_t, \overline\mu^N_t)}$, taking values
in $\mathbb U^{\{\bullet,{\darkgreybullet}\}}$ with the
following features:
\begin{longlist}[(iii)]
\item[(i)] for $k\in\mathcal V^N$, set $u_k(0) = \bullet$, for
$k\notin\mathcal V^N$, set $u_k(0)={\darkgreybullet}$,\vadjust{\goodbreak}
\item[(ii)] $\chi(\bullet)=1, \chi({\darkgreybullet})=0$,
that is, only $\bullet$ can use events in $\eta^{\mathrm{sel}}$,
\item[(iii)] $\vartheta=0$, that is, mutation is absent.
\end{longlist}
For the dynamics of $\overline{\mathcal U}{}^N$, use the same Poisson
processes $\eta^{\mathrm{res}}$ and $\eta^{\mathrm{sel}}$ as
$\mathcal U^N$.
Note that $(X_t^N)_{t\geq0}$, given by $X_t^N =
\overline{\mu}_t(D_t(\mathcal V^N,0))$ is a Markov jump process with
transitions
\begin{eqnarray*}
% \label{eq931}
&\displaystyle \mbox{from } x \mbox{ to }x+\frac1N \qquad
\mbox{at rate }\frac\gamma2
N^2x(1-x) + \alpha N x(1-x),&\\
&\displaystyle \mbox{from } x \mbox{ to }x-\frac1N
\qquad\mbox{at rate }\frac\gamma2
N^2x(1-x).&
\end{eqnarray*}
In particular, $(X^N_t)_{t\geq0}$ converges weakly (with respect to
the Skorohod topology) to the solution $(X_t)_{t\geq0}$ of the SDE
%
%e6.30 #&#
%
\begin{equation}
\label{pp1002}
dX = \alpha X(1-X)\,dt + \sqrt{\gamma X(1-X)}\,dW.
\end{equation}
By construction of $\overline{\mathcal U}{}^N$, we find that
$\mu_t(D_t(\mathcal V^N,0)) \leq X_t^N$, and hence, if
$\limsup_{N\to\infty} \mu_0^N(\mathcal V^N) \leq\delta$ for some
$\delta>0$, then
%
%e6.31 #&#
%
\begin{eqnarray}
\label{eq925}
\limsup_{N\to\infty} \mathbf{P}\Bigl(\sup_{0 \le t \le T}
\mu_t(D_t(\mathcal V^N,0)) > \varepsilon\Bigr) & \le &
\limsup_{N\to\infty} \mathbf{P}\Bigl(\sup_{0 \le t \le T}
X_s^N>\varepsilon\Bigr) \nonumber\\[-8pt]\\[-8pt]
& \le &\mathbf{P}\Bigl(\sup_{0 \le t
\le T}X_s>\varepsilon\big| X_0=\delta\Bigr).
\nonumber
\end{eqnarray}
By Doob's maximal inequality, for each $\varepsilon>0$, we find
$\delta>0$ such that
%
%e6.32 #&#
%
\begin{equation}\label{eq926}
\mathbf{P}\Bigl(\sup_{0 \le t \le T} X_s\big| X_0=\delta\Bigr) \leq
\varepsilon,
\end{equation}
and the result follows.
\end{pf}

The next result is a corollary of the previous lemma and
Proposition~\ref{Pupp1}.

%co6.13 #&#
%
\begin{corollary}[(Tightness of pairwise distances)]%
\label{corwtight}
$\!\!\!$Assume that $(\mathcal U^N_0)_{N \in\N}$ is
tight. Let $R_{12}^N(t)=\langle\nu^{N, \mathcal U^N_t},
r_{12}\rangle$. For any $\varepsilon>0$, there is $C=C(\varepsilon)
< \infty$ such that for all $t \ge0$,
%
%e6.33 #&#
%
\begin{equation} \label{eq19} \limsup_{N \to\infty} \mathbf
P\bigl(R_{12}^N(t) > C\bigr) \le\varepsilon.
\end{equation}
\end{corollary}
\begin{pf}
Let $\varepsilon>0$ be given. For the process $\mathcal J$ from
Definition~\ref{sschar} with \mbox{$J_0=2$}, let $T_1=\inf\{t>0\dvtx
J_t=1\}$. As a birth and death process with quadratic death and
linear birth rates, $\mathcal J$ is recurrent and irreducible.
Choose $C_1>0$ so that
%
%e6.34 #&#
%
\begin{equation}
\label{eq34}
\mathbf P\biggl(T_1 > \frac{C_1}2 \biggr) \le\varepsilon.
\end{equation}
For $C_2>0$ and $\mathcal U_0^N = \overline{(U^N_0, r^N_0,
\mu^N_0)}$, consider the family of subsets of
$U_0^N$
\[
\mathcal W_{C_2}^N:= \{W\subseteq U_0^N\dvtx r(g_1, g_2)\leq C_2
\mbox{ for all }g_1,g_2 \in W\}.\vadjust{\goodbreak}
\]
Clearly, $\mathcal W_{C_2}^N$ contains maximal elements (with
respect to ``$\subseteq$''), and we denote by $W_{C_2}^N$ an arbitrary
maximal element of $\mathcal W_{C_2}^N$. Set $V^N_{C_2}= U^N_0
\setminus W^N_{C_2}$. By the tightness assumption and
Lemma~\ref{lemfr-desc}, we may choose $C_2$ and $\delta>0$ such that
\[
% \label{eq36}
\limsup_{N \to\infty} \mu_s^N (V^N_{C_2}) \le\delta
\quad\Longrightarrow\quad\limsup_{N \to\infty} \mathbf{P}\Bigl(\sup_{s \le
t \le C_1/2} \mu_t^N(D_t^N(V^N_{C_2})) > \varepsilon\Bigr) \le
\varepsilon.
\]
To continue we have to distinguish whether $t \in[0, C_1/2]$ or
not.

For $t \in[0, C_1/2]$ the event $\{R^N_{12}(t) > C_1+C_2\}$ means
that the ancestral lines of a pair of individuals drawn at time $t$
did not coalesce in the time interval $[0,t]$ and that the distance
of their ancestors at time $0$ is at least $C_1+C_2-2t \ge C_2$. By
the choice of $C_1$ and $C_2$, we have
%
%e6.35 #&#
%
\begin{equation} \label{eq20}
\limsup_{N \to\infty} \mathbf P \bigl(R_{12}^N(t) > C_1+C_2 \bigr) <
\varepsilon \qquad\mbox{for all $t \in[0,C_1/2]$}.
\end{equation}
In the case $t > C_1/2$ the event $\{R^N_{12}(t) > C_1\}$ means that
a randomly chosen pair of ancestral lines did not coalesce in the
time interval $[t-C_1/2,t]$, that is,
%
%e6.36 #&#
%
\begin{equation}
\label{pp2010}
\{R^N_{12}(t) > C_1\} = \{A^{2,N}_{t-C_1/2,t}=2 \}.
\end{equation}
By Proposition~\ref{Pupp1} and the choice of $C_1$ it follows that
for $t >
C_1/2$ (independent of $N$),
%
%e6.37 #&#
%
\begin{equation} \label{eq35} \mathbf P\bigl( R^N_{12}(t) >C_1\bigr) = \mathbf
P(A^{2,N}_{t-C_1/2,t}=2)\leq\mathbf P\biggl(T_1 >
\frac{C_1}2 \biggr) \le\varepsilon.
\end{equation}
Combining (\ref{eq20}) and (\ref{eq35}) we obtain (\ref{eq19})
with $C=C_1+C_2$.
\end{pf}

%s7 #&#
\section{\texorpdfstring{Proofs of Theorems \protect\ref{T1}, \protect\ref{T3} and \protect\ref{T4}}
{Proofs of Theorems 1, 3 and 4}}
\label{SPPP2}
Now we have all ingredients for the proofs of our main Theorems
\ref{T1},~\ref{T3} and~\ref{T4}.

%s7.1 #&#
\subsection{\texorpdfstring{Proof of Theorems \protect\ref{T1} and \protect\ref{T3}}
{Proof of Theorems 1 and 3}}
\label{subproof1}
We prove Theorems~\ref{T1} and~\ref{T3} simultaneously. The main step
in the proof is to show that the family of processes $\{\mathcal
U^N\dvtx
N\in\mathbb N\}$ is tight and that all limit points solve the
$(\mathbf P_0, \Omega, \Pi^1)$-martingale problem and fulfill (b) of
Theorem~\ref{T1}. Uniqueness of the solution of the $(\mathbf P_0,
\Omega, \Pi^1)$-martingale problem
% follows either from Theorem~\ref{T2} or alternatively
is a consequence of the duality relation given by
Proposition~\ref{propdual}(2) [see
\citet{EthierKurtz86}, Proposition 4.4.7]. Note that the set of
duality functions from (\ref{eq905}) is separating on $\mathbb M^I$
by Proposition~\ref{propdual}(1). Finally, properties (a) and (e) from
Theorem~\ref{T1} are direct consequences of Propositions~\ref{Pqv}
and~\ref{PfirstSecond}.

In order to establish tightness of $\{\mathcal U^N\dvtx N\in\mathbb N\}$
and property (b) of Theorem~\ref{T1}, we use Lemma 4.5.1 and
Remark 4.5.2 of \citet{EthierKurtz86}, requiring us to check two
conditions: a convergence relation for generators and a compact
containment condition. To verify the first, recall that we showed
convergence of generators of TMMMS to the generator of TFVMS
in Proposition~\ref{propconv2TVFV}.\vadjust{\goodbreak}

Hence, we have to verify the second condition amounting to show the
following compact containment conditions: for all $\varepsilon, T>0$
there exist sets $\Gamma_{\varepsilon,T}\subseteq\mathbb U^I_c$,
relatively compact in $\mathbb U^I_c$ and $\widetilde
\Gamma_{\varepsilon,T}\subseteq\mathbb U^I$, relatively compact in
$\mathbb U^I$, such that
%
%e7.1 #&#
%
\begin{eqnarray}
\label{eq902}
\inf_{N\in\mathbb N} \mathbf P(\mathcal U_t^N\in
\Gamma_{\varepsilon, T} \mbox{ for all } \varepsilon\leq t\leq T)
&>& 1-\varepsilon,\nonumber\\[-8pt]\\[-8pt]
\inf_{N\in\mathbb N} \mathbf P(\mathcal U_t^N\in\widetilde
\Gamma_{\varepsilon, T} \mbox{ for all } 0\leq t\leq T) &>&
1-\varepsilon.
\nonumber
\end{eqnarray}
For $\smallx= \overline{(X,r,\mu)}$, we set
$\pi_1(\smallx):=\overline{(X,r,(\pi_X)_\ast\mu)}$. Since $I$ is
compact, it is a consequence of Theorem 3 in \citet{DGP-topo2011}, that
$\Gamma_{\varepsilon, T}\subseteq\mathbb U^I$
($\widetilde\Gamma_{\varepsilon, T}\subseteq\mathbb U_c^I$) is
relatively compact in $\mathbb U^I$ ($\mathbb U_c^I$) if and only if
$\pi_1(\Gamma_{\varepsilon, T})$
[$\pi_1(\widetilde\Gamma_{\varepsilon, T})$] is relatively compact in
$\mathbb U$ ($\mathbb U_c$).

In order to check existence of $\Gamma_{\varepsilon, T}$
($\widetilde\Gamma_{\varepsilon, T}$) such that (\ref{eq902}) holds
with $\mathcal U_t^N$ replaced by $\pi_1(\mathcal U_t^N)$ and
$\Gamma_{\varepsilon, T}$ ($\widetilde\Gamma_{\varepsilon, T}$)
replaced by $\pi_1(\Gamma_{\varepsilon, T})$
[$\pi_1(\widetilde\Gamma_{\varepsilon, T})$], we use Proposition 2.22
of \citet{GrevenPfaffelhuberWinter2011}. This result gives a
condition for (\ref{eq902}), based on estimates on the number of
ancestors time $\varepsilon>0$ in the past and in terms of frequencies
of descendants of rare ancestors. First, we note that $(\pi_1(\mathcal
U_t^N))_{t\geq0}$ fits the definition of a tree-valued version of a
population model from Proposition 2.18 of
\citet{GrevenPfaffelhuberWinter2011}. For (i) of that proposition,
the required bound on the frequency of descendants is given in
Lemma~\ref{lemfr-desc}. Moreover, (ii) of that proposition is a
consequence of Corollary~\ref{corupp2}. Hence, (\ref{eq902}) follows.

Except for (c) and (d) of Theorem~\ref{T1} the proof of
Theorems~\ref{T1} and~\ref{T3} is complete by the above arguments. To
prove the Feller property of $\mathcal U$, part (c) of
Theorem~\ref{T1}, we use duality. Let $\mathcal U^\smallu=(\mathcal
U_t^\smallu)_{t\geq0}$ be the TFVMS started in $\mathcal U_0^\smallu
= \smallu$ and $\smallu, \smallu_1,\smallu_2,\ldots\in\mathcal
U^I$ be
such that $\smallu_n \stackrel{n\to\infty}{\longrightarrow}\smallu
$ in the
Gromov-weak topology and let $t>0$ be fixed. First we note that for
$\Phi= \Phi^{n,\phi}\in\Pi^1$,
\[
\mathbf E[\Phi(\mathcal U_t^{\smallu_n})] = \mathbf E[\langle
\nu^{\mathcal U_t^{\smallu_n}}, \phi\rangle] = \mathbf E[\langle
\nu^{\smallu_n}, \Xi_t\rangle] \stackrel{n\to\infty
}{\longrightarrow} \mathbf
E[\langle\nu^{\smallu}, \Xi_t\rangle] = \mathbf E[\Phi(\mathcal
U_t^{\smallu})]
\]
by Proposition~\ref{propdual}, where $\Xi= (\Xi_t)_{t\geq0}$ is the
dual process from Definition~\ref{defdual} with $\xi_0 =
\phi$. Hence, by Theorem 5 in \citet{DGP-topo2011}, $\mathcal
U_t^{\smallu_n} \xRightarrow{n\to\infty} \mathcal U_t^{\smallu}$ and
the Feller property follows.

For (d) in Theorem~\ref{T1} notice that the strong Markov property follows
from the Feller property by standard theory [e.g., Theorem 4.2.7 in
\citet{EthierKurtz86}, and note that local compactness of the state
space is not
used in the proof].

% \blu{For (f) we note that due to the well-posedness it suffices to
% give a solution with this property. Consider first the case for
% $\alpha=0$ and then use Girsanov to get the general case.}

%s7.2 #&#
\subsection{\texorpdfstring{Proof of Theorem \protect\ref{T4}}{Proof of Theorem 4}}
As observed before Theorem~\ref{T4}, a unique equilibrium for
$\mathcal U$ implies a unique equilibrium for $\widetilde\zeta$, so we
are left with showing the converse.

If we have convergence from every initial point to a limiting law, then
this law is the unique invariant measure of the process. In order to
see that the limiting law is invariant, consider\vadjust{\goodbreak}
$f\in\overline{\mathcal C}(\mathbb U^I)$, and let $(S_t)_{t\geq0}$ be
the semigroup of the TFVMS. Since the map $\smallu\mapsto
(S_tf)(\smallu)$ is continuous by Theorem~\ref{T1}.c, the limiting law
is invariant using the same argument as in
Proposition 1.8(d) of \citet{Liggett85}. Hence we have to
establish the
convergence statement. Recall that the family $\{\smallu\mapsto
\langle\nu^\smallu, \xi\rangle\dvtx \xi\in\Upsilon\}$ is separating
$\mathcal M_1(\mathbb U^I)$; see Proposition~\ref{propdual}. Hence
we have to show two assertions [see, e.g., \citet{EthierKurtz86}, Lemma
3.4.3]:
\begin{longlist}[(ii)]
\item[(i)] The family $\{\mathcal U_t\dvtx t>1\}$ is tight in $\mathbb
U_c^I$.
\item[(ii)] For all $\xi\in\Upsilon$, $\lim_{t\to\infty}\mathbf
E_\smallu[\langle\nu^{\mathcal U_t}, \xi\rangle]$ exists and does
not depend on $\smallu$.
\end{longlist}
When these two properties hold, we conclude from (i) that there are
convergent subsequences of $(\mathcal U_t)_{t\geq0}$. Let
$\smallu\in\mathbb M_I$ and $t_1,t_2,\ldots$ be such that $\mathcal
U_{\infty}$ is the weak limit of $(\mathcal U_{t_n})_{n=1,2,\ldots}$,
started in $\smallu$. Then, (ii) implies that, for all
$\Phi\in\Pi^1$ with $\Phi(\smallu) = \langle\nu^\smallu, \xi
\rangle$
and $\xi\in\Upsilon$
%
%e7.2 #&#
%
\begin{equation}
\label{eqdual2}
\mathbf E_{\smallu}[\Phi(\mathcal U_{\infty})] = \lim_{n\to\infty}
\mathbf E_{\smallu}[\langle\nu^{\mathcal U_{t_n}}, \xi\rangle] =
\lim_{t\to\infty} \mathbf E_{\smallu}[\langle\nu^{\mathcal
U_{t}}, \xi\rangle]
\end{equation}
exists and is independent of $\smallu$.

We start by proving (i). By Theorem 4 in \citet{DGP-topo2011}, we need
to show that $\{\pi_1(\mathcal U_t)\dvtx t>1\}$ is tight in $\mathbb
U_c$. For this, we use Proposition 6.2
of \citet{GrevenPfaffelhuberWinter2011}. In particular, we have to
check that:
\begin{longlist}[(2)]
\item[(1)] $\{R^N_{12}(t)\dvtx t>1\}$ is tight,
\item[(2)] $\{A_{t-\varepsilon,t}\dvtx t>1\}$ is tight for $0<\varepsilon
<1$, where $A_{t-\varepsilon,t}$ from Definition~\ref{defgraph} is
the number of ancestors of $\mathcal U_t$ at time $t-\varepsilon$,
or, equivalently, the number of $2\varepsilon$-balls needed to cover
$\mathcal U_t$.
\end{longlist}
Once (1) and (2) are shown, let $\delta>0$. It is
straightforward to construct a set $\Gamma_\delta\subseteq\mathbb
U_c$ which fulfills (i) and (ii) of Proposition 6.2
of \citet{GrevenPfaffelhuberWinter2011} with $\inf_{t>1}\mathbf
P(\mathcal U_t \in\Gamma_\delta)>1-\delta$. While (1) is true by
Corollary~\ref{corwtight}, (2)~holds according to
Corollary~\ref{corupp2}.

We now show (ii) if $\widetilde\zeta$ has a unique
equilibrium. Consider the process $\Xi= (\Xi_t)_{t\geq0}$ from
Definition~\ref{defdual}. Recall from Proposition~\ref{lnoninc}(3)
that there is an almost surely finite $T$ such that $\Xi_T$ does not
depend on $\underline{\underline r}$. We use the duality relation
from Proposition~\ref{propdual} and the strong Markov property of
$\Xi$ to see that for $\Xi_0 = \xi\in\Upsilon$,
%
%e7.3 #&#
%
\begin{eqnarray}
\label{eq1211}
\lim_{t\to\infty} \mathbf E_\smallu[\langle\nu^{\mathcal U_t},
\xi\rangle] & = & \lim_{t\to\infty} \mathbf E_\xi[\langle
\nu^{\smallu},\Xi_t\rangle] = \lim_{t\to\infty} \mathbf
E_\xi[\mathbf E_{\Xi_T} [\langle\nu^{\smallu},\Xi_t\rangle]]
\nonumber\\[-8pt]\\[-8pt]
& = & \lim_{t\to\infty} \int\mathbf E_\smallu[\langle\nu^{\mathcal
U_t}, \widetilde\xi\rangle] \mathbf P_\xi(\Xi_T \in
d\widetilde\xi)
\nonumber
\end{eqnarray}
exists and does not depend on $\smallu$. This holds since for
$\widetilde\xi\in\Upsilon$, not depending on $\underline
{\underline
r}$, the limit $\lim_{t\to\infty}\mathbf E_\smallu[\langle
\nu^{\mathcal U_t}, \widetilde\xi\rangle]$ exists and is independent
of $\smallu$ since $\widetilde\zeta$ has a unique equilibrium. Note
that $t\mapsto\|{\Xi_t}\|_{\infty}$ is nonincreasing by Proposition~\ref{lnoninc}(1)
and therefore, all expectations in (\ref{eq1211})
are well defined.\vadjust{\goodbreak}

Next, we show that (ii) holds if $\vartheta>0, \alpha>0$ and mutation
has a parent-independent component, again using the dual process $\Xi
= (\Xi_t)_{t\geq0}$ from Definition~\ref{defdual}. From Proposition
\ref{lnoninc}(2) we know that $\Xi$ converges almost surely to a
(random) constant function $\Xi_{\infty}$ taking values in
$\overline{\mathcal C}_0^1$. Hence, for $\Xi_0 = \xi\in\Upsilon$,
%
%e7.4 #&#
%
\begin{equation}
\label{eqdual2p}
\lim_{t\to\infty} \mathbf E_{\smallu}[\langle\nu^{\mathcal
U_{t}}, \xi\rangle] = \lim_{t\to\infty} \mathbf E_{\xi}[\langle
\nu^{\smallu}, \Xi_{t}\rangle] = \mathbf
E_\xi[\langle\nu^{\smallu},\Xi_{\infty}\rangle] = \mathbf
E_\xi[\Xi_{\infty}],
\end{equation}
where the expression on the right-hand side does not depend on
$\smallu$. Again, note that $t\mapsto\|{\Xi_t}\|_{\infty}$ is
nonincreasing by Proposition~\ref{lnoninc}(1) and therefore, all
expectations in (\ref{eqdual2p}) are well defined. Hence, (ii)
follows if either $\widetilde\zeta$ is ergodic or if mutation has an
independent part, and this completes the proof of Theorem~\ref{T4}.

%s7.3 #&#
\subsection{\texorpdfstring{Proof of Theorem \protect\ref{T2}}{Proof of Theorem 2}}
\label{SPPP3}
Before we turn to the proof of Theorem~\ref{T2}, we recall the
Girsanov transform for continuous semimartingales from
\citet{Kallenberg2002}, Theorems 18.19 and 18.21.
%
%le7.1 #&#
%
\begin{lemma}[(The Girsanov theorem for continuous semimartingales)]%
\label{lGirs}
Let $\mathcal M=(M_t)_{t\geq0}$ be a
continuous $\mathbf P$-martingale for some probability measure
$\mathbf P$, and assume $\mathcal Z = (Z_t)_{t\geq0}$, given by $Z_t
= e^{M_t - (1/2) [\mathcal M]_t}$, is a martingale. If $\mathcal
N = (N_t)_{t\geq0}$ is a local $\mathbf P$-martingale, and $\mathbf
Q$ is defined via its Radon--Nikodym derivative with respect to
$\mathbf P$, that is, $\frac{d\mathbf Q}{d\mathbf P}|_{\mathcal
F_t} = Z_t$, then $\mathcal N - [\mathcal M, \mathcal N]$ is a
local $\mathbf Q$-martingale. (Here, $[\mathcal M, \mathcal N]$ is
the covariation process between $\mathcal M$ and $\mathcal N$ and
$[\mathcal M] = [\mathcal M, \mathcal M]$.)
\end{lemma}
\begin{pf*}{Proof of Theorem~\ref{T2}}
Since $|{\alpha' -\alpha}| <\infty$, $\mathcal M$ is bounded, and
therefore the right-hand side of (\ref{eqdQdP}) is a
martingale. Thus $\mathbf Q$ is well defined.

By Theorem 5 from \citet{DGP-topo2011}, $\Pi^1$ contains an algebra
that separates points, so the TFVMS fulfills the assumptions of
Proposition~\ref{Pqv}. The generator $\Omega_\alpha$ is second
order by Proposition~\ref{PfirstSecond}, and its only second order
term is $\Omega^{\mathrm{res}}$. In particular, we can use
Corollary~\ref{corcov}. This is important since the additional
drift term introduced by the Girsanov change of measure is given by
a covariation; see Lemma~\ref{lGirs}. We have to compute
$[\Phi(\mathcal U), \Psi(\mathcal U)]$ for $\Phi(\mathcal U) =
(\Phi(\mathcal U_t))_{t\geq0}, \Psi(\mathcal U) = (\Psi(\mathcal
U_t))_{t\geq0}$ for $\Phi\in\Pi^1$ and $\Psi$ from
(\ref{eqPsi}). We take $\Phi\in\Pi^1_{n}$ and compute, using the
symmetry of $\chi'$,
\begin{eqnarray*}%\label{dgpGir17}
&&\Omega^{\mathrm{res}} \bigl(\Phi(\smallu) \cdot\Psi(\smallu)\bigr) -
\Psi(\smallu) \cdot\Omega^{\mathrm{res}} \Phi(\smallu) -
\Phi(\smallu) \cdot\Omega^{\mathrm{res}} \Psi(\smallu) \\
&&\qquad =
\frac{\alpha'-\alpha} \gamma\bigl( \Omega^{\mathrm{res}} \langle
\nu^\smallu, \phi\cdot(\chi'_{1,2}\circ\rho_1^n) \rangle-
\langle\nu^\smallu, \chi'_{1,2}\rangle\cdot
\Omega^{\mathrm{res}}\langle\nu^\smallu, \phi\rangle\\
&&\hspace*{140.5pt}\qquad\quad{}
- \langle\nu^\smallu, \phi\rangle\cdot
\Omega^{\mathrm{res}} \langle\nu^\smallu, \chi'_{1,2}\rangle
\bigr) \\
&&\qquad = \frac{\alpha'-\alpha} 2 \Biggl( \sum_{k,l=1}^{n+2} \langle
\nu^\smallu, \phi\cdot(\chi'_{1,2}\cdot\rho_1^n)\circ
\theta_{k,l} - \phi\cdot(\chi'_{1,2}\circ\rho_1^n)\rangle\\
&&\hspace*{36pt}\qquad\quad{}
- \sum_{k,l=1}^n \langle\nu^\smallu, (\phi\circ
\theta_{k,l})\cdot(\chi'_{1,2}\circ\rho_1^n) -
\phi\cdot(\chi'_{1,2}\circ\rho_1^n)\rangle\\
&&\hspace*{58pt}\qquad\quad{}
- 2\langle\nu^\smallu, \phi\cdot( \chi'_{1,2}
\circ\theta_{1,2} \circ\rho_1^n)
- \phi\cdot(\chi'_{1,2}\circ\rho_1^n)\rangle\Biggr) \\
&&\qquad = (\alpha'-\alpha) \sum_{k=1}^{n} \langle\nu^\smallu,
(\phi\cdot\chi'_{n+1,n+2})\circ\theta_{k,n+1} - \phi\cdot
\chi'_{n+1,n+2}\rangle\\
&&\qquad = (\alpha'-\alpha) \sum_{k=1}^{n}
\langle\nu^\smallu, \phi\cdot\chi'_{k,n+1} - \phi\cdot
\chi'_{n+1,n+2} \rangle\\
&&\qquad = \Omega^{\mathrm{sel}}_{\alpha'}
\Phi(\smallu) - \Omega^{\mathrm{sel}}_{\alpha} \Phi(\smallu).
\end{eqnarray*}
Hence, Corollary~\ref{corcov} implies that
%
%e7.5 #&#
%
\begin{equation}\label{eqad16}
[\Phi(\mathcal U), \mathcal M]_t = [\Phi(\mathcal U),
\Psi(\mathcal U)]_t = \int_0^t \bigl(
\Omega_{\alpha'}^{\mathrm{sel}}\Phi(\mathcal U_s) -
\Omega^{\mathrm{sel}}_{\alpha} \Phi(\mathcal U_s) \bigr) \,ds,
\end{equation}
where $\mathcal U=(\mathcal U_t)_{t\geq0}$ is a solution of the
$(\mathbf P_0, \Omega_\alpha, \Pi^1)$-martingale problem. For any
$\Phi\in\Pi^1$,
%
%e7.6 #&#
%
\begin{equation}\label{dgpGir18a}
\mathcal N_\Phi:= \biggl(\Phi(\mathcal U_t) - \int_0^t
\Omega_\alpha\Phi(\mathcal U_s)\,ds\biggr)_{t\geq0}
\end{equation}
is a continuous $\mathbf P$-martingale. Thus, by Girsanov's theorem
for continuous semimartingales, Lemma~\ref{lGirs} and
(\ref{eqad16}), we see that
\begin{eqnarray*}
&&\biggl( \Phi(\mathcal U_t) - \int_0^t \Omega_\alpha\Phi(\mathcal
U_s)\,ds - [\Phi(\mathcal U), \Psi(\mathcal U)]_t \biggr)_{t\geq0}\\
&&\qquad =
\biggl( \Phi(\mathcal U_t) - \int_0^t \Omega_{\alpha'} \Phi
(\mathcal
U_s)\,ds\biggr)_{t\geq0}
\end{eqnarray*}
is a $\mathbf Q$-martingale for $\mathbf Q$ defined by
(\ref{eqdQdP}). Since $\Phi\in\Pi^1$ was arbitrary, it follows that
$\mathbf Q$ solves the $(\mathbf P_0, \Omega_{\alpha'}, \Pi^1)$-martingale
problem.
\end{pf*}

%s8 #&#
\section{\texorpdfstring{Proof of Theorem \protect\ref{T5}}{Proof of Theorem 5}}
\label{SproofAppl}
If $\mathcal U_{\infty}^\alpha$ is as in Theorem~\ref{T5}, the proof is
based on the fact that
%
%e8.1 #&#
%
\begin{equation}\label{eq707}
\mathbf E[\Omega_\alpha\Phi(\mathcal U_{\infty}^\alpha)] = 0
\end{equation}
for $\Phi\in\Pi^1$. (This follows easily from the
$\Omega_\alpha$-martingale problem.) Moreover, for small $\alpha>0$,
the equilibrium $\mathcal U_{\infty}^\alpha$ is close to the equilibrium
without selection, and the equilibrium under neutrality is well
understood. In order to use this knowledge for the neutral case, the
following fact is fundamental.

%le8.1 #&#
%
\begin{lemma}[(Continuity of $\alpha\mapsto\mathcal U_{\infty}^\alpha$)]
\label{l94}
Let $\mathcal U_{\infty}^\alpha$ be as in
Theorem~\ref{T5}. Then, for $\Phi\in\Pi^1$,
%
%e8.2 #&#
%
\begin{equation}\label{agx5}
\mathbf E[\Phi(\mathcal U_{\infty}^\alpha)] - \mathbf E[\Phi(\mathcal
U_{\infty}^0)] = \mathcal O(\alpha) \qquad\mbox{as } \alpha\to0.
\end{equation}
\end{lemma}
\begin{pf}
First, note that mutation is parent-independent here, $z=1$. Let
$\Phi(\smallu) = \langle\nu^{\smallu}, \phi\rangle$ with
$\phi\in\widebar{\mathcal C}{}^1_n$. Recall from the proof of
Theorem~\ref{T4} [see (\ref{eqdual2p})] that $\mathbf
E[\Phi(\mathcal U_{\infty}^\alpha)] = \mathbf
E_\phi[\Xi_{\infty}^\alpha]$, where $(\Xi_t^\alpha)_{t\geq0}$ is the
dual process with selection coefficient $\alpha$ and $\Xi_0^\alpha=
\phi$. For the proof, we couple the dual processes for selection
coefficients $\alpha$ and $0$ using the same transitions as given by
(\ref{eqdua0}), (\ref{eqdua1}) and (\ref{eqdua3}). Recall that
there is a random time $T<\infty$ such that $\Xi_t^\alpha=
\Xi_T^\alpha$ for $t\geq T$ and $\Xi^\alpha_{\infty}=
\Xi^\alpha_T$. The only difference between $(\Xi_t^\alpha)_{t\geq
0}$ and $(\Xi_t^0)_{t\geq0}$ is that only the former process can
make transitions given by (\ref{eqdua3a}) or
(\ref{eqdua3b}). Hence, for the coupled process, we get
$\Xi^\alpha_{\infty}= \Xi^0_{\infty}$ if no\vspace*{1pt} such transition occurs
before time $T$. Consider a time $s$ when $\Xi_s^\alpha\in
\widebar{\mathcal C}{}^1_k$. By (\ref{eq200pp}), the chance that a
selective event occurs until time $t$ when $\Xi_t^\alpha\in
\widebar{\mathcal C}{}^1_{k-1}$ is (recall $z=1$) $\alpha k/(
\alpha k + \gamma{ k\choose2} + \vartheta k)$. Hence, for some
finite $C, C'>0$, depending only on $\phi$ and
$\vartheta$,\looseness=-1
%
%e8.3 #&#
%
\begin{eqnarray}
\label{eq841}
| \mathbf E[\Phi(\mathcal U_{\infty}^\alpha)] - \mathbf
E[\Phi(\mathcal U_{\infty}^0)]|
&=& |\mathbf
E_\phi[\Xi_{\infty}^\alpha] - \mathbf E_\phi[\Xi_{\infty}^0]|
\nonumber\\
&\leq&
C \cdot\mathbf P[\Xi_{\infty}^\alpha\neq\Xi_{\infty}^0]\\
&\leq&
C \sum_{k=1}^n \frac{\alpha}{(\alpha+ \vartheta)+ (\gamma/2)(k-1)}
\leq
C'\alpha\nonumber
\end{eqnarray}\looseness=0
for small $\alpha$ and the result follows.
\end{pf}

We start more generally than needed in the proof of Theorem~\ref{T5}. In
particular, given $\underline{\underline r}$ is the distance matrix of an
ultrametric tree, we define tree lengths for subtrees of any finite
number of
leaves.
%
%de8.2 #&#
%
\begin{definition}[(Tree lengths and test functions)]
\label{defPhiij}
(1) For $\underline{\underline r}\in\mathbb R_+^{{\mathbb
N\choose2}}$, we define
%
%e8.4 #&#
%
\begin{equation}\label{eqelln}
\ell_n(\underline{\underline r}) =
\inf_{\sigma\in\widetilde\Sigma_n} \sum_{i=1}^n r_{i\sigma(i)},
\end{equation}
where $\widetilde\Sigma_n \subseteq\Sigma_n$ is the set of
permutations of $\mathbb N$ leaving $n+1,n+2,\ldots$ constant and
having exactly one cycle on $1,\ldots,n$.

(2) For fixed $\lambda\geq0$ let $\phi_{ij}^n\in\overline
{\mathcal
C}{}^1_{n+j}$ be of the form
%
%e8.5 #&#
%
\begin{equation}
\phi_{ij}^n(\underline{\underline r}, \underline u) = e^{ -
\lambda\cdot\ell_n(\underline{\underline r})} \cdot\ind{u_1
= \bullet}\cdots\ind{u_i=\bullet}\cdot\ind{u_{n+1}=\bullet
}\cdots
\ind{u_{n+j}=\bullet}.
\end{equation}
For consistency, we define $\phi_{00}^1:=1$. Moreover we set $\Phi
_{ij}^n:=
\Phi^{n+j,\phi_{ij}^n}$.\vadjust{\goodbreak}
\end{definition}
%
%re8.3 #&#
%
\begin{remark}[(Interpretation)]
(1) If $\underline{\underline r}$ is the distance matrix arising
by sampling points $x_1,x_2,\ldots$ from an ultrametric space
$(U,r,\mu)$, it was shown in Lem\-ma~3.1 of
\citet{GrevenPfaffelhuberWinter2011} that
$\ell_n(\underline{\underline r})$ gives the subtree length of the
subtree spanned by $x_1,\ldots,x_n$.

(2) Considering $\Phi_{ij}^n(\smallu)$ as a function of
$\lambda$
gives the Laplace transform of the subtree length of $n$ sampled
points from $\smallu$ on the set where $i$ points within the
subtree and an additional number $j$ outside the subtree carry
allele $\bullet$. In particular, $\phi_{ij}^n$ depends on the
first $n+j$ points, and hence $\Phi_{ij}^n\in\overline{\mathcal
C}{}^1_{n+j}$.
\end{remark}

%s8.1 #&#
\subsection{Equilibrium distances under neutrality}
The action of $\Omega$ on functions $\Phi_{ij}^n$ given in
Definition~\ref{defPhiij} has a particularly nice form for
$\alpha=0$. Recall that $\Omega_\alpha$ denotes the generator given in
(\ref{dgpGir2}) for $\alpha\geq0$.
%
%le8.4 #&#
%
\begin{lemma}[(Action of $\Omega_0$ on $\Phi_{ij}^n$)]%
\label{l94b}
Let $\alpha=0$ and $\Phi_{ij}^n$ be as in
Definition~\ref{defPhiij}. Then
%
%e8.6 #&#
%
\begin{eqnarray}
\label{eqmain}
\Omega_0 \Phi_{ij}^n & = & -n\lambda
\Phi^{n}_{ij}\mathbh{1}_{n\geq2} + i \frac12
(\vartheta_{{\darkgreybullett}} \Phi_{i-1,j}^{n-1} -
\overline\vartheta\Phi_{ij}^{n}) + j \frac12 (
\vartheta_{{\darkgreybullett}} \Phi_{i,j-1}^{n} -
\overline\vartheta\Phi_{ij}^{n}) \nonumber\\
&&{} + \gamma\biggl(
\pmatrix{ i\cr2} \bigl( \Phi_{i-1,j}^{(n-1)} - \Phi_{ij}^{n}\bigr) +
i(n-i) \bigl( \Phi_{ij}^{(n-1)} - \Phi_{ij}^{n}\bigr)
\nonumber\\[-8pt]\\[-8pt]
&&\hspace*{26.5pt}{} + \pmatrix{n-i\cr2} \bigl( \Phi_{ij}^{(n-1)} -
\Phi_{ij}^{n}\bigr) + ij ( \Phi_{i,j-1}^{n} -
\Phi_{ij}^{n}) \nonumber\\
&&\hspace*{26.5pt}{} +
(n-i)j (
\Phi_{i+1,j-1}^{n} - \Phi_{ij}^{n}) + \pmatrix{ j\cr2} (
\Phi_{i,j-1}^{n} - \Phi_{ij}^{n}) \biggr).\nonumber
\end{eqnarray}
\end{lemma}
\begin{pf}
First, observe that for $n\geq2$
%
%e8.7 #&#
%
\begin{equation}\label{eq701}
\bigl\langle\nabla_{\underline{\underline r}}e^{-\lambda\cdot
\ell_n(\underline{\underline r})}, \underline{\underline
2}\bigr\rangle= -n \lambda\cdot e^{-\lambda\cdot
\ell_n(\underline{\underline r})},
\end{equation}
which explains the first term on the right-hand side of
(\ref{eqmain}). Mutation to ${\darkgreybullet}$ occurs at
rate $\frac{\vartheta_\bullet}{2}$ and to $\bullet$ with rate
$\frac{\vartheta_{\darkgreybullett}}{2}$. Hence, for
$\phi\in\mathcal B(I)$
%
%e8.8 #&#
%
\begin{equation}\label{eq702}\quad
\frac{\overline\vartheta}{2}B \phi(u) = \frac{\vartheta_\bullet}
2\ind{u=\bullet} \bigl(\phi({\darkgreybullet})-\phi({\bullet})\bigr)
+ \frac{\vartheta_{\darkgreybullett}}
2\bigl(1-\ind{u={\bullet}}\bigr)\bigl(\phi({\bullet})-\phi({\darkgreybullet})\bigr).
\end{equation}
In particular,
%
%e8.9 #&#
%
\begin{equation}\label{eq703}
B \ind{u={\bullet}} = -
\frac{\vartheta_{\bullet}}{2}\ind{u={\bullet}} +
\frac{\vartheta_{\darkgreybullett}}{2}\bigl(1-\ind{u={\bullet}}\bigr)
= \frac{\vartheta_{\darkgreybullett}}{2} -
\frac{\overline\vartheta}{2}\ind{u={\bullet}}.
\end{equation}
Since the mutation operator acts on all components in $\phi_{ij}^n$
separately, we obtain the second and third term in
(\ref{eqmain}). Finally, resampling can happen between any of the
${n+j\choose2}$ with different results within and outside the
subtree and the result follows.
\end{pf}
%
%pr8.5 #&#
%
\begin{proposition}[($\Phi_{ij}^n$ under neutrality)]\label{Pneu}
Let $\mathcal U^0_{\infty}$ be distributed as in
Theorem~\ref{T4} with $\alpha=0$ and the mutation given by
(\ref{eq499}). Then
%
%e8.10 #&#
%e8.11 #&#
%e8.12 #&#
%e8.13 #&#
%e8.14 #&#
%e8.15 #&#
%e8.16 #&#
%e8.17 #&#
%
\begin{eqnarray}
\label{eql11}
\mathbf E[\Phi_{00}^{1}(\mathcal U^0_{\infty})] & = & 1,\\[-2pt]
\label{eql12}
\mathbf E[\Phi_{10}^{1}(\mathcal U^0_{\infty})] & = & \mathbf
E[\Phi_{01}^{1}(\mathcal U^0_{\infty})] =
\frac{\vartheta_{{\darkgreybullett}}}{\vartheta_\bullet
+ \vartheta_{{\darkgreybullett}}},\\[-2pt]
\label{eql13}
\mathbf E[\Phi_{02}^{1}(\mathcal U^0_{\infty})] & = & \mathbf
E[\Phi_{11}^{1}(\mathcal U^0_{\infty})] =
\frac{\vartheta_{{\darkgreybullett}} +
\gamma}{\vartheta_{{\darkgreybullett}} +
\vartheta_\bullet+ \gamma}\cdot
\frac{\vartheta_{{\darkgreybullett}}}{\vartheta_\bullet
+ \vartheta_{{\darkgreybullett}}},\\[-2pt]
\label{eql14}
\mathbf E[\Phi_{00}^{2}(\mathcal U^0_{\infty})] & = &
\frac{\gamma}{\gamma+2\lambda},\\[-2pt]
\label{eql15}
\mathbf E[\Phi_{10}^{2}(\mathcal U^0_{\infty})] & = & \mathbf
E[\Phi_{01}^{2}(\mathcal U^0_{\infty})] =
\frac{\gamma}{\gamma+2\lambda} \cdot\frac{
\vartheta_{{\darkgreybullett}}}{\vartheta_\bullet
+ \vartheta_{{\darkgreybullett}}},\\[-2pt]
\label{eql16}
\mathbf E[\Phi_{20}^{2}(\mathcal U^0_{\infty})] & = &
\frac{\vartheta_{{\darkgreybullett}}}{\vartheta_\bullet+
\vartheta_{{\darkgreybullett}}} \cdot
\frac{\gamma}{\gamma+ 2\lambda} \cdot\frac{ \gamma+ 2\lambda+
\vartheta_{{\darkgreybullett}}}{\gamma
+ 2\lambda+ \vartheta_\bullet+
\vartheta_{{\darkgreybullett}}},\\[-2pt]
\label{eql17}
\mathbf E[\Phi_{11}^{2}(\mathcal U^0_{\infty})] & = &
\frac{\vartheta_{{\darkgreybullett}}}{\vartheta_\bullet+
\vartheta_{{\darkgreybullett}}}\cdot
\frac{\gamma}{3\gamma+2\lambda+\vartheta_\bullet+
\vartheta_{{\darkgreybullett}}}\nonumber\\[-2pt]
&&{}\times\biggl( \frac{\gamma+
\vartheta_{{\darkgreybullett}}}{\gamma+ 2\lambda}
 + \frac{\gamma+
\vartheta_{{\darkgreybullett}}}{\gamma+\vartheta_\bullet
+\vartheta_{{\darkgreybullett}}}\\[-2pt]
&&\hspace*{17.5pt}{}+ \frac{\gamma}{\gamma+
2\lambda} \cdot\frac{ \gamma+ 2\lambda+
\vartheta_{{\darkgreybullett}}}{\gamma+ 2\lambda+
\vartheta_\bullet+ \vartheta_{{\darkgreybullett}} } \biggr),
\nonumber\\[-2pt]
\label{eql18}
\mathbf E[\Phi_{02}^{2}(\mathcal U^0_{\infty})] & = &
\frac{\vartheta_{{\darkgreybullett}}}{\vartheta_\bullet+
\vartheta_{{\darkgreybullett}}} \cdot
\frac{\gamma}{6\gamma+
2\lambda+\vartheta_\bullet+\vartheta_{{\darkgreybullett}}}
\nonumber\\[-2pt]
&&{} \times\biggl( \frac{\gamma+
\vartheta_{{\darkgreybullett}}}{\gamma+ 2\lambda}
 + \frac{\gamma+
\vartheta_{{\darkgreybullett}}}{\gamma+\vartheta_\bullet
+ \vartheta_{{\darkgreybullett}}}\nonumber\\[-2pt]
&&\hspace*{18.2pt}{}+ \frac{4\gamma}{3\gamma
+ 2\lambda+ \vartheta_\bullet+
\vartheta_{{\darkgreybullett}}}\\[-2pt]
&&\hspace*{29pt}{}\times\biggl( \frac{\gamma+
\vartheta_{{\darkgreybullett}}}{\gamma+ 2\lambda}
+ \frac{\gamma+
\vartheta_{{\darkgreybullett}}}{\gamma+
\vartheta_\bullet- \vartheta_{{\darkgreybullett}}}\nonumber\\[-2pt]
&&\hspace*{29pt}\hspace*{17.5pt}{} +
\frac{\gamma}{\gamma+ 2\lambda} \cdot\frac{\gamma+ 2\lambda+
\vartheta_{{\darkgreybullett}}}{\gamma+ 2\lambda+
\vartheta_\bullet+
\vartheta_{{\darkgreybullett}}}\biggr)\biggr).\nonumber
\end{eqnarray}
\end{proposition}
\begin{pf}
The proof is based on (\ref{eq707}) for the special choice of
functions as in Definition~\ref{defPhiij}. Clearly, (\ref{eql11})
holds since $\Phi_{00}^1(\mathcal U^0_{\infty})=1$ by definition. The
left and the middle expression in (\ref{eql12}) both give the
probability that a single chosen individual has the
$\bullet$-allele. This is
$\frac{\vartheta_{{\darkgreybullett}}}{\vartheta_\bullet+
\vartheta_{{\darkgreybullett}}}$, as can, for example, be
seen from
competing Poisson processes along the ancestral line of the one
chosen individual (or a generator calculation).

In the rest of the proof, we abbreviate
%
%e8.18 #&#
%
\begin{equation}\label{eq705}
\Phi_{ij}^n:= \mathbf E[\Phi_{ij}^n(\mathcal U^0_{\infty})].\vadjust{\goodbreak}
\end{equation}
We have, using Lemma~\ref{l94b} for $\Phi_{11}^1$ and $\Phi_{02}^1$
%
%e8.19 #&#
%
\begin{eqnarray}
\label{eq706}
0 & = & \tfrac12 (\vartheta_{{\darkgreybullett}}
\Phi_{01}^{1} - \overline\vartheta\Phi_{11}^{1} +
\vartheta_{{\darkgreybullett}} \Phi_{10}^{1} -
\widebar\vartheta\Phi_{11}^{1}) + \gamma( \Phi_{10}^{1}
-
\Phi_{11}^{1}),\nonumber\\[-8pt]\\[-8pt]
0 & = & (\vartheta_{{\darkgreybullett}} \Phi_{01}^{1} -
\overline\vartheta\Phi_{02}^{1}) + \gamma(
\Phi_{01}^{1} - \Phi_{02}^{1}),\nonumber
\end{eqnarray}
which implies (\ref{eql13}). For (\ref{eql14}), the only
nonvanishing resampling term in (\ref{eqmain}) is the one with
rate ${n-i\choose2}$; hence, applying Lemma~\ref{l94b} for
$\Phi_{00}^2$,
%
%e8.20 #&#
%
\begin{equation}\label{eq708}
0 = -2 \lambda\Phi_{00}^{2} + \gamma(1 - \Phi_{00}^{2}),
\end{equation}
and the result follows. [Of course, (\ref{eql14}) can also be shown
by the fact that the MRCA of two sampled individuals in equilibrium
has a coalescent time which is exponential with rate $\gamma$.]

Let us turn to (\ref{eql15}). We find from (\ref{eqmain}),
%
%e8.21 #&#
%
\begin{eqnarray}
\label{eq709}\qquad
0 & = & -2 \lambda\Phi_{10}^{2} + \tfrac12 (
\vartheta_{{\darkgreybullett}} \Phi_{00}^{2} -
\overline\vartheta\Phi_{10}^{2}) +
\gamma(\Phi_{10}^{1} - \Phi_{10}^{2}),\nonumber\\[-8pt]\\[-8pt]
0 & = & -2 \lambda\Phi_{01}^{2} + \tfrac12 (
\vartheta_{{\darkgreybullett}} \Phi_{00}^{2} -
\overline\vartheta\Phi_{01}^{2}) + \gamma( \Phi_{01}^{1} -
\Phi_{01}^{2} + 2\Phi_{10}^{2} - 2\Phi_{01}^{2}).
\nonumber
\end{eqnarray}
From the difference of the last two equations, the first equality in
(\ref{eql15}) follows. Solving the first equations for
$\Phi_{10}^{2}$ and using (\ref{eql12}) and (\ref{eql13}) then
gives the second equality in (\ref{eql15}). [Again, we remark that
(\ref{eql15}) is not surprising: $\Phi_{10}^{2}$ as well as
$\Phi_{01}^{2}$ give the Laplace transform for two randomly chosen
points, given one of the points or a third point has type
$\bullet$. Following back the ancestral line of the latter point
shows that the Laplace transform is independent of the type of the
other chosen individual.]

Next, we have
%
%e8.22 #&#
%
\begin{equation}\label{eq710}
0 = -2\lambda\Phi_{20}^{2} + (
\vartheta_{{\darkgreybullett}} \Phi_{10}^{2} -
\overline\vartheta\Phi_{20}^{2}) + \gamma( \Phi_{10}^{1} -
\Phi_{20}^{2}),
\end{equation}
which shows (\ref{eql16}). For (\ref{eql17}) and (\ref{eql18}),
we have the pair of equations
%
%e8.23 #&#
%
\begin{eqnarray}
\label{eq711}
0 & = & -2\lambda\Phi_{11}^{2} + \tfrac12
(\vartheta_{{\darkgreybullett}} \Phi_{01}^{2} -
\overline\vartheta\Phi_{11}^{2} +
\vartheta_{{\darkgreybullett}} \Phi_{10}^{2} -
\overline\vartheta\Phi_{11}^{2} ) \nonumber\\
&&{}
+ \gamma( \Phi_{11}^{1} -
\Phi_{11}^{2} + \Phi_{10}^{2} - \Phi_{11}^{2} +
\Phi_{20}^{2} - \Phi_{11}^{2}),\nonumber\\[-8pt]\\[-8pt]
0 & = & -2\lambda\Phi_{02}^{2} + (
\vartheta_{{\darkgreybullett}}\Phi_{01}^{2} -
\overline\vartheta\Phi_{02}^{2}) \nonumber\\
&&{}
+ \gamma( \Phi_{02}^{1} - \Phi_{02}^{2} +
\Phi_{01}^{2} - \Phi_{02}^{2} + 4 \Phi_{11}^{2} - 4
\Phi_{02}^{2}).\nonumber
\end{eqnarray}
Solving this linear system (e.g., by using \textsc{Mathematica}) gives
the assertions.
\end{pf}

%s8.2 #&#
\subsection{\texorpdfstring{Proof of Theorem \protect\ref{T5}}{Proof of Theorem 5}}
First, by Lemma~\ref{l94}, $\mathbf E[\Phi(\mathcal U^\alpha_{\infty})]
= \mathbf E[\Phi(\mathcal U^{0}_{\infty})] + \mathcal O(\alpha)$ for
$\alpha\to0$. Hence, by applying (\ref{eq707}) to the function
$\Phi_{00}^2$ from Definition~\ref{defPhiij},
%
%e8.24 #&#
%
\begin{eqnarray}
\label{eq712}
0 &=& -2 \lambda\mathbf E[\Phi_{00}^{2}(\mathcal
U^\alpha_{\infty})] + \gamma\cdot\mathbf E[1 -
\Phi_{00}^{2}(\mathcal U^\alpha_{\infty})]\nonumber\\[-8pt]\\[-8pt]
&&{} + 2\alpha\mathbf E[
\Phi_{10}^{2}(\mathcal U^\alpha_{\infty}) - \Phi_{01}^{2}(\mathcal
U^\alpha_{\infty}) ].\nonumber\vadjust{\goodbreak}
\end{eqnarray}
Since
%
%e8.25 #&#
%
\begin{equation}
\label{eq713}
\mathbf E[ \Phi_{10}^{2}(\mathcal U^\alpha_{\infty}) -
\Phi_{01}^{2}(\mathcal U^\alpha_{\infty}) ] = \mathbf E[
\Phi_{10}^{2}(\mathcal U^0_{\infty}) - \Phi_{01}^{2}(\mathcal
U^0_{\infty}) ] + \mathcal O(\alpha) = \mathcal O(\alpha)\hspace*{-28pt}
\end{equation}
by Lemma~\ref{l94} and Lemma~\ref{l94b}, we find that
%
%e8.26 #&#
%
\begin{equation}
\label{eq714}
\mathbf E[\Phi_{00}^{2}(\mathcal U^\alpha_{\infty})] =
\frac{\gamma}{\gamma+2\lambda} + \mathcal O(\alpha^2).
\end{equation}
Now, in order to compute $\mathbf E[ \Phi_{10}^{2}(\mathcal
U^{\alpha}_{\infty}) - \Phi_{01}^{2}(\mathcal U^{\alpha}_{\infty}) ]$ more
accurately, up to second order in $\alpha$, we apply the equilibrium
condition (\ref{eq707}) on $\Phi_{10}^{2}-\Phi_{01}^{2}$ and obtain,
since $\Phi_{10}^1=\Phi_{01}^1$,
%
%e8.27 #&#
%
\begin{eqnarray}
\label{eq715}
0 & = & -2\lambda\mathbf E[ \Phi_{10}^{2}(\mathcal
U^{\alpha}_{\infty}) - \Phi_{01}^{2}(\mathcal U^{\alpha}_{\infty}) ]
\nonumber\\[-2pt]
&&{} + \frac{\vartheta_{\darkgreybullett}}{2}
\mathbf E[\Phi_{00}^2(\mathcal U^{\alpha}_{\infty})] -
\frac{\overline\vartheta}{2} \mathbf E[\Phi_{10}^2(\mathcal
U^{\alpha}_{\infty})] -
\frac{\vartheta_{\darkgreybullett}}{2} \mathbf
E[\Phi_{00}^2(\mathcal U^{\alpha}_{\infty})] +
\frac{\overline\vartheta}{2} \mathbf E[\Phi_{01}^2(\mathcal
U^{\alpha}_{\infty})] \nonumber\hspace*{-35pt}\\[-2pt]
&&{} + \gamma\bigl( \mathbf
E[\Phi_{10}^1(\mathcal U^{\alpha}_{\infty})-\Phi_{10}^2(\mathcal
U^{\alpha}_{\infty})] \nonumber\\[-2pt]
&&\hspace*{23.4pt}{} - \mathbf
E[\Phi_{01}^1(\mathcal U^{\alpha}_{\infty})-\Phi_{01}^2(\mathcal
U^{\alpha}_{\infty})] - 2\mathbf E[\Phi_{10}^2(\mathcal
U^{\alpha}_{\infty})-\Phi_{01}^2(\mathcal U^{\alpha}_{\infty})]\bigr)
\nonumber\\[-2pt]
&&{} + \alpha\bigl(\mathbf E[\Phi_{10}^2(\mathcal
U^{\alpha}_{\infty}) +\Phi_{20}^2(\mathcal
U^{\alpha}_{\infty})-2\Phi_{11}^2(\mathcal U^{\alpha}_{\infty}) ]
\nonumber\\[-8.5pt]\\[-8.5pt]
&&\hspace*{21.5pt}{}
- \mathbf E[2\Phi_{11}^2(\mathcal
U^{\alpha}_{\infty})+ \Phi_{01}^2(\mathcal
U^{\alpha}_{\infty})-3\Phi_{02}^2(\mathcal
U^{\alpha}_{\infty})]\bigr)\nonumber\\[-2pt]
& = & \biggl(-2\lambda-
\frac{\overline\vartheta}{2} - 3\gamma\biggr)\mathbf E[
\Phi_{10}^{2}(\mathcal U^{\alpha}_{\infty}) - \Phi_{01}^{2}(\mathcal
U^{\alpha}_{\infty}) ] + \gamma\mathbf E[\Phi_{10}^1(\mathcal
U^{\alpha}_{\infty}) - \Phi_{01}^1(\mathcal U^{\alpha}_{\infty}
)]\nonumber\\[-2pt]
&&{}
+ \alpha\mathbf E[\Phi_{10}^{2}(\mathcal
U^{\alpha}_{\infty})-\Phi_{01}^{2}(\mathcal U^{\alpha}_{\infty}) +
\Phi_{20}^{2}(\mathcal U^{\alpha}_{\infty}) -
4\Phi_{11}^{2}(\mathcal U^{\alpha}_{\infty}) +
3\Phi_{02}^{2}(\mathcal U^{\alpha}_{\infty})]\nonumber\\[-2pt]
& = &
\biggl(-2\lambda-
\frac{\overline\vartheta}{2} - 3\gamma\biggr)\mathbf E[
\Phi_{10}^{2}(\mathcal U^{\alpha}_{\infty}) - \Phi_{01}^{2}(\mathcal
U^{\alpha}_{\infty}) ] \nonumber\\[-2pt]
&&{} + \alpha\mathbf
E[\Phi_{20}^{2}(\mathcal U^{0}_{\infty}) - 4\Phi_{11}^{2}(\mathcal
U^{0}_{\infty}) + 3\Phi_{02}^{2}(\mathcal U^{0}_{\infty})]+ \mathcal
O(\alpha^2).\nonumber
\end{eqnarray}
In particular, under neutrality, by Proposition~\ref{Pneu},
%
%e8.28 #&#
%
\begin{eqnarray}
\label{eq716}
&&
\mathbf E[ \Phi_{20}^{2}(\mathcal U^{0}_{\infty})  -
4\Phi_{11}^{2}(\mathcal U^{0}_{\infty}) + 3\Phi_{02}^{2}(\mathcal
U^{0}_{\infty})] \nonumber\\[-8pt]\\[-8pt]
&&\qquad = \frac{2\gamma
\vartheta_\bullet\vartheta_{{\darkgreybullett}}(2\gamma+
2\lambda+ \overline\vartheta)}{\overline\vartheta(\gamma+
\overline\vartheta)(\gamma+ 2\lambda+ \overline\vartheta)(6\gamma+
2\lambda+ \overline\vartheta)}\lambda.\nonumber
\end{eqnarray}
Now, combining (\ref{eq712}), (\ref{eq715}) and (\ref{eq716}), we
see that
\begin{eqnarray*}
% \label{eq717}
&&\mathbf E [ \Phi_{00}^{2} (\mathcal U^\alpha_{\infty})] -
\frac{\gamma}{\gamma+ 2\lambda} \\[-2pt]
&&\qquad= \frac{2\alpha}{\gamma+
2\lambda} \mathbf E[\Phi_{10}^{2} (\mathcal U^\alpha_{\infty}) -
\Phi_{01}^{2} (\mathcal U^\alpha_{\infty})]\\[-2pt]
&&\qquad =
\frac{2\alpha^2}{(\gamma+ 2\lambda)(3\gamma+ 2\lambda+
(1/2)\overline\vartheta)} \mathbf E[\Phi_{20}^{2}(\mathcal
U^{0}_{\infty}) - 4\Phi_{11}^{2}(\mathcal U^{0}_{\infty}) +
3\Phi_{02}^{2}(\mathcal U^{0}_{\infty})]\\[-2pt]
&&\qquad\quad{} + \mathcal O(\alpha^3) \\[-2pt]
&&\qquad
= \frac{8\gamma
\vartheta_\bullet\vartheta_{{\darkgreybullett}}(2\gamma+
2\lambda+ \overline\vartheta)}{\overline\vartheta(\gamma+
\overline\vartheta)(\gamma+ 2\lambda+
\overline\vartheta)(6\gamma+ 2\lambda+
\overline\vartheta)(\gamma+ 2\lambda)^2(6\gamma+ 4\lambda+
\overline\vartheta)}\lambda\alpha^2\\[-2pt]
&&\qquad\quad{} + \mathcal O(\alpha^3),
\end{eqnarray*}
and the assertion follows.

%apA #&#
%
\begin{appendix}\label{app}
\section*{Appendix: Notation}

We collect the most important notation here:
\begin{longlist}[$\blacktriangleright$]
\item[$\blacktriangleright$] $N$: population size of Moran model (Section
\ref{ssbasicdiscrete}),
\item[$\blacktriangleright$] $I$: type space, compact metric space (Section
\ref{ssbasicdiscrete}),
\item[$\blacktriangleright$] $U_N:= \{1,\ldots,N\}$ (Definition \ref
{defgraph}),%\eqref{eqRN}
\item[$\blacktriangleright$] $S_N:=R_N\times[0,\infty)$ (Definition
\ref{defgraph}),
\item[$\blacktriangleright$] $A_s(l,t)\in U_N$: ancestor of
individual $l$ at time $s$
(Definition~\ref{defgraph}),
\item[$\blacktriangleright$] $\eta$: Poisson processes
(Definition~\ref{defgraph}), %\eqref{eq303}
\item[$\blacktriangleright$] $\gamma$: resampling rate (\ref{eqgamma}),
\item[$\blacktriangleright$] $\vartheta$: mutation rate (\ref
{eqvarthetabeta}),
\item[$\blacktriangleright$] $\beta(u,dv)$ transition kernel on $I$
for mutation
(\ref{eqvarthetabeta}),
\item[$\blacktriangleright$] $\widebar\beta, \widetilde\beta$:
two components of $\beta$ for
a parent-independent part (\ref{eq1001}),
\item[$\blacktriangleright$] $\alpha$: selection coefficient (\ref
{eqalphachi}),
\item[$\blacktriangleright$] $\chi(u)$, $\chi'(u,v)$: haploid
fitness of type $u$ and diploid
of $\{u,v\}$ (\ref{eqalphachi}),~(\ref{eqalphawtchi}),
\item[$\blacktriangleright$] $\widehat\chi, \widehat\chi'$:
fitness functions for
measure-valued process (\ref{eq206b}), (\ref{eq208}),
% \item$\preceq_t$: partial order on a tree
% (Definition~\ref{defTMMMS})
%
\item[$\blacktriangleright$] $\mathbb M^I$: set of marked metric
measure spaces
(\ref{eqPP001}),
\item[$\blacktriangleright$] $\mathbb U^I, \mathbb U^I_c$: state
space of the processes
(Definition~\ref{defmmm}),\vspace*{1pt}
\item[$\blacktriangleright$] $\smallx= \overline{(X,r,\mu)},
\smallu=
\overline{(U,r,\mu)}$: generic elements of $\mathbb U^I$
(Definition~\ref{defmmm}),
\item[$\blacktriangleright$] $\widehat\Omega$: generator of the
measure-valued Fleming--Viot
process (\ref{eq202}),
\item[$\blacktriangleright$] $\Omega$: generator of the TFVMS, also
$\Omega_\alpha$
(\ref{eq202}),
\item[$\blacktriangleright$] $\mathcal U = (\mathcal U_t)_{t\geq0}$:
the TFVMS
(Theorem~\ref{T1}),\vspace*{1pt}
\item[$\blacktriangleright$] $\mathcal U^N = (\mathcal U^N_t)_{t\geq
0}$: the TMMMS
(Definition~\ref{defTMMMS2})
\item[$\blacktriangleright$] $\mathcal U_{\infty}$: long-time limit of
$\mathcal U$
(Theorem~\ref{T4}),
\item[$\blacktriangleright$] $\zeta^N$: measure-valued Moran model
(\ref{eq495}),
\item[$\blacktriangleright$] $\zeta$: measure-valued Fleming--Viot
process (Example
\ref{exmFV}),
\item[$\blacktriangleright$] $\varphi\dvtx E\to E'$: embedding (Remark
\ref{remnotation}),
\item[$\blacktriangleright$] $\nu^\smallx$: distance matrix
distribution (\ref{eq23}),
\item[$\blacktriangleright$] $\Sigma$: set of permutations
(\ref{eqSigma}),
\item[$\blacktriangleright$] $\theta$: resampling operator
(\ref{gen2}),
\item[$\blacktriangleright$] $R_\sigma$: map exchanging indices
according to permutation
$\sigma$ (\ref{eq29}),
\item[$\blacktriangleright$] $\Phi= \Phi^{n,\phi}$: polynomial
(\ref{eqphi}),
\item[$\blacktriangleright$] $\Pi, \Pi^1$: set of polynomials
(\ref{eq201}),
\item[$\blacktriangleright$] $\sigma_k$, $\overline\sigma_k$:
shift operators
(\ref{eqsig1}), (\ref{eqsig2}),
\item[$\blacktriangleright$] $\rho_1^n$: shift operator (\ref{eq222}),
\item[$\blacktriangleright$] $\Pi_N$: polynomials for finite
populations (\ref{eq7}),
\item[$\blacktriangleright$] $R_{12}$: distance of two randomly
chosen points (Remark~\ref{remR12}),
\item[$\blacktriangleright$] $\Upsilon$: state space of
function-valued dual process
(\ref{eq904}),
\item[$\blacktriangleright$] $\Xi$: dual process (Definition
\ref{defdual}),
\item[$\blacktriangleright$] $\ell_n$: tree length for $n$
individuals (\ref{eqelln}).
\end{longlist}
\end{appendix}

\section*{Acknowledgments}

We thank Anton Wakolbinger for fruitful discussion and Steve Evans for
pointing us to the paper of \citet{BakryEmery1985}.
Part of this work has been
carried out when A. Depperschmidt was taking part in the Junior
Trimester Program Stochastics at the Hausdorff Center in Bonn:
hospitality and financial support are gratefully acknowledged.

%suskaldyti doi

% imsref loaded by lrinkeviciute, 2012-04-30 11:06:25
% imsref loaded by lrinkeviciute, 2012-04-30 12:35:26
%

\printaddresses

\end{document}